\documentstyle[12pt]{article} 
\makeatletter 
\def\mineappendix{ 
        \setcounter{section}{1} 
        \setcounter{subsection}{0} 
        \def\thesection{\Alph{section}} 
        \def\sectionap{\@startsection  {section}{1}{\z@} 
                        {-3.5ex plus-1ex minus-.2ex} {0ex plus.2ex} 
                        {\reset@font\Large\bf  Appendix:  \, } 
                        } 
        } 
\makeatother 
\def\Proclaim #1. #2\par{\bigbreak\noindent{\sc#1.\enspace}{\it#2}\par} 
 
\font\Bbbfont=msbm10 
\newfam\msbfam 
\textfont\msbfam=\Bbbfont 
\def\Bbb#1{{\fam\msbfam\relax#1}}

\newcommand{\eqref}[1]{equation~(\ref{#1})}

\newcommand{\gwii}[1]{\left< \hspace{-2pt} \left< \, #1 \, 
        \right>  \hspace{-2pt} \right>_{0}} 
 
\newcommand{\gwiione}[1]{\left< \hspace{-2pt} \left< \, #1 \, 
        \right> \hspace{-2pt} \right>_{1}} 
 
\newcommand{\gwiitwo}[1]{\left< \hspace{-2pt} \left< \, #1 \, 
        \right> \hspace{-2pt} \right>_{2}} 
 
\newcommand{\gwig}[1]{\left< \, #1 \, \right>_{g}} 
\newcommand{\gwiig}[1]{\left< \hspace{-2pt} \left< \, #1 \, 
    \right> \hspace{-2pt} \right>_{g}} 
\newcommand{\gwiih}[2]{\left< \hspace{-2pt} \left< \, #2 \, \right> \hspace{-2pt} \right>_{#1}} 
 
\newcommand{\grav}[2]{\tau_{#1}(\gamma_{#2})} 
\newcommand{\grava}[1]{\tau_{#1}(\gamma_{\alpha})} 
\newcommand{\gravua}[1]{\tau_{#1}(\gamma^{\alpha})} 
\newcommand{\gravb}[1]{\tau_{#1}(\gamma_{\beta})} 
 
\newcommand{\gravm}[1]{\tau_{#1}(\gamma_{\mu})} 
 
\newcommand{\gravn}[1]{\tau_{#1}(\gamma_{\nu})}

\newcommand{\ga}{\gamma_{\alpha}} 
\newcommand{\gua}{\gamma^{\alpha}} 
\newcommand{\gb}{\gamma_{\beta}} 
\newcommand{\gub}{\gamma^{\beta}} 
\newcommand{\gm}{\gamma_{\mu}} 
\newcommand{\gum}{\gamma^{\mu}} 
\newcommand{\gn}{\gamma_{\nu}} 
\newcommand{\gun}{\gamma^{\nu}}

\newcommand{\ba}{b_{\alpha}} 
\newcommand{\bb}{b_{\beta}}

\newcommand{\vs}{{\cal S}} 
\newcommand{\vx}{{\cal X}} 
\newcommand{\vd}{{\cal D}} 
\newcommand{\vw}{{\cal W}} 
\newcommand{\vv}{{\cal V}} 
\newcommand{\vl}{{\cal L}} 
\newcommand{\vg}{{\cal G}} 
\newcommand{\vz}{{\cal Z}} 
\newcommand{\vy}{{\cal Y}} 
 
\newcommand{\bvs}{{\, \overline{\cal S}\,}} 
\newcommand{\bvx}{{\, \overline{\cal X}\,}} 
\newcommand{\bvd}{{\, \overline{\cal D} \,}} 
\newcommand{\bvw}{{\, \overline{\cal W} \,}} 
 
\newcommand{\bvl}{{\, \overline{\cal L} \,}}

\newcommand{\qp}[2]{#1 \bullet #2} 
\newcommand{\qpc}[2]{{\cal #1} \bullet #2} 
\newcommand{\qpcc}[2]{{\cal #1} \bullet {\cal #2}}

\newtheorem{lem}{Lemma}[section] 
\newtheorem{cor}[lem]{Corollary} 
\newtheorem{thm}[lem]{Theorem}

\newtheorem{defi}[lem]{Definition}

\title{Quantum product on the big phase space and the Virasoro conjecture} 
\author{Xiaobo Liu \thanks{Research partially supported by 
            Alfred P. Sloan Research Fellowship and NSF Postdoctoral 
        Research Fellowship}} 
\date{} 
 
\begin{document} 
\maketitle

Quantum cohomology is a family of new ring structures on the space 
of cohomology classes of a compact symplectic manifold (or a 
smooth projective variety) $V$. The quantum products are defined 
by third order partial derivatives of the generating function of 
primary Gromov-Witten invariants of $V$ (cf. \cite{RT1}). In a 
similar way, using the generating function of all descendant 
Gromov-Witten invariants, we can define products on an infinite 
dimensional vector space, called the  big phase space, which can 
be thought of as a product of infinite copies of the small phase 
space $H^{*}(V; {\Bbb C})$. It seems that the products on the big 
phase space have not gotten enough attention in the literature so 
far. In this paper we will study some basic structures of such 
products and apply them to the study of topological recursion 
relations and the Virasoro conjecture. 
 
The Virasoro conjecture predicts that the generating function 
of the Gromov-Witten invariants is annihilated by infinitely 
many differential operators, denoted by $\{ L_{n} \mid n \geq -1\}$, 
 which form 
a half branch of the Virasoro algebra. 
This conjecture was proposed by Eguchi, Hori and Xiong \cite{EHX} 
and also by S. Katz (cf. \cite{CK} \cite{EJX}). 
It is a natural generalization of a conjecture of Witten (cf. \cite{W2} 
\cite{Kon} \cite{W2}) and provides a powerful tool in the computation 
of Gromov-Witten invariants. 
The genus-0 Virasoro conjecture was proved in \cite{LT} 
(cf. \cite{DZ2} and \cite{G2} for alternative proofs). 
The genus-1 Virasoro conjecture for manifolds with semisimple 
quantum cohomology was proved in \cite{DZ2}. Without assuming semisimplicity, 
the genus-1 Virasoro conjecture was reduced to the genus-1 
$L_{1}$-constraint on the small phase space in \cite{L1}. 
It was also proved in \cite{L1} \cite{L2} that the genus-1 Virasoro conjecture 
holds if the quantum cohomology is not too degenerate (a condition 
weaker than semisimplicity). 
 
The genus-$g$ Virasoro conjecture with $g \geq 1$ can be 
formulated in a way which computes the derivatives of the 
genus-$g$ generating function along a sequence of vector fields, 
called the Virasoro vector fields (see 
Section~\ref{sec:VirConj}). The study of the properties of these 
vector fields will be important in both proving and applying the 
Virasoro conjecture in all genera. In this paper we will give a 
simple recursive description of the Virasoro vector fields (see 
\eqref{eqn:VirVF} and Theorem~\ref{thm:LkRec}). This recursive 
description enables us to understand the relations between the 
Virasoro vector fields and 
 the quantum powers of the Euler vector field defined by 
\eqref{eqn:EulerDef}. The 
action of the Virasoro vector fields on the generating function of 
genus-$g$ Gromov-Witten invariants is equivalent to the action of 
a sequence of vector fields constructed from quantum powers of the 
Euler vector field. To prove this  fact, we need to use the 
quantum product on the big phase space to reinterpret genus-$g$ 
topological recursion relations. The most important difference 
between the quantum product on the big phase space and the one on 
the small phase space is that there is no identity element for 
the product on the big phase space. 
The best candidate for the identity is the  string 
vector field defined by \eqref{eqn:StringDef}. 
However, this vector field is not an identity 
in the usual sense. How close this vector field is to an identity 
is reflected through various topological recursion relations. Such 
interpretation of topological recursion relations will be very 
useful in the study of Virasoro conjecture. For example, this will 
enable us to represent all the vector fields obtained from the 
Virasoro vector fields under certain naturally defined operations 
in terms of twisted quantum powers of the Euler vector field (cf. 
Theorem~\ref{thm:Stau-Lk} and the comments afterwards). This 
explains why Virasoro constraints are very powerful for manifolds 
with semisimple quantum cohomology, as in this case the quantum 
powers of the Euler vector field span the space of primary vector 
fields. 
 
We believe that the structures defined in 
this paper will be very useful in the study of the Virasoro 
conjecture for all genera. As a demonstration, we will apply these 
structures 
 to the genus-2 Virasoro conjecture. 
The study of the genus-2 Virasoro conjecture is important because 
this is the first case that we do not have a formula to reduce the 
problem to the small phase space. The behavior of the Virasoro 
conjecture in this case will provide much needed insight in what 
we should expect in the higher genus cases.  Moreover the 
techniques developed here could be easily adapted to the study of 
higher genus Virasoro conjecture. In this paper, we will prove 
that for any manifold, the genus-2 Virasoro conjecture holds if 
and only if the genus-1 and genus-2 $L_{1}$-constraints hold(see 
Theorem~\ref{thm:L1->L2}). 
The main reason for this result is that for $n \geq 2$, 
the genus-2 part of the 
genus-2 $L_{n}$-constraint can be recursively computed from 
the genus-2 part of the genus-2 
$L_{1}$-constraint (See Corollary~\ref{cor:Recpsi} 
for the precise recursion formula). In the case that the quantum 
cohomology of the underlying manifold is not too degenerate (in particular 
is semisimple), such recursion formula also uniquely determines 
the genus-2 part of $L_{1}$-constraint (see Theorem~\ref{thm:VirNondeg}). 
Therefore the 
genus-2 Virasoro conjecture for such manifolds can now be reduced to a genus-1 
problem. To complete the proof of the genus-2 Virasoro conjecture 
for this case would involve more detailed analysis of the complicated 
tensors $A_{1}$, $A_{2}$ and $B$ used in 
equations (\ref{eqn:TRR1}) - (\ref{eqn:BP}). We will do 
this in a separate paper. 
We notice that the first 3 Virasoro 
operators $\{ L_{-1}, L_{0}, L_{1}\}$ 
 form a 3-dimensional subalgebra of the Virasoro 
algebra which is isomorphic to $sl(2)$. 
So Theorem~\ref{thm:L1->L2} in particular 
 implies that for any manifold, the $sl(2)$ symmetry of the 
genus-2 generating function is sufficient to deduce the genus-2 
Virasoro conjecture. Similar situation also occurred in the 
genus-1 case (cf. \cite{L1}). We wonder whether the same pattern 
will continue for all genera. 
 
The major application of the genus-$g$ Virasoro conjecture is to 
compute the genus-$g$ generating function of the Gromov-Witten 
invariants in terms of data with genus less than $g$. Although in 
principle this can be done if the quantum cohomology is 
semisimple, it is not easy to solve the genus-$g$ generating 
function explicitly  from the Virasoro constraints for general 
manifolds with semisimple quantum cohomology.  Recently 
Dubrovin-Zhang and Eguchi-Getzler-Xiong computed the genus-2 
generating function for Frobenius manifolds with two primary 
fields  assuming that the genus-2 Virasoro constraints hold (cf. 
\cite{EGX}). Note that in the Gromov-Witten theory, only ${\Bbb C}P^{1}$ has 
two primary fields. In this paper,  we will prove (without 
assuming the Virasoro constraints) that the genus-2 generating 
function can be expressed explicitly in terms of genus-0 and 
genus-1 data if the quantum cohomology of the underlying manifold 
is not too degenerate (in particular, if the quantum cohomology is 
semisimple) (see Theorem~\ref{thm:VirNondeg}).  As proved in \cite{DZ1}, 
the genus-1 generating function can be expressed explicitly in terms 
of genus-0 data when the quantum cohomology is semisimple. Therefore 
Theorem~\ref{thm:VirNondeg} also implies that the genus-2 generating function 
can be expressed explicitly in terms of genus-0 data in the semisimple case.

It will be interesting to see how many ingredients are needed in 
the study of the genus-2 Virasoro conjecture. For this purpose, 
we need to study relations among 
 the genus-2 topological 
 recursion relations in \cite{G1} and \cite{BP} 
(see equations (\ref{eqn:TRR1-old}) - (\ref{eqn:BP-old})) . 
In this paper, 
we will show that \eqref{eqn:TRR2-old} 
implies \eqref{eqn:TRR1-old} 
(see Theorems~\ref{thm:TRR2to1}), 
and  \eqref{eqn:TRR1-old} together with \eqref{eqn:BP-old} 
also implies \eqref{eqn:TRR2-old} (see 
Theorem~\ref{thm:BPTRR1toTRR2}). This tells us that, at least for 
manifolds whose quantum cohomology is not too degenerate, 
 the number of ingredients needed in the study of the genus-2 Virasoro 
conjecture is the same as that for the genus-1 case. It will be 
interesting to investigate whether this phenomenon will continue 
in higher genus cases. 
 
This paper is organized as follows. In Section~\ref{sec:QuantProd}, 
we define the quantum product on the big phase space. In 
Section~\ref{sec:TRR}, we first re-formulate topological recursion relations 
using an operator $T$ which measures the difference between the string 
vector field and an identity of the quantum product, 
and then apply it to study 
relations among genus-2 topological recursion relations. 
In Section~\ref{sec:EulerVF}, we study properties of quantum powers of 
the Euler vector field. Most properties of these vector fields 
can be derived from the quasi-homogeneity equation 
(\ref{eqn:quasiallgenus}). 
In particular, the recursive operator $R$ used to describe the 
Virasoro vector fields arrive naturally in the study of the quasi-homogeneity 
equation (see Theorem~\ref{thm:Xprod}). 
Section~\ref{sec:VirVF} is devoted to the study of 
properties of the Virasoro vector fields, in particular their relation 
with quantum powers of the Euler vector field. Using the operator $T$ and its 
right inverse, we can produce some Lie algebras which contains the Lie algebra of 
Virasoro vector fields as a proper subalgebra. In particular, we will give a realization 
of the Lie algebra of integral pseudo-differential operators on the unit circle. 
In Section~\ref{sec:AppVir}, we first formulate the Virasoro conjecture 
using recursive operators, and then study the genus-2 Virasoro conjecture. 
For manifolds whose quantum cohomology is not too degenerate, 
we explicitly solve the genus-2 generating function in terms of 
 genus-0 and genus-1 data in Section~\ref{sec:F2}. In the proof of 
Theorem~\ref{thm:VirNondeg}, we need a lemma which is proven in the appendix.

Part of this paper was written when the author visited MIT in the 
spring 2001. The author would like to thank MIT for hospitality 
and G. Tian for helpful discussions.

\section{Quantum Product on the big phase space} 
\label{sec:QuantProd} 
 
For simplicity, we assume that $V$ is a smooth projective variety 
with $H^{\rm odd}(V; {\Bbb C})= 0$.  All results in this paper 
are also true for compact symplectic manifolds except those 
in Section~\ref{sec:L1->L2}. Choose a 
basis $\{ \gamma_{1}, \gamma_{2}, \ldots, \gamma_{N} \}$ of 
$H^{*}(V; {\Bbb C})$ with $\gamma_{1}$ equal to the identity of 
the ordinary cohomology ring. Let 
\[ \gwig{\grav{n_{1}}{\alpha_{1}} \, \grav{n_{2}}{\alpha_{2}} \, 
    \ldots \, \grav{n_{k}}{\alpha_{k}}} \] 
be the genus-$g$ descendant Gromov-Witten invariant associated 
to $\gamma_{\alpha_{1}}, \ldots, \gamma_{\alpha_{k}}$ and nonnegative 
integers $n_{1}, \ldots, n_{k}$ (which represent the powers 
of the first Chern classes of certain tautological line bundles 
over the moduli space of stable maps from genus-$g$ curves 
to $V$ with $k$ marked points) (cf. \cite{W1} \cite{RT2}). 
The genus-$g$ generating function is defined to be 
\[ F_{g} =  \sum_{k \geq 0} \frac{1}{k!} 
         \sum_{ \begin{array}{c} 
        {\scriptstyle \alpha_{1}, \ldots, \alpha_{k}} \\ 
                {\scriptstyle  n_{1}, \ldots, n_{k}} 
                \end{array}} 
                t^{\alpha_{1}}_{n_{1}} \cdots t^{\alpha_{k}}_{n_{k}} 
    \gwig{\grav{n_{1}}{\alpha_{1}} \, \grav{n_{2}}{\alpha_{2}} \, 
        \ldots \, \grav{n_{k}}{\alpha_{k}}}, \] 
where $\{t^{\alpha}_{n} \mid n \in {\Bbb Z}_{+}, \alpha = 1, 
\cdots, N\}$ is an infinite set of parameters. We can think of 
these parameters as coordinates on an infinite dimensional vector 
space, called the {\it big phase space}. The finite dimensional 
subspace defined by $\{t^{\alpha}_{n} = 0 \, \, \, {\rm if} \, \, 
\, 
        n > 0\}$ 
is called the {\it small phase space}. The function $F_{g}$ is understood 
as a formal power series of $t_{n}^{\alpha}$. 
As in \cite{LT}, it is convenient to 
introduce a $k$-tensor 
 $\left< \left< \right. \right. \underbrace{\cdot \cdots \cdot}_{k} 
        \left. \left. \right> \right> $ 
defined by 
\[ \gwiig{{\cal W}_{1} {\cal W}_{2} \cdots {\cal W}_{k}} \, \, 
         := \sum_{m_{1}, \alpha_{1}, \ldots, m_{k}, \alpha_{k}} 
                f^{1}_{m_{1}, \alpha_{1}} \cdots f^{k}_{m_{k}, \alpha_{k}} 
        \, \, \, \frac{\partial^{k}}{\partial t^{\alpha_{1}}_{m_{1}} 
            \partial t^{\alpha_{2}}_{m_{k}} \cdots 
            \partial t^{\alpha_{k}}_{m_{k}}} F_{g}, 
 \] 
for (formal) vector fields 
${\cal W}_{i} = \sum_{m, \alpha} 
        f^{i}_{m, \alpha} \, \frac{\partial}{\partial t_{m}^{\alpha}}$ where 
$f^{i}_{m, \alpha}$ are (formal) functions on the big phase space. 
We can also view this tensor as the $k$-th 
covariant derivative of $F_{g}$. 
This tensor is called the {\it $k$-point (correlation) function}. 
We will always identify 
$\grav{n}{\alpha}$ with the tangent vector field 
$\frac{\partial}{\partial t^{\alpha}_{n}}$ on the big phase space 
and abbreviate $\grav{0}{\alpha}$ as $\gamma_{\alpha}$. 
We also consider $\grav{n}{\alpha}$ with $n<0$ as a zero vector field. 
Let $\eta_{\alpha \beta} = \int_{V} \gamma_{\alpha} \cup \gamma_{\beta}$ 
be the Poincare pairing. We use $\eta = (\eta_{\alpha \beta})$ 
and its inverse $\eta^{-1} = (\eta^{\alpha \beta})$ to lower and 
raise indices. So 
$\tau_{n}(\gamma^{\alpha}) = \eta^{\alpha \beta} \grav{n}{\beta}$. 
Here we adopt the convention of summing over repeated indices. 
We call a vector field 
${\cal W} = \sum_{m, \alpha} f_{m, \alpha} \grava{m}$ a 
{\it primary vector field} if $f_{m, \alpha} = 0$ whenever $m >0$, 
a {\it descendant vector field} if $f_{m, \alpha} = 0$ whenever $m=0$. 
 
For any two vector fields ${\cal U}$ and ${\cal W}$ on the big phase, 
define the {\bf quantum product} of ${\cal U}$ and ${\cal W}$ by 
\[ \qpcc{U}{W} := 
    \gwii{{\cal U} {\cal W} \, \gamma^{\alpha}} \gamma_{\alpha}. \] 
By definition, the quantum product of two vector fields is always 
a primary vector field. This product is apparently commutative. It 
is also associative due to the {\it generalized WDVV equation} 
\[ \gwii{{\cal W}_{1} {\cal W}_{2} \gamma^{\alpha}} 
    \gwii{\gamma_{\alpha} {\cal W}_{3} {\cal W}_{4}} 
   = \gwii{{\cal W}_{1} {\cal W}_{3} \gamma^{\alpha}} 
    \gwii{\gamma_{\alpha} {\cal W}_{2} {\cal W}_{4}}, 
\] 
which follows in turn from the genus-0 topological recursion 
relation (cf. \cite{W1}). When restricted to tangent vector fields 
on the small phase, this is precisely the product in the quantum 
cohomology of $V$ (also called the big quantum cohomology by some 
authors). For any vector field ${\cal W}$ on the big phase space, 
we define ${\cal W}^{k}$ to be the $k$-th quantum power of ${\cal W}$. 
i.e., 
\[ {\cal W}^{k} = \underbrace{{\cal W} \bullet {\cal W} \bullet \cdots 
        {\cal W}}_{k}, \] 
for $k >0$.

For the quantum product on the small phase space, 
the constant vector field $\gamma_{1}$, which was chosen to be 
the identity of the ordinary cohomology ring, is also the identity 
for the quantum cohomology. However on the big phase space, 
there is no identity vector field for the quantum product. 
A vector field which is close to an identity is the {\it string vector field} 
\begin{equation} \label{eqn:StringDef} 
 {\cal S} = - \sum_{m, \alpha} \tilde{t}^{\alpha}_{m} 
        \grav{m-1}{\alpha}, 
\end{equation} 
where $\tilde{t}^{\alpha}_{m} = t^{\alpha}_{m} - 
    \delta_{m}^{1} \delta_{\alpha}^{1}$. 
This vector field 
 can be considered 
as a sort of identity in the following sense: 
\begin{lem} \label{lem:stringid} 
{\rm (i)}  \hspace{5pt} 
$\qpcc{S}{W} = {\cal W}$ {\rm for any primary vector field} $\cal W$. 
 
\hspace{48pt} {\rm (ii)} \hspace{3pt} 
$\qpcc{\qpcc{S}{U}}{W} = \qpcc{U}{W}$ 
{\rm for all vector fields} $\cal U$ {\rm and} $\cal W$. 
\end{lem} 
{\bf Proof}: This is a consequence of 
the {\it string equation} 
\[ \gwiig{\vs} = \frac{1}{2} \delta_{g, 0} \eta_{\alpha \beta} 
t_{0}^{\alpha} t_{0}^{\beta}. \] 
In fact, taking second derivatives of the genus-0 string equation, we obtain 
$\gwii{{\cal S} \ga \gub} = \delta_{\alpha}^{\beta}$ for all 
$\alpha$ and $\beta$ (cf. \cite[Lemma 1.1]{LT}). This is equivalent to (i). 
The second formula follows from the associativity of the quantum product 
and (i) since $\qpcc{U}{W}$ is a primary vector field by definition. 
$\Box$ 
 
We define an equivalence relation between two vector fields on the 
big phase space: 
\begin{defi} {\rm We say two vectors $\cal U$ and $\cal W$ are 
equivalent, denoted by ${\cal U} \sim {\cal W}$, 
if $\qpcc{U}{V} = \qpcc{W}{V}$ for 
all vector fields $\cal V$.} 
\end{defi} 
The first formula in Lemma~\ref{lem:stringid} implies that 
two primary vector fields are equivalent if and only if they 
are equal. The second formula in Lemma~\ref{lem:stringid} is equivalent to 
\[ \qpcc{S}{W} \sim {\cal W} \] 
for any ${\cal W}$. 
Therefore we can say that $S$ is 
an identity in this weak sense. Since $\qpcc{S}{W}$ is always a primary 
vector field, this relation also tells us that any vector field on the 
big phase space is equivalent to a primary vector field. 
Moreover if ${\cal U}$ is a primary vector field, then 
for any vector field ${\cal W}$, 
\begin{equation} \label{eqn:equivPrimary} 
{\cal W} \sim {\cal U} \, \, \, \Longleftrightarrow 
    \, \, \, \qpcc{S}{W} = {\cal U}. 
\end{equation} 
For the convenience of later computations, we define 
\begin{equation} \label{eqn:primarycorr} 
\bvw := \qpcc{S}{W} 
\end{equation} 
for any vector field $\vw$. In particular, Lemma~\ref{lem:stringid} 
implies that $\bvs$ is the identity of the subalgebra of all primary 
vector fields with respect to the quantum product. Note that 
when restricted to the small phase space, both $\vs$ and $\bvs$ 
are equal to the identity of the ordinary (and the quantum) cohomology 
ring.

It is also convenient to introduce the following (functionally) 
linear transformations on the space of vector fields on the big phase 
space: 
\begin{defi} 
{\rm For any vector field} 
${\cal W} = \sum_{n, \alpha} f_{n, \alpha} \grava{n}$, 
{\rm define} 
\[ \tau_{+}({\cal W}) := \sum_{n, \alpha} f_{n, \alpha} \grava{n +1}, \hspace{50pt} 
   \tau_{-}({\cal W}) := \sum_{n, \alpha} f_{n, \alpha} \grava{n -1},  \] 
\[ T({\cal W}) := \tau_{+}({\cal W}) - \gwii{{\cal W} \gua} \ga,  \hspace{50pt} 
   \pi({\cal W}) := \sum_{\alpha} f_{0, \alpha} \ga. \] 
\end{defi} 
Then 
\[ \tau_{+} \tau_{-} ({\cal W}) = {\cal W} - \pi({\cal W}), \hspace{20pt} 
   \tau_{-} \tau_{+} ({\cal W}) = \tau_{-} T ({\cal W}) = {\cal W}. \] 
Moreover the genus-0 string equation implies that (cf. \cite[Lemma 1.1]{LT}) 
\begin{equation} \label{eqn:StringProd} 
 \bvw = \gwii{\tau_{-}({\cal W}) \gua} \ga + \pi({\cal W}). 
\end{equation} 
An immediate consequence of this formula is the following: 
\begin{equation} \label{eqn:T} 
 T({\cal W}) = \tau_{+}({\cal W}) - \qpc{S}{\tau_{+}({\cal W})} 
\end{equation} 
for any vector field $\vw$. Based on this formula and the 
fact that $\vw - \vs \bullet \vw =0$ if $\vw$ is a primary vector field, 
we can interpret the operator $T$ as a measurement of the difference between 
$\vs$ and a true identity of the quantum product on the big phase space. 
We have the following characterization for the vectors which are equivalent 
to 0. 
\begin{lem} \label{lem:0vector} 
For any vector field $\cal W$, 
\begin{eqnarray*} 
 {\cal W} \sim 0 & \Longleftrightarrow & {\cal W} = T(\tau_{-}({\cal W})) \\ 
    & \Longleftrightarrow & {\cal W} = T({\cal U}) 
    \, \, \, {\rm for \, \, \, some} \, \, \,  {\cal U}. 
\end{eqnarray*} 
\end{lem} 
{\bf Remark}: 
First derivatives of the string equation gives us 
$\gwii{ \vs \gua} = t^{\alpha}_{0}$ (cf. \cite[Lemma 1.1 (2)]{LT}). Therefore 
\[ T({\cal S}) = {\cal D} \] 
where $\cal D$ is the {\it dilaton vector field} 
defined by 
$ {\cal D}= - \sum_{m, \alpha} \tilde{t}^{\alpha}_{m} \, \, \grava{m}.$ 
This lemma in particular implies that 
${\cal D} \sim 0, $ 
which also follows from the {\it dilaton equation} 
\[ \left<\left< {\cal D} \right>\right>_{g} = 
        (2g-2) F_{g} + \frac{1}{24} \, \chi(V)\delta_{g, 1} \] 
(cf. \cite[Lemma 1.2]{LT}). 
Moreover $T({\cal D}) \sim 0$ is equivalent to \cite[Lemma 5.2 (2)]{LT}, 
which is also equivalent to the genus-0 $\widetilde{\cal L}_{1}$ constraint. 
Similar reasoning also 
applies to the genus-0 $\widetilde{\cal L}_{2}$ constraint by considering the 
vector field $T(R({\cal D}))$ where $R$ is defined in Definition~\ref{def:R}. 
 
{\bf Proof of Lemma~\ref{lem:0vector}}: 
By equation~(\ref{eqn:T}), for any vector field $\cal W$, 
\begin{eqnarray*} 
 {\cal W} & = &\tau_{+}(\tau_{-}({\cal W})) + \pi({\cal W}) \\ 
    & = & T(\tau_{-}({\cal W})) + \qpc{S}{\tau_{+}(\tau_{-}({\cal W}))} 
      + \pi({\cal W}) \\ 
    & = & T(\tau_{-}({\cal W})) + \qpc{S}{({\cal W} -\pi({\cal W}))} 
     + \pi({\cal W}). 
\end{eqnarray*} 
Since $\qpc{S}{\pi({\cal W})} = \pi({\cal W})$, we have 
\begin{equation} \label{eqn:WT} 
{\cal W} = T(\tau_{-}({\cal W})) + \bvw. 
\end{equation} 
In particular, ${\cal W} = T(\tau_{-}({\cal W}))$ if ${\cal W} \sim 0$. 
On the other hand, equation~(\ref{eqn:T}) implies $T({\cal U}) \sim 0$ 
for any $\cal U$ since $\bvw \sim {\cal W}$ for all $\cal W$. 
$\Box$ 
 
Let $\nabla$ be the covariant derivative defined by 
\[ \nabla_{\cal V} {\cal W} = 
    \sum_{m, \alpha} \, ({\cal V} f_{m, \alpha}) \, \grava{m} \] 
for any vector fields ${\cal V}$ and 
${\cal W} = \sum_{m, \alpha} \, f_{m, \alpha} \, \grava{m}$. 
Then we have 
\begin{equation} \label{eqn:DerCorr} 
 {\cal V} \gwiig{ {\cal W}_{1} \, {\cal W}_{2} \, \ldots \, {\cal W}_{k}} 
= \gwiig{{\cal V} \, {\cal W}_{1} \, \ldots \, {\cal W}_{k}} 
  + \sum_{j =1}^{k} 
    \gwiig{ {\cal W}_{1} \, \ldots \, (\nabla_{\cal V}{\cal W}_{j}) \, 
        \ldots \, {\cal W}_{k}} 
\end{equation} 
for any vector fields ${\cal V}$ and ${\cal W}_{j}$. 
In particular, we have 
\begin{equation} \label{eqn:DerProd} 
 \nabla_{\cal V} (\qpcc{W}{U}) 
= \gwii{{\cal V} \, {\cal W} \, {\cal U} \, \gua} \ga 
    + \qp{(\nabla_{\cal V} {\cal W})}{\cal U} 
    + \qpc{W}{(\nabla_{\cal V}{\cal U})} 
\end{equation} 
for any vector fields $\cal U$, $\cal V$ and $\cal W$. The 
following formulas will be useful in the study of topological 
recursion relations: 
\begin{lem} \label{lem:DerT} 
For all vector fields ${\cal V}$ and ${\cal W}$, 
\begin{eqnarray*} 
& (1) & \nabla_{\cal V} \, \tau_{+}({\cal  W}) 
    = \tau_{+}( \nabla_{\cal V} {\cal W}). \\ 
& (2) & \nabla_{\cal V} \, \tau_{-}({\cal  W}) 
    = \tau_{-}( \nabla_{\cal V} {\cal W}). \\ 
& (3) & \nabla_{\cal V} \, T({\cal  W}) 
= T(\nabla_{\cal V} {\cal  W})  - \qpcc{V}{W}. 
\end{eqnarray*} 
\end{lem} 
{\bf Proof}: 
It is straightforward to check (1) and (2). 
Hence 
\begin{eqnarray*} 
 \nabla_{\cal V} \, T({\cal  W}) 
&=& \nabla_{\cal V} \, \tau_{+}({\cal  W}) 
    - ({\cal V} \gwii{{\cal  W} \, \gua}) \ga \\ 
&=& \tau_{+}(\nabla_{\cal V} {\cal  W}) 
    -\gwii{(\nabla_{\cal V} {\cal  W}) \, \gua} \ga 
    - \gwii{{\cal V} \, {\cal  W} \, \gua} \ga \\ 
&=& T(\nabla_{\cal V} {\cal  W})  - \qpcc{V}{W}. 
\end{eqnarray*} 
This proves (3). 
$\Box$ 
 
\begin{cor} \label{cor:4ptT} 
For any vector fields ${\cal W}_{1}$, ${\cal W}_{2}$, and 
${\cal V}$, 
\[ 
 \gwii{{\cal W}_{1} \, {\cal W}_{2} \, T({\cal V}) \, \gua} \ga 
= {\cal W}_{1} \bullet {\cal W}_{2} \bullet {\cal V}. 
\] 
\end{cor} 
{\bf Proof}: 
By equation~(\ref{eqn:DerProd}), 
\begin{eqnarray*} 
&& \gwii{{\cal W}_{1} \, {\cal W}_{2} \, T({\cal V}) \, \gua} \ga 
= \nabla_{{\cal W}_{1}} ({\cal W}_{2} \bullet T({\cal V})) 
        - (\nabla_{{\cal W}_{1}} {\cal W}_{2}) \bullet T({\cal V}) 
        - {\cal W}_{2} \bullet (\nabla_{{\cal W}_{1}} T({\cal V})). 
\end{eqnarray*} 
The first two terms on the right hand side vanish since 
$T({\cal V}) \sim 0$. The corollary then follows from Lemma~\ref{lem:DerT}. 
$\Box$ 
 
As a consequence of \eqref{eqn:DerProd}, we also have 
\begin{cor} \label{cor:TderProd} 
For any vector fields ${\cal W}_{1}$, ${\cal W}_{2}$, and 
${\cal V}$, 
\[ 
 \nabla_{T(\vv)} ({\cal W}_{1} \bullet {\cal W}_{2}) 
=  (\nabla_{T(\vv)} {\cal W}_{1}) \bullet {\cal W}_{2} 
   +  {\cal W}_{1} \bullet (\nabla_{T(\vv)}  {\cal W}_{2}) 
    +{\cal W}_{1} \bullet {\cal W}_{2} \bullet {\cal V}. 
\] 
\end{cor} 
Since $T(\vv) \sim 0$, this corollary implies that 
$\nabla_{T^{2}(\vv)}$ is a derivation of the quantum product 
 for any vector field $\vv$. 
Since covariant derivatives preserve the space of primary 
vector fields, $\nabla_{T^{2}(\vv)}$ and the quantum product 
produce a holomorphic vertex algebra structure on this space 
(cf. \cite{Kac}). 
 
Taking derivatives of the equation in Corollary~\ref{cor:4ptT}, we have 
\begin{eqnarray} 
 \gwii{{\cal W}_{1} \, {\cal W}_{2} \, \vw_{3} \, T({\cal V}) \, \gua} \ga 
& = & \gwii{ \{\vv \bullet \vw_{1} \} \, \vw_{2} \, \vw_{3} \, \gua} \ga 
\nonumber \\ 
&& + \gwii{ \{\vv \bullet \vw_{2} \} \, \vw_{1} \, \vw_{3} \, \gua} \ga 
\nonumber \\ 
&&  + \gwii{ \vv \, \vw_{1} \, \vw_{2} \, \{\vw_{3} \bullet \gua\} } \ga 
\label{eqn:5ptT} 
\end{eqnarray} 
and 
\begin{eqnarray} 
 \gwii{{\cal W}_{1} \, {\cal W}_{2} \, \vw_{3} \, \vw_{4} \, 
        T({\cal V}) \, \gua} \ga 
& = & \sum_{i=1}^{3} 
\gwii{ \{\vv \bullet \vw_{i} \} \, \vw_{1} \, \ldots \, \widehat{\vw_{i}} 
    \, \ldots \, \vw_{4} \, \gua} \ga 
\nonumber \\ 
&&  + \gwii{ \vv \, \vw_{1} \, \vw_{2} \, \vw_{3} \, 
         \{\vw_{4} \bullet \gua\} } \ga 
\nonumber \\ 
&& + \gwii{ \vv \, \vw_{1} \, \vw_{2} \, \gub} 
    \gwii{ \gb \, \vw_{3} \, \vw_{4} \, \gua } \ga 
\nonumber \\ 
&& + \gwii{ \vv \, \vw_{1} \, \vw_{3} \, \gub} 
    \gwii{ \gb \, \vw_{2} \, \vw_{4} \, \gua } \ga 
\nonumber \\ 
&& + \gwii{ \vv \, \vw_{2} \, \vw_{3} \, \gub} 
    \gwii{ \gb \, \vw_{1} \, \vw_{4} \, \gua } \ga 
\label{eqn:6ptT} 
\end{eqnarray} 
for all vector fields. These formulas will be useful later when we 
study the genus-2 Virasoro conjecture. 
For $k \geq 1$, an induction on $k$ shows that 
\begin{equation} \label{eqn:derTRRg0} 
 \gwii{ T^{k}(\vw) \, \vv_{1} \, \cdots \, \vv_{k} \, \gua} \ga = 0. 
\end{equation}

We collect some useful formulas for the string vector field in 
the following 
\begin{lem} \label{lem:String} 
\begin{eqnarray*} 
& (1) & \nabla_{\cal W} {\cal S} = - \tau_{-}({\cal W}). \\ 
& (2) & \gwii{ {\cal S} \, {\cal V} \, {\cal W} \, \gua} 
    = \gwii{ \tau_{-}({\cal V}) \, {\cal W} \, \gua} 
     + \gwii{ {\cal V} \, \tau_{-}({\cal W}) \, \gua}. \\ 
& (3) & {\cal W} \gwii{ {\cal S} \, {\cal S} \, \gua} 
    =  \gwii{ {\cal W} \, \tau_{-}({\cal S}) \, \gua} 
       -\gwii{ \tau_{-}({\cal W}) \, {\cal S} \, \gua}. 
\end{eqnarray*} 
\end{lem} 
{\bf Proof}: The first formula follows from the definition 
of $\cal S$. The genus-0 string equation implies that 
\[ \gwii{ {\cal S} \, \gravm{m} \, \gravn{n} \, \gua} 
    = \gwii{ \gravm{m-1} \, \gravn{n} \, \gua} 
        + \gwii{ \gravm{m} \, \gravn{n-1} \, \gua} \] 
for all $\gravm{m}$ and $\gravn{n}$. By linearity of the correlation 
functions, this implies the second formula. 
The last formula follows from the first two formulas and 
equation~(\ref{eqn:DerCorr}). 
$\Box$ 
 
Lemma~\ref{lem:String} (2) and \eqref{eqn:DerProd} implies 
$\nabla_{\vs}$ is a derivation of the quantum product 
on the space of primary vector fields and therefore 
also produce a holomorphic vertex algebra structure. 
 
Using Lemma~\ref{lem:String} (1) and Lemma~\ref{lem:DerT} (3), 
we obtain 
\[ \nabla_{T^{m}(\vs)} T^{k}(\vs) = - T^{k+m-1}(\vs). \] 
Therefore, 
\[ [T^{k}(\vs), \, T^{m}(\vs)] = 0 \] 
for all $k, m \geq 0$. 
 
Another application of Lemma~\ref{lem:String} (1) 
is that 
\begin{equation} \label{eqn:DerStr} 
 \gwiig{ \vs \, \vw_{1} \, \cdots \, \vw_{k} } 
    = \sum_{i=1}^{k} \gwiig{ \vw_{1} \, \cdots \, \left\{ \tau_{-}(\vw_{i}) \right\} 
            \, \cdots \, \vw_{k} } 
            + \delta_{g, 0} \nabla^{k}_{\vw_{1}, \cdots, \vw_{k}} 
                \left( \frac{1}{2} \eta_{\alpha \beta} t_{0}^{\alpha} t_{0}^{\beta} \right) 
\end{equation} 
for any vector fields $\vw_{1}, \ldots, \vw_{k}$, where $\nabla^{k}$ is the 
$k$-th covariant derivative. This formula is obtained by repeatedly 
taking derivatives on both sides of the string equation and applying 
\eqref{eqn:DerCorr}.

\section{Applications to topological recursion relations} 
\label{sec:TRR} 
 
It seems that the quantum product on the big phase space provides 
appropriate algebraic machinery for studying problems on the big 
phase space, e.g. the Virasoro conjecture 
and topological recursion relations. We first consider 
the topological recursion relations in this section and postpone 
the Virasoro conjecture to later sections. 
 
Lemma~\ref{lem:0vector} indicates that the transformation $T$ 
trivializes vector fields at the genus-0 level. In some sense, 
all known topological recursion relations seem to indicate 
how this operator trivialize vector fields for various genera. 
Let us first see the cases of genus less than or equal to 2. 
 
Since 
\[ T(\grava{n}) = \grava{n+1} - \gwii{ \grava{n} \, \gum} \gm, \] 
the coefficient of $\gn$ for $T(\grava{n}) \bullet \gravb{m}$ is 
\[ \gwii{\grava{n+1} \, \gravb{m} \, \gun} - 
    \gwii{\grava{n} \, \gum} \gwii{\gm \, \gravb{m} \, \gun}. \] 
The vanishing of this quantity is the most important case of the 
genus-0 topological recursion relation (in particular, it implies 
the associativity of the quantum product on the big phase space). 
Therefore $T(\grava{n}) \sim 0$ is equivalent to this 
special case of the genus-0 topological recursion relation. 
 
The genus-1 topological recursion relation is the following: 
\[ \gwiione{\grava{n+1}} = \gwii{\grava{n} \, \gum} \gwiione{\gm} 
    + \frac{1}{24} \gwii{ \grava{n} \, \gum  \, \gm}. \] 
This formula is equivalent to 
\begin{equation} \label{eqn:TRRg1} 
\gwiione{T({\cal W})} \, \, = \, \, 
 \frac{1}{24} \gwii{{\cal W} \, \gum  \, \gm} 
\end{equation} 
for any vector field ${\cal W}$. For $g > 0$, 
we call a vector field ${\cal W}$ {\it trivial at the genus-$g$ level} 
if $\gwiig{\cal W}$ can be represented by data of genera less than $g$. 
Then the genus-1 topological recursion relation just means 
that $T({\cal W})$ is trivial at the genus-1 level for all $\cal W$. 
 
Taking derivatives of \eqref{eqn:TRRg1} and using Lemma~\ref{lem:DerT}, 
we have 
\begin{equation} 
\gwiione{T({\cal W}) \, \vv} \, \, = \, \, 
    \gwiione{ \{\vw \bullet \vv\} } + 
     \frac{1}{24} \gwii{{\cal W} \, \vv \, \gum  \, \gm} 
\label{eqn:derTRRg1} 
\end{equation} 
and 
\begin{eqnarray} 
\gwiione{T({\cal W}) \, \vv_{1} \, \vv_{2}} 
& = &   \gwiione{ \{\vw \bullet \vv_{1}\} \, \vv_{2} } + 
    \gwiione{ \{\vw \bullet \vv_{2}\} \, \vv_{1} } + 
    \gwii{ \vw \, \vv_{1} \, \vv_{2} \, \gua } \gwiione{\ga} 
\nonumber \\ 
&&  + \frac{1}{24} \gwii{{\cal W} \, \vv_{1} \, \vv_{2} \, \gum  \, \gm} 
\label{eqn:2derTRRg1} 
\end{eqnarray} 
for all vector fields. These formulas will be used later. 
 
The genus-2 topological recursion relations are much more complicated 
than genus-0 and genus-1 topological recursion relations. The following 
two recursion relations were found in \cite{G1}: 
\begin{eqnarray} 
&& \gwiitwo{\tau_{i+2}(x)} \, \, = \, \, 
 \gwii{\tau_{i+1}(x) \, \gua} \gwiitwo{\ga} 
    + \gwii{\tau_{i}(x) \, \gua} \gwiitwo{\tau_{1}(\ga)} 
    \nonumber \\ 
&& \hspace{30pt} 
    - \gwii{\tau_{i}(x) \, \gua} \gwii{\ga \, \gub} \gwiitwo{\gb} 
    + \frac{7}{10} \gwii{\tau_{i}(x) \, \gua \, \gub} 
        \gwiione{\ga} \gwiione{\gb} 
    \nonumber \\ 
&& \hspace{30pt} 
    + \frac{1}{10} \gwii{\tau_{i}(x) \, \gua \, \gub} 
        \gwiione{\ga \, \gb} 
    - \frac{1}{240} \gwiione{\tau_{i}(x) \, \ga} 
        \gwii{\gua \, \gb \, \gub} 
    \nonumber \\ 
&&  \hspace{30pt} 
    + \frac{13}{240} \gwii{\tau_{i}(x) \, \ga \, \gua \, \gub} 
        \gwiione{\gb} 
    + \frac{1}{960} \gwii{ \tau_{i}(x) \, \gua \, \ga \, \gub \, \gb} 
\label{eqn:TRR1-old} 
\end{eqnarray} 
for $i \geq 0$ and $x \in \{ \gamma_{1}, \ldots, \gamma_{N}\}$, and 
\begin{eqnarray} 
&& \hspace{-15pt} \gwiitwo{\tau_{i+1}(x) \, \tau_{j+1}(y)} \, \, = \, \, 
    \gwiitwo{ \tau_{i+1}(x) \, \gua} \gwii{ \ga \, \tau_{j}(y)} 
    + \gwii{\tau_{i}(x) \, \gua} \gwiitwo{\ga \, \tau_{j+1}(y)} 
    \nonumber \\ 
&&  \hspace{15pt} 
    - \gwii{\tau_{i}(x) \, \gua} \gwii{\tau_{j}(y) \, \gub} 
        \gwiitwo{\ga \, \gb} 
    + 3 \gwii{\tau_{i}(x) \, \tau_{j}(y) \, \gua} 
        \gwiitwo{ \tau_{1}(\ga)} 
    \nonumber \\ 
&&  \hspace{15pt} 
    -3 \gwii{\tau_{i}(x) \, \tau_{j}(y) \, \gua} 
        \gwii{\ga \, \gub} \gwiitwo{ \gb} 
    + \frac{13}{10} \gwii{\tau_{i}(x) \, \tau_{j}(y) \, \gua \, \gub} 
        \gwiione{\ga} \gwiione{\gb} 
    \nonumber \\ 
&&  \hspace{15pt} 
    + \frac{4}{5} \gwiione{\tau_{i}(x) \, \ga} 
        \gwii{ \gua \, \tau_{j}(y) \, \gub} \gwiione{\gb} 
    + \frac{4}{5} \gwii{\tau_{i}(x) \, \gua \, \gub} 
        \gwiione{\tau_{j}(y) \, \ga} \gwiione{\gb} 
    \nonumber \\ 
&&  \hspace{15pt} 
    - \frac{4}{5} \gwii{\tau_{i}(x) \, \tau_{j}(y) \, \gua} 
        \gwiione{ \ga \, \gub} \gwiione{\gb} 
    + \frac{23}{240} \gwii{ \tau_{i}(x) \, \tau_{j}(y) \, 
        \gua \, \ga \, \gub} \gwiione{\gb} 
    \nonumber \\ 
&&  \hspace{15pt} 
    + \frac{1}{48}\gwii{ \tau_{i}(x) \, \gua \, \ga \, \gub} 
        \gwiione{ \gb \,\tau_{j}(y)} 
    + \frac{1}{48}\gwii{ \tau_{j}(y) \, \gua \, \ga \, \gub} 
        \gwiione{ \gb \,\tau_{i}(x)} 
    \nonumber \\ 
&&  \hspace{15pt} 
    - \frac{1}{80} \gwiione{\tau_{i}(x) \, \tau_{j}(y) \, \gua} 
        \gwii{ \ga \, \gub \, \gb} 
    + \frac{7}{30} \gwii{\tau_{i}(x) \, \tau_{j}(y) \, \gua \, \gub} 
        \gwiione{ \ga \, \gb} 
    \nonumber \\ 
&&  \hspace{15pt} 
    + \frac{1}{30} \gwii{\tau_{i}(x) \, \gua \, \gub} 
        \gwiione{\ga \, \gb \, \tau_{j}(y)} 
    + \frac{1}{30} \gwii{\tau_{j}(y) \, \gua \, \gub} 
        \gwiione{\ga \, \gb \, \tau_{i}(x)} 
    \nonumber \\ 
&&  \hspace{15pt} 
    - \frac{1}{30} \gwii{\tau_{i}(x) \, \tau_{j}(y) \, \gua} 
        \gwiione{\ga \, \gub \, \gb} 
    + \frac{1}{576}\gwii{\tau_{i}(x) \, \tau_{j}(y) \, \gua \, \ga 
        \, \gub \, \gb} 
\label{eqn:TRR2-old} 
\end{eqnarray} 
for $i, j \geq 0$ and $x, y \in \{ \gamma_{1}, \ldots, \gamma_{N}\}$. 
Another genus-2 recursion relation was found in \cite{BP}: 
\begin{eqnarray} 
&& \hspace{-15pt} 0 = \sum_{\sigma \in S_{3}} 
    - 2 \gwii{\tau_{i_{\sigma(1)}}(x_{\sigma(1)}) \, \, \ga \, \gub} 
        \gwii{ \gua \, \, \tau_{i_{\sigma(2)}}(x_{\sigma(2)}) \, \, 
            \tau_{i_{\sigma(3)}}(x_{\sigma(3)})} 
        \gwiitwo{\gb} 
    \nonumber \\ 
&&  \hspace{30pt} 
+ 2 \gwii{\tau_{i_{\sigma(1)}}(x_{\sigma(1)}) \, \, 
        \tau_{i_{\sigma(2)}}(x_{\sigma(2)}) \,  \, 
        \tau_{i_{\sigma(3)}}(x_{\sigma(3)}) \, \, \gua} 
    \gwiitwo{\tau_{1}(\ga)} 
    \nonumber \\ 
&&  \hspace{30pt} 
- 2 \gwii{\tau_{i_{\sigma(1)}}(x_{\sigma(1)}) \, \, 
        \tau_{i_{\sigma(2)}}(x_{\sigma(2)}) \, \, 
        \tau_{i_{\sigma(3)}}(x_{\sigma(3)}) \, \, \gua} 
    \gwii{\ga \, \gub} \gwiitwo{\gb} 
    \nonumber \\ 
&&  \hspace{30pt} 
+ 3 \gwiitwo{ \tau_{1+i_{\sigma(1)}}(x_{\sigma(1)}) \, \, \ga} 
    \gwii{ \gua \, \, \tau_{i_{\sigma(2)}}(x_{\sigma(2)}) \, \, 
            \tau_{i_{\sigma(3)}}(x_{\sigma(3)})} 
    \nonumber \\ 
&&  \hspace{30pt} 
- 3 \gwii{ \tau_{i_{\sigma(1)}}(x_{\sigma(1)}) \, \, \gua} 
    \gwii{\tau_{i_{\sigma(2)}}(x_{\sigma(2)}) \, \, 
        \tau_{i_{\sigma(3)}}(x_{\sigma(3)}) \, \, \gub} 
    \gwiitwo{\ga \gb} 
    \nonumber \\ 
&&  \hspace{30pt} 
- 3  \gwiitwo{ \tau_{i_{\sigma(1)}}(x_{\sigma(1)}) \, \, \tau_{1}(\ga)} 
    \gwii{ \gua \, \, \tau_{i_{\sigma(2)}}(x_{\sigma(2)}) \, \, 
            \tau_{i_{\sigma(3)}}(x_{\sigma(3)})} 
    \nonumber \\ 
&&  \hspace{30pt} 
+ 3 \gwiitwo{ \tau_{i_{\sigma(1)}}(x_{\sigma(1)}) \, \, \ga} 
    \gwii{\gua \gb} \gwii{\gub \, \, \tau_{i_{\sigma(2)}}(x_{\sigma(2)}) \, 
        \, \tau_{i_{\sigma(3)}}(x_{\sigma(3)})} 
    \nonumber \\ 
&&  \hspace{30pt} 
+ \frac{1}{5} \gwii{\tau_{i_{\sigma(1)}}(x_{\sigma(1)}) \, \, 
        \tau_{i_{\sigma(2)}}(x_{\sigma(2)}) \, \, 
        \tau_{i_{\sigma(3)}}(x_{\sigma(3)}) \, \, \gua \, \gub} 
    \gwiione{\ga} \gwiione{\gb} 
    \nonumber \\ 
&&  \hspace{30pt} 
- \frac{6}{5} \gwii{\tau_{i_{\sigma(1)}}(x_{\sigma(1)}) \, \, 
        \tau_{i_{\sigma(2)}}(x_{\sigma(2)}) \, \, \gua \, \gub} 
    \gwiione{\ga} \gwiione{\gb \, \, \tau_{i_{\sigma(3)}}(x_{\sigma(3)})} 
    \nonumber \\ 
&&  \hspace{30pt} 
+ \frac{12}{5}\gwii{\tau_{i_{\sigma(1)}}(x_{\sigma(1)}) \, \, \gua \, \gub} 
    \gwiione{\ga} \gwiione{\gb \, \, \tau_{i_{\sigma(2)}}(x_{\sigma(2)}) \, 
        \, \tau_{i_{\sigma(3)}}(x_{\sigma(3)})} 
    \nonumber \\ 
&&  \hspace{30pt} 
- \frac{18}{5}\gwiione{\tau_{i_{\sigma(1)}}(x_{\sigma(1)}) \, \, \ga} 
    \gwii{\gua \, \, \tau_{i_{\sigma(2)}}(x_{\sigma(2)}) \, \, \gub} 
    \gwiione{\gb \, \, \tau_{i_{\sigma(3)}}(x_{\sigma(3)})} 
    \nonumber \\ 
&&  \hspace{30pt} 
- \frac{6}{5} \gwii{\tau_{i_{\sigma(1)}}(x_{\sigma(1)}) \, \, 
        \tau_{i_{\sigma(2)}}(x_{\sigma(2)}) \, \, 
        \tau_{i_{\sigma(3)}}(x_{\sigma(3)}) \, \, \gua} 
    \gwiione{\ga \, \gub} \gwiione{\gb} 
    \nonumber \\ 
&&  \hspace{30pt} 
+ \frac{9}{5}\gwiione{\tau_{i_{\sigma(1)}}(x_{\sigma(1)}) \, \, \ga} 
    \gwiione{\gua \, \gb} 
    \gwii{\gub \, \, \tau_{i_{\sigma(2)}}(x_{\sigma(2)}) \, \, 
            \tau_{i_{\sigma(3)}}(x_{\sigma(3)})} 
    \nonumber \\ 
&&  \hspace{30pt} 
- \frac{6}{5}\gwiione{\ga} 
    \gwiione{\gua \, \, \tau_{i_{\sigma(1)}}(x_{\sigma(1)}) \, \, \gb} 
    \gwii{\gub \, \, \tau_{i_{\sigma(2)}}(x_{\sigma(2)}) \,  \, 
            \tau_{i_{\sigma(3)}}(x_{\sigma(3)})} 
    \nonumber \\ 
&&  \hspace{30pt} 
+ \frac{1}{120} \gwii{\tau_{i_{\sigma(1)}}(x_{\sigma(1)}) \, \, 
        \tau_{i_{\sigma(2)}}(x_{\sigma(2)}) \, \, 
        \tau_{i_{\sigma(3)}}(x_{\sigma(3)}) \, \, \gua \, \ga \, \gub} 
    \gwiione{\gb} 
    \nonumber \\ 
&&  \hspace{30pt} 
- \frac{3}{40} \gwii{\tau_{i_{\sigma(1)}}(x_{\sigma(1)}) \, \, 
        \tau_{i_{\sigma(2)}}(x_{\sigma(2)}) \, \, \gua \, \ga \, \gub} 
    \gwiione{\gb \, \, \tau_{i_{\sigma(3)}}(x_{\sigma(3)})} 
    \nonumber \\ 
&&  \hspace{30pt} 
+ \frac{3}{40} \gwii{\tau_{i_{\sigma(1)}}(x_{\sigma(1)}) \, \, 
            \gua \, \ga \, \gub} 
    \gwiione{\gb \, \, \tau_{i_{\sigma(2)}}(x_{\sigma(2)}) \, \, 
            \tau_{i_{\sigma(3)}}(x_{\sigma(3)})} 
    \nonumber \\ 
&&  \hspace{30pt} 
- \frac{1}{120} \gwiione{\tau_{i_{\sigma(1)}}(x_{\sigma(1)}) \, \, 
        \tau_{i_{\sigma(2)}}(x_{\sigma(2)}) \, \, 
        \tau_{i_{\sigma(3)}}(x_{\sigma(3)}) \, \, \ga} 
    \gwii{\gua \, \gub \, \gb} 
    \nonumber \\ 
&&  \hspace{30pt} 
+ \frac{1}{10} \gwii{\tau_{i_{\sigma(1)}}(x_{\sigma(1)}) \, \, 
        \tau_{i_{\sigma(2)}}(x_{\sigma(2)}) \, \, 
        \tau_{i_{\sigma(3)}}(x_{\sigma(3)}) \, \, \gua \, \gub} 
    \gwiione{\ga \, \gb} 
    \nonumber \\ 
&&  \hspace{30pt} 
- \frac{3}{10} \gwiione{\tau_{i_{\sigma(1)}}(x_{\sigma(1)}) \, \, \ga \, \gb} 
        \gwii{\gua \, \gub \, \, \tau_{i_{\sigma(2)}}(x_{\sigma(2)}) \, 
        \, \tau_{i_{\sigma(3)}}(x_{\sigma(3)})} 
    \nonumber \\ 
&&  \hspace{30pt} 
+ \frac{1}{10} \gwiione{\tau_{i_{\sigma(1)}}(x_{\sigma(1)}) \, \, 
            \tau_{i_{\sigma(2)}}(x_{\sigma(2)}) \, \, \ga \, \gb} 
        \gwii{\gua \, \gub \, \, \tau_{i_{\sigma(3)}}(x_{\sigma(3)})} 
    \nonumber \\ 
&&  \hspace{30pt} 
- \frac{1}{20} \gwii{\tau_{i_{\sigma(1)}}(x_{\sigma(1)}) \, \, 
        \tau_{i_{\sigma(2)}}(x_{\sigma(2)}) \, \, 
        \tau_{i_{\sigma(3)}}(x_{\sigma(3)}) \,  \, \gua} 
    \gwiione{\gua \, \gub \, \gb} 
    \nonumber \\ 
&&  \hspace{30pt} 
- \frac{1}{20} \gwiione{\tau_{i_{\sigma(1)}}(x_{\sigma(1)}) \, \, 
        \gua \, \ga \, \gb} 
    \gwii{\gub \, \, \tau_{i_{\sigma(2)}}(x_{\sigma(2)}) \, \, 
        \tau_{i_{\sigma(3)}}(x_{\sigma(3)})}, 
\label{eqn:BP-old} 
\end{eqnarray} 
where $i_{1}, i_{2}, i_{3} \geq 0$, 
$x_{1}, x_{2}, x_{3} \in \{ \gamma_{1}, \ldots, \gamma_{N}\}$, 
and $S_{3}$ is the symmetry group of three elements. 
 
Unlike the genus-1 case, $T({\cal W})$ is no longer trivial at the genus-2 
level in general. In fact, the genus-2 dilaton equation implies 
\[ \gwiitwo{T({\cal S})} = \gwiitwo{\cal D} = 2 F_{2}. \] 
Unless that $F_{2}$ can be expressed as a function of $F_{0}$ and 
$F_{1}$, $T({\cal S})$ is not trivial at the genus-2 level. 
However, the topological recursion relation (\ref{eqn:TRR1-old}) implies 
that $T^{2}({\cal W}) := T(T({\cal W}))$ is trivial at the genus-2 level for 
all vector field $\cal W$. 
 
The genus-2 topological recursion relations 
(\ref{eqn:TRR1-old}) - (\ref{eqn:BP-old}) 
 can be represented in the following forms respectively: 
For any vector fields $\cal W$, $\cal V$, ${\cal W}_{1}$, 
${\cal W}_{2}$, and ${\cal W}_{3}$ on the big phase space, 
\begin{eqnarray} 
&& \gwiitwo{T^{2}({\cal W})} = A_{1}(\vw), \hspace{200pt} \label{eqn:TRR1} \\ 
&& 
\gwiitwo{T({\cal W}) \, T({\cal V})} - 3 \gwiitwo{ T(\qpcc{W}{V})} 
    = A_{2}(\vw, \vv), 
\label{eqn:TRR2} 
\end{eqnarray} 
and 
\begin{eqnarray} 
&& 2 \gwiitwo{\{{\cal W}_{1} \bullet {\cal W}_{2} \bullet {\cal W}_{3} \}} 
- 2 \gwii{{\cal W}_{1} \, {\cal W}_{2} \, {\cal W}_{3} \, \gua} 
    \gwiitwo{T(\ga)} 
\nonumber \\ 
&-& \gwiitwo{ T({\cal W}_{1}) \, \{{\cal W}_{2} \bullet {\cal W}_{3}\}} 
    + \gwiitwo{{\cal W}_{1} \, T({\cal W}_{2} \bullet {\cal W}_{3})} 
\nonumber \\ 
&-& \gwiitwo{ T({\cal W}_{2}) \, \{{\cal W}_{1} \bullet {\cal W}_{3}\}} 
    + \gwiitwo{{\cal W}_{2} \, T({\cal W}_{1} \bullet {\cal W}_{3})} 
\nonumber \\ 
&-& \gwiitwo{ T({\cal W}_{3}) \, \{{\cal W}_{1} \bullet {\cal W}_{2}\}} 
    + \gwiitwo{{\cal W}_{3} \, T({\cal W}_{1} \bullet {\cal W}_{2})} 
\, \, = \, \, B(\vw_{1}, \vw_{2}, \vw_{3}), 
\label{eqn:BP} 
\end{eqnarray} 
where 
\begin{eqnarray*} 
A_{1}(\vw) &=& 
    \frac{7}{10} \gwiione{\ga} \gwiione{\{\gua \bullet \vw\}} 
    + \frac{1}{10} \gwiione{\ga \, \{\gua \bullet \vw\}} \\ 
&&     - \frac{1}{240} \gwiione{\vw \, \{\ga \bullet \gua\}} 
    + \frac{13}{240} \gwii{\vw \, \ga \, \gua \, \gub} 
        \gwiione{\gb} \\ 
&&    + \frac{1}{960} \gwii{ \vw \, \gua \, \ga \, \gub \, \gb}, \\ 
A_{2}(\vw, \vv) &=& 
     \frac{13}{10} \gwii{\vw \, \vv \, \gua \, \gub} 
        \gwiione{\ga} \gwiione{\gb} 
    + \frac{4}{5} \gwiione{\vw \, \gua} 
            \gwiione{\{\ga \bullet \vv\}} \\ 
&&    + \frac{4}{5} 
        \gwiione{\vv \, \gua} \gwiione{\{\ga \bullet \vw\}} 
    - \frac{4}{5} 
        \gwiione{ \{\vw \bullet \vv\}  \, \gua} \gwiione{\ga} \\ 
&&    + \frac{23}{240} \gwii{ \vw \, \vv \, 
        \gua \, \ga \, \gub} \gwiione{\gb} 
    + \frac{1}{48}\gwii{ \vw \, \gua \, \ga \, \gub} 
        \gwiione{ \gb \,\vv}  \\ 
&&    + \frac{1}{48}\gwii{ \vv \, \gua \, \ga \, \gub} 
        \gwiione{ \gb \,\vw} 
    - \frac{1}{80} \gwiione{\vw \, \vv \, \{\gua \bullet \ga\}} \\ 
&&    + \frac{7}{30} \gwii{\vw \, \vv \, \gua \, \gub} 
        \gwiione{ \ga \, \gb} 
    + \frac{1}{30} 
        \gwiione{\ga \, \{\gua \bullet \vw\} \, \vv} \\ 
&&    + \frac{1}{30} 
        \gwiione{\ga \, \{\gua \bullet \vv\} \, \vw} 
    - \frac{1}{30} 
        \gwiione{\{ \vw \bullet \vv\} \ga \, \gua} \\ 
&&    + \frac{1}{576}\gwii{\vw \, \vv \, \gua \, \ga 
        \, \gub \, \gb}, 
\end{eqnarray*} 
and 
\begin{eqnarray*} 
&& B(\vw_{1}, \vw_{2}, \vw_{3}) \\ 
&=& \frac{1}{5} \gwii{\vw_{1} \, \vw_{2} \, \vw_{3} \, \, \gua \, \gub} 
    \gwiione{\ga} \gwiione{\gb} 
   - \frac{6}{5} \gwii{\vw_{1} \, \vw_{2} \, \vw_{3} \, \gua} 
    \gwiione{\ga \, \gub} \gwiione{\gb} 
    \nonumber \\ 
&& 
+ \frac{1}{120} \gwii{\vw_{1} \,   \vw_{2} \, \vw_{3} 
         \, \gua \, \ga \, \gub} 
    \gwiione{\gb} 
- \frac{1}{120} \gwiione{\vw_{1} \, \vw_{2} \, \vw_{3} 
       \, \{ \gua \bullet \ga \}} 
    \nonumber \\ 
&& 
+ \frac{1}{10} \gwii{\vw_{1} \, \vw_{2} \, \vw_{3} 
         \, \gua \, \gub} 
    \gwiione{\ga \, \gb} 
- \frac{1}{20} \gwii{\vw_{1} \, \vw_{2} \, \vw_{3} \, \gua} 
    \gwiione{\ga \, \gub \, \gb} 
    \nonumber \\ 
&& 
- \frac{1}{5} \sum_{\sigma \in S_{3}} 
     \gwii{\vw_{\sigma(1)} \, \vw_{\sigma(2)} \, \gua \, \gub} 
    \gwiione{\ga} \gwiione{\gb \, \, \vw_{\sigma(3)}} 
    \nonumber \\ 
&& 
+ \frac{2}{5}\sum_{\sigma \in S_{3}} 
    \gwiione{\{\vw_{\sigma(1)} \bullet \ga \}} 
    \gwiione{\gua \, \, \vw_{\sigma(2)} \, 
        \, \vw_{\sigma(3)}} 
    \nonumber \\ 
&& 
- \frac{3}{5} \sum_{\sigma \in S_{3}} 
    \gwiione{\vw_{\sigma(1)} \, \, \{ \vw_{\sigma(2)} \bullet \gua \}} 
    \gwiione{\ga \, \, \vw_{\sigma(3)}} 
    \nonumber \\ 
&& 
+ \frac{3}{10}\sum_{\sigma \in S_{3}} 
    \gwiione{\vw_{\sigma(1)} \, \, \ga} 
    \gwiione{\gua \, \{\vw_{\sigma(2)} \bullet \vw_{\sigma(3)} \}} 
    \nonumber \\ 
&& 
- \frac{1}{5}\sum_{\sigma \in S_{3}} 
    \gwiione{\ga} 
    \gwiione{\gua \, \, \vw_{\sigma(1)} \, \, 
        \{\vw_{\sigma(2)} \bullet \vw_{\sigma(3)}\}} 
    \nonumber \\ 
&& 
- \frac{1}{80} \sum_{\sigma \in S_{3}} 
    \gwii{\vw_{\sigma(1)} \, \, 
        \vw_{\sigma(2)} \, \, \gua \, \ga \, \gub} 
    \gwiione{\gb \, \, \vw_{\sigma(3)}} 
    \nonumber \\ 
&& 
+ \frac{1}{80} \sum_{\sigma \in S_{3}} 
    \gwii{\vw_{\sigma(1)} \, \, 
            \gua \, \ga \, \gub} 
    \gwiione{\gb \, \, \vw_{\sigma(2)} \, \, 
            \vw_{\sigma(3)}} 
    \nonumber \\ 
&& 
- \frac{1}{20} \sum_{\sigma \in S_{3}} 
     \gwiione{\vw_{\sigma(1)} \, \, \ga \, \gb} 
        \gwii{\gua \, \gub \, \, \vw_{\sigma(2)} \, 
        \, \vw_{\sigma(3)}} 
    \nonumber \\ 
&& 
+ \frac{1}{60} \sum_{\sigma \in S_{3}} 
    \gwiione{\vw_{\sigma(1)} \, \, 
            \vw_{\sigma(2)} \, \, \ga \, \{\gua \bullet \vw_{\sigma(3)} \}} 
    \nonumber \\ 
&& 
- \frac{1}{120} \sum_{\sigma \in S_{3}} 
     \gwiione{\vw_{\sigma(1)} \, \, 
        \gua \, \ga \, \{ \vw_{\sigma(2)} \bullet \vw_{\sigma(3)} \} }. 
\end{eqnarray*} 
In particular, 
$A_{1}$, $A_{2}$, and $B$ are symmetric tensors only depending on 
genus-0 and genus-1 data.

We now study relations among these genus-2 recursion relations. 
We first observe that equation 
(\ref{eqn:TRR1}) is equivalent to 
\begin{equation} \label{eqn:TRR1-2} 
\gwiitwo{T({\cal W})} \, = \, \gwiitwo{T(\bvw)} 
    + A_{1}(\tau_{-}(\vw)). 
\end{equation} 
for any vector field ${\cal W}$. 
The following two lemmas are formal consequences of 
(\ref{eqn:TRR1}) or its equivalent form (\ref{eqn:TRR1-2}). 
\begin{lem} \label{lem:T2VW} 
Equation (\ref{eqn:TRR1}) implies 
\[ \gwiitwo{T^{2}({\cal V}) \, \, {\cal W}} 
    \, = \,  \gwiitwo{T(\qpcc{V}{W})} 
    +  \left(\nabla_{\vw} A_{1} \right)(\vv) \] 
for any vector fields ${\cal V}$ and ${\cal W}$, where $\nabla_{\vw} A_{1}$ 
is the covariant derivative of $A_{1}$ defined by 
\[ \left(\nabla_{\vw} A_{1} \right)(\vv) = \vw \, A_{1}(\vv) - A_{1}(\nabla_{\vw} \vv). \] 
\end{lem} 
{\bf Proof}: 
By equation~(\ref{eqn:DerCorr}), 
\[ \gwiitwo{T^{2}({\cal V}) \, \, {\cal W}} 
    = {\cal W} \gwiitwo{T^{2}({\cal V})} - 
        \gwiitwo{ \left\{\nabla_{\cal W} T^{2}({\cal V}) \right\}}. 
\] 
By equation~(\ref{eqn:TRR1}) and Lemma~\ref{lem:DerT}, we have 
\[ \gwiitwo{T^{2}({\cal V}) \, \, {\cal W}} 
    \, = \, \vw \, A_{1}(\vv) -\gwiitwo{ T(\nabla_{\cal W} T({\cal V}))} 
      + \gwiitwo{ \left\{ {\cal W} \bullet T({\cal V}) \right\}}. \] 
The last term on the right hand side vanishes since $T({\cal V}) \sim 0$. 
Using Lemma~\ref{lem:DerT} again, we have 
\[ \gwiitwo{T^{2}({\cal V}) \, \, {\cal W}} 
    \, = \, \vw \, A_{1}(\vv) -\gwiitwo{ T^{2}(\nabla_{\cal W} {\cal V})} 
      + \gwiitwo{ T({\cal W} \bullet {\cal V}) }. \] 
The lemma then follows from equation~(\ref{eqn:TRR1}). 
$\Box$ 
 
\begin{cor} \label{cor:TVTSW} 
Equation (\ref{eqn:TRR1}) implies 
\[ \gwiitwo{T({\cal V}) \, \, T(\bvw)} 
    = \gwiitwo{T({\cal V}) \, \, T({\cal W})} 
     - \left(\nabla_{T(\vv)} \, A_{1}\right)(\tau_{-}(\vw)) \] 
for any vector field ${\cal V}$ and ${\cal W}$. 
\end{cor} 
{\bf Proof}: 
Since $\vw - \bvw = T(\tau_{-}(\vw))$, 
\[ \gwiitwo{T(\vv) \, \, T(\vw)} - \gwiitwo{T(\vv) \, \, T(\bvw)} 
    =  \gwiitwo{T(\vv) \, \, T^{2}(\tau_{-}(\vw))}. \] 
Applying Lemma~\ref{lem:T2VW} and using the fact $T(\vv) \sim 0$, 
we obtain the desired formula. 
$\Box$ 
 
Since $\overline{T(\vw)} = 0$, if $\vw$ is replaced by $T(\vw)$ in 
Corollary~\ref{cor:TVTSW}, 
we obtain 
\[ \gwiitwo{T({\cal V}) \, \, T^{2}(\vw)} 
    = \left(\nabla_{T(\vv)} \, A_{1} \right) (\vw) \] 
for any vector field ${\cal V}$ and ${\cal W}$. 
On the other hand, since $T(\vv) \sim 0$, if $\vw$ is replaced by $T(\vw)$ in \eqref{eqn:TRR2}, 
we obtain 
\[ \gwiitwo{T({\cal V}) \, \, T^{2}(\vw)} 
    = A_{2}(\vv, T(\vw)). \] 
There are two consequences from these equations: 
\begin{cor} \label{cor:A1A2} 
Equations (\ref{eqn:TRR1}) and (\ref{eqn:TRR2}) imply 
\[ A_{2}(\vv, T(\vw)) = \left(\nabla_{T(\vv)} \, A_{1} \right)(\vw)  \] 
for any vector fields $\vv$ and $\vw$. 
\end{cor} 
{\bf Remark}: Note that both sides of the equation in Corollary~\ref{cor:A1A2} only involve 
genus-0 and genus-1 data. It might be possible to prove this relation 
by known genus-1 relations. Since our main interest in this paper is to study higher genus cases, 
 we will not do this here. 
 
\begin{cor} \label{cor:TRR1->TRR2} 
Modulo the genus-1 relation in Corollary~\ref{cor:A1A2}, 
equation~(\ref{eqn:TRR2}) follows from equation~(\ref{eqn:TRR1}) if 
one of the vector fields is equivalent to $0$. 
\end{cor} 
 
A special case of Corollary~\ref{cor:A1A2} is the following 
\begin{cor} \label{cor:A1A2str} 
\[ A_{2}(\vs, \vw) = 3 A_{1}(\tau_{-}(\vw))  \] 
for any vector field $\vw$. 
\end{cor} 
{\bf Proof}: 
This formula can be proved directly by using the derivatives of 
the string equation (\ref{eqn:DerStr}). 
On the other hand, it follows from derivatives of the string equation and the 
dilaton equation that 
\[ A_{2}(\vs, \vw) = 0 \] 
if $\vw$ is a primary vector field and 
\[ \nabla_{\vd} \, A_{1} = 3 \, A_{1}. \] 
Therefore the formula in this corollary also follows from Corollary~\ref{cor:A1A2} 
by using \eqref{eqn:WT}. 
$\Box$ 
 
In the opposite direction of Corollary~\ref{cor:TRR1->TRR2}, we have the following 
\begin{thm} \label{thm:TRR2to1} 
The topological recursion relation (\ref{eqn:TRR1}) is a formal 
consequence of (\ref{eqn:TRR2}). 
\end{thm} 
{\bf Proof}: 
The first derivatives of the genus-2 string equation and dilaton equation 
have the following form 
\begin{equation} \label{eqn:StrDil2pt} 
\gwiitwo{{\cal S} \, \vw} = \gwiitwo{\tau_{-}(\vw)}, \hspace{30pt} 
\gwiitwo{{\cal D} \, \vw} = 3 \gwiitwo{\vw} 
\end{equation} 
for any vector field $\vw$. 
Since $T({\cal S}) = {\cal D}$, applying equation~(\ref{eqn:TRR2}) 
for ${\cal V} = {\cal S}$, we obtain 
\begin{eqnarray*} 
A_{2}(\vs, \vw) & = & 
\gwiitwo{T({\cal W}) \, {\cal D}} - 3 \gwiitwo{ T(\bvw)} \\ 
&=& 3 \gwiitwo{T({\cal W})} - 3 \gwiitwo{ T(\bvw)} \\ 
&=& 3 \gwiitwo{T({\cal W}-\bvw)} \, = \, 3 \gwiitwo{T^{2}(\tau_{-}(\vw))} 
\end{eqnarray*} 
for any vector field $\cal W$. 
This is equivalent to equation~(\ref{eqn:TRR1}) 
because of the Corollary~\ref{cor:A1A2str}. 
$\Box$ 
 
Next, we study the relations between equations (\ref{eqn:TRR2}) and (\ref{eqn:BP}). 
If we replace $\vw_{3}$ by $T(\vv)$ in equation~(\ref{eqn:BP}), we obtain 
\begin{eqnarray*} 
&& - 2 \gwii{{\cal W}_{1} \, {\cal W}_{2} \, T({\cal V}) \, \gua} 
    \gwiitwo{T(\ga)} 
- \gwiitwo{ T^{2}({\cal V}) \, \{{\cal W}_{1} \bullet {\cal W}_{2}\}} 
    + \gwiitwo{T({\cal V}) \, T({\cal W}_{1} \bullet {\cal W}_{2})} 
\nonumber \\ 
 &=& B(\vw_{1}, \vw_{2}, T(\vv)). 
\end{eqnarray*} 
By Corollary~\ref{cor:4ptT}, this is equivalent to 
\begin{eqnarray*} 
 && - 2 \gwiitwo{T({\cal W}_{1} \bullet {\cal W}_{2} \bullet {\cal V})} 
- \gwiitwo{ T^{2}({\cal V}) \, \{{\cal W}_{1} \bullet {\cal W}_{2}\}} 
    + \gwiitwo{T({\cal V}) \, T({\cal W}_{1} \bullet {\cal W}_{2})} \\ 
& = & B(\vw_{1}, \vw_{2}, T(\vv)). 
\end{eqnarray*} 
On the other hand, by Theorem~\ref{thm:TRR2to1} and 
Lemma~\ref{lem:T2VW}, we know that equation~(\ref{eqn:TRR2}) implies 
\begin{eqnarray*} 
&& - 2 \gwiitwo{T({\cal W}_{1} \bullet {\cal W}_{2} \bullet {\cal V})} 
- \gwiitwo{ T^{2}({\cal V}) \, \{{\cal W}_{1} \bullet {\cal W}_{2}\}} 
    + \gwiitwo{T({\cal V}) \, T({\cal W}_{1} \bullet {\cal W}_{2})}  \\ 
&=&  - \left(\nabla_{\vw_{1} \bullet \vw_{2}} A_{1} \right) (\vv) 
    + A_{2}(\vw_{1} \bullet \vw_{2}, \vv). 
\end{eqnarray*} 
Therefore we have the following two consequences 
\begin{lem} \label{lem:AB} 
Equations (\ref{eqn:TRR2}) and (\ref{eqn:BP}) imply 
\[ B(\vw_{1}, \vw_{2}, T(\vv)) = A_{2}(\vw_{1} \bullet \vw_{2}, \vv) 
    - \left(\nabla_{\vw_{1} \bullet \vw_{2}} A_{1} \right) (\vv) \] 
for all vector fields. 
\end{lem} 
{\bf Remark}: Both sides of the equation in this lemma only involve genus-0 
and genus-1 data. It might be possible to prove this equation directly 
using known genus-1 relation. Again we will not do this in this paper. 
 
\begin{lem} \label{lem:Inverse} 
Modulo the genus-1 relation in Lemma~\ref{lem:AB}, 
equation~(\ref{eqn:BP}) 
follows from equation~(\ref{eqn:TRR2}) if 
one of the vector fields is equivalent to $0$. 
\end{lem}

In the opposite direction of Lemma~\ref{lem:Inverse}, we have 
\begin{thm} \label{thm:BPTRR1toTRR2} 
Modulo the genus-1 relation in Lemma~\ref{lem:AB}, 
the topological recursion relation (\ref{eqn:TRR2}) is a formal consequence 
of the topological recursion relation (\ref{eqn:TRR1}) and (\ref{eqn:BP}). 
\end{thm} 
{\bf Proof}: 
Since $T({\cal V}) \sim 0$, 
applying equation~(\ref{eqn:BP}) for ${\cal W}_{1} = {\cal S}$, 
${\cal W}_{2} = T({\cal V})$ and ${\cal W}_{3} = {\cal W}$, we 
obtain 
\[-2 \gwii{ {\cal S} \, T({\cal V}) \, {\cal W} \, \gua} \gwiitwo{T(\ga)} 
- \gwiitwo{ T^{2}({\cal V}) \,  \bvw } 
    + \gwiitwo{T({\cal V}) \, T(\bvw)} 
= B(\vs, \vw, T(\vv)). 
\] 
By Lemma~\ref{lem:String}, we have 
\[ \gwii{ {\cal S} \, T({\cal V}) \, {\cal W} \, \gua} \gwiitwo{T(\ga)} 
= \gwiitwo{T(\qpcc{V}{W})}. 
\] 
The theorem then follows from Lemma~\ref{lem:T2VW} and 
Corollary~\ref{cor:TVTSW} 
(which in turn follows from equation~(\ref{eqn:TRR1})), and 
Lemma~\ref{lem:stringid}. 
$\Box$

We now make a remark on higher genus topological recursion relation. 
For genus $g=1,\, 2$, the 
topological recursion relations (\ref{eqn:TRRg1}) and (\ref{eqn:TRR1}) 
are 
derived by using a formula for expressing the tautological class 
$\psi_{1}^{g}$ (i.e. the $g$-th power of 
the first Chern class of the tautological line 
bundle defined by the cotangent space of each curve at the marked point) 
on the moduli space of stable curves 
$\overline{\cal M}_{g, 1}$ in terms of boundary classes. For general $g$, 
 it was 
conjectured in \cite{G2} that polynomials of degree $g$ in the 
tautological classes $\psi_{i}$ are boundary classes on 
$\overline{\cal M}_{g, n}$. This conjecture was proved in 
\cite{Io}. A somewhat stronger version for a special case of this 
conjecture would be $\psi_{1}^{g}$ is equal to a boundary class in 
$\overline{\cal M}_{g, 1}$ without genus-$g$ components. This 
would imply the following 
\begin{equation} T^{g}(\vw) \,\,\, {\rm is \,\,\, trivial \,\,\, 
at \,\,\, the \,\,\, genus-}g \,\,\, {\rm level \,\,\, for} 
\,\,\, g \geq 1. 
\label{eqn:TRRconj} 
\end{equation} 
 
The following topological recursion relation for all $g \geq 1$ was 
derived in \cite{EX} under the assumption that the genus-$g$ 
generating function is a function of derivatives of genus-0 
generating function: 
\[ \gwiig{ \tau_{n+3g-1}(\ga)} = 
    \sum_{j=0}^{3g-2} \gwii{ \grava{n+3g-2-j} \, \, \gb} 
    \Gamma_{j}^{\beta}\] 
where $\Gamma_{0}^{\beta} = \gwiig{\gub}$ and 
\[ \Gamma_{j}^{\beta} = \gwiig{\tau_{j}(\gub)} 
    -\sum_{k=0}^{j-1} \gwii{ \tau_{k}(\gub) \, \, \gm} 
    \Gamma_{j-1-k}^{\mu}. \] 
An easy induction on $j$ shows \[ \gwiig{T^{j}(\gub)} = 
\Gamma_{j}^{\beta}\] and the similar induction also shows that this 
topological recursion relation is precisely 
\begin{equation} \label{eqn:TRREX} 
\gwiig{T^{3g-1}(\vw)} = 0 
\end{equation} 
for all vector field $\vw$. For $g=1$ and $2$, this equation 
follows from equations (\ref{eqn:derTRRg0}), 
(\ref{eqn:derTRRg1}), and (\ref{eqn:TRR1}). For general $g$, it 
follows from $\psi_{1}^{3g-1} = 0$ on $\overline{\cal M}_{g, 1}$ 
since the complex dimension of $\overline{\cal M}_{g, 1}$ is 
$3g-2$. This was first observed by Getzler.

To apply these topological recursion relations, we observe that 
\begin{equation} \label{eqn:WTW} 
\vw = T^{k} (\tau_{-}^{k}(\vw)) + \sum_{i=0}^{k-1} T^{i} 
(\overline{\tau_{-}^{i}(\vw)}) 
\end{equation} 
for any vector field $\vw$ and $k \geq 1$. This formula is obtained by 
repeatedly applying \eqref{eqn:WT}. Since 
$\overline{\tau_{-}^{i}(\vw)}$ is a primary vector field, the 
descendant level of $T^{i} (\overline{\tau_{-}^{i}(\vw)})$ is at 
most $i$. In applications, we can choose $k$ large enough to 
apply the suitable topological recursion relations. For example, 
to apply \eqref{eqn:TRRconj}, we would choose $k = g$. To apply 
\eqref{eqn:TRREX}, we would choose $k= 3g-1$.

\section{Quantum powers of the Euler vector field} 
\label{sec:EulerVF} 
 
We now turn to the application of the quantum product on the big phase 
space to the Virasoro conjecture. As pointed in \cite{LT}, the most 
important vector field in studying the Virasoro conjecture is the 
{\it Euler vector field}, which is defined by 
\begin{equation} \label{eqn:EulerDef} 
{\cal X} := - \sum_{m, \alpha} \left(m + b_{\alpha} - b_{1} -1 
                        \right)\tilde{t}^{\alpha}_{m} \, \grava{m} 
        - \sum_{m, \alpha, \beta} 
        {\cal C}_{\alpha}^{\beta}\tilde{t}^{\alpha}_{m} \, \gravb{m-1}, 
\end{equation} 
where 
\[ b_{\alpha} = ({\rm holomorphic \,\,\, dimension \,\,\, of \,\,\,} \ga) 
   - \frac{1}{2} ({\rm complex \,\,\, dimension \,\,\, of \,\,\,} V) 
   + \frac{1}{2}\] 
 and the matrix ${\cal C} = ( {\cal C}_{\alpha}^{\beta})$ 
is defined by $ c_{1}(V) \cup \ga = {\cal C}_{\alpha}^{\beta} \, 
\gb$. For compact symplectic manifolds, the holomorphic dimension 
of $\ga$ can be replaced by a half of the real dimension of $\ga$ 
in the definition of $b_{\alpha}$. Moreover, the basis $\{ 
\gamma_{1}, \ldots, \gamma_{N}\}$ of $H^{*}(V, {\Bbb C})$ can be 
chosen in a way such that the following holds: If $\eta^{\alpha 
\beta} \neq 0$ or $\eta_{\alpha \beta} \neq 0$, then $b_{\alpha} = 
1-b_{\beta}$, ${\cal C}_{\alpha}^{\beta} \neq 0$ implies 
$b_{\beta} = 1 + b_{\alpha}$, and ${\cal C}_{\alpha\beta} \neq 0$ 
implies $b_{\beta} = - b_{\alpha}$. 
 
The Euler vector field satisfies the following {\it quasi-homogeneity 
equation} 
\begin{equation} \label{eqn:quasiallgenus} 
 \left<\left< {\cal X} \right>\right>_{g} = 
        2(b_{1}+1)(1-g) F_{g} 
        + \frac{1}{2} \delta_{g, 0}\sum_{\alpha, \beta} {\cal C}_{\alpha \beta} 
        t^{\alpha}_{0}t^{\beta}_{0} 
        - \frac{1}{24} \delta_{g, 1} \int_{V} c_{1}(V) \cup c_{d-1}(V). 
\end{equation} 
The genus-0 quasi-homogeneity equation implies (cf. \cite[Lemma 
1.4 (3)]{LT}) 
\begin{eqnarray} 
&&\left<\left< \, \tau_{m}(\gamma_{\alpha}) \, {\cal X} \, 
                \tau_{n}(\gamma^{\beta}) \, \right>\right>_{0} 
    \nonumber \\ 
&=&  \delta_{m, 0} \delta_{n,0} {\cal C}_{\alpha}^{\beta} + 
        (m+n+b_{\alpha}+1-b_{\beta}) \left<\left< \, \tau_{m}(\gamma_{\alpha}) 
               \,  \tau_{n}(\gamma^{\beta}) \, \right>\right>_{0} 
    \nonumber \\ 
 &  & \textrm{ \hspace{1pt}}  + \sum_{\mu} {\cal C}_{\alpha}^{\mu} 
        \left<\left< \, \tau_{m-1}(\gamma_{\mu}) \, 
                \tau_{n}(\gamma^{\beta}) \, \right>\right>_{0} 
        + \sum_{\mu} {\cal C}^{\beta}_{\mu} 
        \left<\left< \, \tau_{m}(\gamma_{\alpha}) 
                \tau_{n-1}(\gamma^{\mu}) \, \right>\right>_{0}. 
\label{eqn:X3pt} 
\end{eqnarray} 
A special case of this formula is the following 
\[ \gwii{\ga \, {\cal X} \, \gub} = 
        {\cal C}_{\alpha}^{\beta} + (b_{\alpha} + 1 - b_{\beta}) 
        \gwii{\ga \, \gub}. \] 
Therefore for any vector field ${\cal W}$, 
\begin{equation} \label{eqn:DerEuler3pt} 
 {\cal W} \gwii{\ga \, {\cal X} \, \gub} = 
        (b_{\alpha} + 1 - b_{\beta}) 
        \gwii{\ga \, {\cal W} \, \gub}. 
\end{equation} 
An immediate application of this formula is the following 
Virasoro type relation among the quantum 
powers of the Euler vector fields: 
\begin{thm} \label{thm:VirEuler} 
For $m, k \geq 0$, 
if ${\cal W} \sim {\cal X}^{k}$ and ${\cal V} \sim {\cal X}^{m}$, then 
\[ [{\cal W}, {\cal V}] \sim (m-k) \, {\cal X}^{m+k-1}. \] 
Here (and thereafter) ${\cal X}^{0}$ is understood as ${\cal S}$ . 
\end{thm} 
{\bf Proof}: 
If ${\cal W} \sim {\cal X}^{k}$, then 
\[ \gwii{\ga \, {\cal W} \, \gub} = \gwii{\ga \, {\cal X}^{k} \, \gub} \] 
for all $\alpha, \beta$. 
Let $M$ be the $N \times N$ matrix whose $(\alpha, \beta)$ entry is 
$\gwii{\ga \, {\cal X} \, \gub}$ and $D$ be the diagonal matrix whose 
diagonal entries are $b_{1}, \ldots, b_{N}$. 
Then 
\begin{eqnarray} 
\gwii{\ga \, {\cal X}^{k} \, \gub} 
&=& \gwii{\ga \, {\cal X} \, \gamma^{\mu_{1}} } 
      \gwii{\gamma_{\mu_{1}} \, {\cal X} \, \gamma^{\mu_{2}} } 
        \ldots 
      \gwii{\gamma_{\mu_{k-1}} \, {\cal X} \, \gub } 
\, = \, \left( M^{k} \right)_{\alpha}^{\beta}. 
        \label{eqn:EulerPower} 
\end{eqnarray} 
Therefore, by equation~(\ref{eqn:DerEuler3pt}), we have 
\[ {\cal W} \, M = D M^{k} - M^{k}D + M^{k}, \] 
and consequently for any integer $m \geq 0$, 
\begin{eqnarray*} 
 {\cal W} \, M^{m} 
& = &  \sum_{j=1}^{m} M^{j-1} (D M^{k} - M^{k}D + M^{k}) M^{m-j} \\ 
& = & m M^{m+k-1} + \sum_{j=1}^{m} M^{j-1} D M^{m+k-j} 
      - \sum_{j=1}^{m} M^{k+j-1}D M^{m-j}. 
\end{eqnarray*} 
Hence, 
if ${\cal W} \sim {\cal X}^{k}$ and ${\cal V} \sim {\cal X}^{m}$, then 
\begin{equation} \label{eqn:VirEulerMat} 
 {\cal W} \, M^{m} - {\cal V} \, M^{k} 
    = (m - k) M^{m+k-1}. 
\end{equation} 
Since every vector on the big phase space is equivalent to a 
primary vector field, to prove the theorem, it suffices to show 
\[ \ga \bullet [{\cal W}, {\cal V}] = (m-k) \, \ga \bullet {\cal X}^{m+k-1} \] 
for all $\alpha$. 
Since $\ga$ is a constant vector field, by equation~(\ref{eqn:DerProd}), 
we have 
\begin{eqnarray*} 
 \ga \bullet [{\cal W}, {\cal V}] 
&=&  \ga \bullet \nabla_{\cal W} {\cal V} - 
    \ga \bullet \nabla_{\cal V} {\cal W} \\ 
&=&  \nabla_{\cal W} (\ga \bullet {\cal V}) - 
    \nabla_{\cal V} (\ga \bullet {\cal W}) \\ 
&=&  \nabla_{\cal W} (\ga \bullet {\cal X}^{m}) - 
    \nabla_{\cal V} (\ga \bullet {\cal X}^{k}). 
\end{eqnarray*} 
Since for any $m \geq 0$, 
\[ \ga \bullet {\cal X}^{m} = (M^{m})_{\alpha}^{\beta} \, \gb \] 
where $(M^{m})_{\alpha}^{\beta}$ is the $(\alpha, \beta)$ entry 
of $M^{m}$, the theorem follows from equation~(\ref{eqn:VirEulerMat}). 
$\Box$ 
 
Since $\bvx^{k} \sim {\cal X}^{k}$ for all $k$ 
and Lie brackets of primary 
fields are always primary fields, we have 
\begin{cor} 
\[ [ \bvx^{k}, \, \bvx^{m}] = (m-k) \, \bvx^{m+k-1} \] 
for all $m, k \geq 0$. 
\end{cor} 
We also have 
\begin{cor} \label{cor:VirEuler} 
For $m, k \geq 0$ and $\{m, k\} \neq \{0, 2\}$, 
\[ [ {\cal X}^{k}, \, {\cal X}^{m}] = (m-k) \, {\cal X}^{m+k-1}. \] 
Note that 
$[{\cal S}, {\cal X}^{2}] \neq 2{\cal X}$ since 
$[{\cal S}, {\cal X}^{2}]$ is a primary vector field but 
${\cal X}$ contains descendant vector fields. In this case we have 
\[ [{\cal S}, {\cal X}^{2}] = 2 \, \bvx.\] 
\end{cor} 
{\bf Proof}: 
A straightforward computation shows that 
\[ [{\cal S}, {\cal X}] = {\cal S}.\] 
If $\{m, k\} \neq \{0, 1\}$, then 
 $[ {\cal X}^{k}, \, {\cal X}^{m}]$ is a primary vector field. So 
\[[ {\cal X}^{k}, \, {\cal X}^{m}] 
    = {\cal S} \bullet [ {\cal X}^{k}, \, {\cal X}^{m}] 
    = (m-k) \, \bvx^{m+k-1} \] 
by Theorem~\ref{thm:VirEuler}. 
In particular, 
\[ [{\cal S}, {\cal X}^{2}] = 2 \, \bvx.\] 
If in addition $\{m, k\} \neq \{0, 2\}$, then 
$(m-k) \, {\cal X}^{m+k-1}$ is also a primary vector field. 
So 
\[(m-k) \, \bvx^{m+k-1} 
    = (m-k) \, {\cal X}^{m+k-1}. \] 
This proves the corollary. 
$\Box$ 
 
Another application of the quasi-homogeneity equation is the following 
\begin{lem} \label{lem:TWEulerPower} 
For any vector field ${\cal W}$ and integer $k \geq 0$, 
\[ \nabla_{T({\cal W})} \bvx^{k} 
    = - {\cal W} \bullet \bvx^{k}. \] 
\end{lem} 
{\bf Proof}: 
By formula~(\ref{eqn:DerEuler3pt}) and 
(\ref{eqn:EulerPower}), 
\[ T({\cal W}) \gwii{\ga \, {\cal X}^{k} \, \gub} = 0 \] 
since $T({\cal W}) \sim 0$. 
This formula is also true for $k=0$ due to the string equation. 
Therefore, by Lemma~\ref{lem:String}, 
\begin{eqnarray*} 
\nabla_{T({\cal W})} \bvx^{k} 
& = & \nabla_{T({\cal W})} \left\{ \gwii{{\cal S} \, {\cal S} \, \gua} 
    \gwii{\ga \, {\cal X}^{k} \, \gub} \gb \right\} \\ 
& = &  \left\{ T({\cal W}) \gwii{{\cal S} \, {\cal S} \, \gua} \right\} 
    \gwii{\ga \, {\cal X}^{k} \, \gub} \gb  \\ 
& = &  - \gwii{{\cal W} \, {\cal S} \, \gua} 
    \gwii{\ga \, {\cal X}^{k} \, \gub} \gb  \\ 
& = & - \vw \bullet \bvx^{k}. 
\end{eqnarray*} 
$\Box$ 
\begin{cor} 
\label{cor:BracketTXpower} 
\[ [ T(\bvx^{m}), \, T(\bvx^{k})] = 0 \] 
for all $m, k \geq 0$. 
\end{cor} 
{\bf Proof}: 
By Lemma~\ref{lem:DerT}, 
\begin{eqnarray*} 
 [ T(\bvx^{m}), \, T(\bvx^{k})] 
& = & \nabla_{T(\bvx^{m})} T(\bvx^{k}) - 
    \nabla_{T(\bvx^{k})} T(\bvx^{m})  \\ 
& = & T(\nabla_{T(\bvx^{m})} \bvx^{k}) - 
    T(\nabla_{T(\bvx^{k})} \bvx^{m})  \\ 
& = & T(-\bvx^{m+k}) - 
    T(-\bvx^{k+m}) = 0. 
\end{eqnarray*} 
$\Box$ 
 
\begin{cor} \label{cor:TWderTau-S} 
For any vector field ${\cal W}$ and integer $k \geq 0$, 
\[ \nabla_{T({\cal W})} (\tau_{-}({\cal S}) \bullet \bvx^{k}) 
    = - \tau_{-}({\cal W}) \bullet \bvx^{k}. \] 
\end{cor} 
{\bf Proof}: By Lemma~\ref{lem:DerT} and \ref{lem:String}, 
\[ \nabla_{T({\cal W})} (\tau_{-}({\cal S})) 
   = \tau_{-}(\nabla_{T({\cal W})} {\cal S}) 
   =  - \tau_{-}({\cal W}). \] 
Therefore by equation~(\ref{eqn:DerProd}), 
Corollary~\ref{cor:4ptT} 
 and Theorem~\ref{lem:TWEulerPower}, 
\begin{eqnarray*} 
\nabla_{T({\cal W})} (\tau_{-}({\cal S}) \bullet 
    \bvx^{k}) 
& = & (\nabla_{T({\cal W})} (\tau_{-}({\cal S}))) \bullet 
    \bvx^{k} 
    +\tau_{-}({\cal S}) \bullet \nabla_{T({\cal W})} 
    \bvx^{k} \\ 
&&  + \gwii{T({\cal W}) \, \, \tau_{-}({\cal S}) \, \, 
    \bvx^{k} \, \gua} \ga \\ 
& = & (-\tau_{-}({\cal W})) \bullet 
    \bvx^{k} 
    +\tau_{-}({\cal S}) \bullet (-{\cal W} 
    \bullet \bvx^{k}) \\ 
&&  + {\cal W} \bullet \tau_{-}({\cal S}) \bullet 
    \bvx^{k} \\ 
& = & -\tau_{-}({\cal W}) \bullet \bvx^{k} 
\end{eqnarray*} 
$\Box$ 
 
So far we have only used very special cases of 
equation~(\ref{eqn:X3pt}). To take full advantage of this 
equation, we need to introduce two operations on the space of 
vector fields. 
\begin{defi} 
For ${\cal W} = \sum_{m, \alpha} f_{m, \alpha} \grava{m}$ 
and ${\cal V} = \sum_{m, \alpha} g_{m, \alpha} \grava{m}$, define 
\[ {\cal W} * {\cal V} = \sum_{m, \alpha} f_{m, \alpha} 
        g_{m, \alpha} \grava{m} \] 
and 
\[ C({\cal W}) = \sum_{m, \alpha, \beta} f_{m, \alpha} 
        {\cal C}_{\alpha}^{\beta} \gravb{m}. \] 
\end{defi} 
We can think of ``$*$'' as a commutative and associative product 
on the space of vector fields with the identity 
\begin{equation} 
\label{eqn:starID} 
 {\cal Z} := \sum_{m, \alpha} \grava{m}. 
\end{equation} 
Another important vector field for this product is 
\begin{equation} 
\label{eqn:gradingVF} 
 {\cal G} := \sum_{m, \alpha} (m + b_{\alpha}) \grava{m}. 
\end{equation} 
Using this vector field, we can define another linear transformation 
on the space of vector fields: 
\begin{defi} \label{def:R} 
For any vector field ${\cal W}$, define 
\[ R({\cal W}) = {\cal G} * T({\cal W}) + C({\cal W}). \] 
\end{defi} 
We then have the following 
\begin{thm} \label{thm:Xprod} 
For any vector field ${\cal W}$, 
\[ \qpcc{W}{X} = \qpc{S}{R({\cal W})} = \overline{R(\vw)}. \] 
\end{thm} 
{\bf Proof}: 
Let ${\cal W} = \sum_{m, \alpha} f_{m, \alpha} \grava{m}$. 
Then 
\begin{eqnarray*} 
\qpcc{W}{X} 
& = & \sum_{m, \alpha, \mu} f_{m, \alpha} 
    \gwii{\grava{m} \, {\cal X} \, \gum} \gm. 
\end{eqnarray*} 
Applying equation~(\ref{eqn:X3pt}) to each 3-point functions on the 
right hand side, we have 
\begin{eqnarray*} 
\qpcc{W}{X} 
& = & \gwii{\{{\cal G} * {\cal W} \} \, \gum} \gm 
    + (1-b_{\mu}) \gwii{{\cal W} \, \gum} \gm \\ 
& & 
    + \gwii{\{ \tau_{-}(C({\cal W})) \} \, \gum} \gm 
    + \sum_{\alpha} f_{0, \alpha} {\cal C}_{\alpha}^{\mu} \gm. 
\end{eqnarray*} 
By equation~(\ref{eqn:StringProd}), we have 
\begin{eqnarray*} 
\qpcc{W}{X} 
& = & {\cal S} \bullet \{\tau_{+}({\cal G} * {\cal W}) \} 
    + (1-b_{\mu}) \gwii{{\cal S} \, \tau_{+}({\cal W}) \, \gum} \gm 
    + {\cal S} \bullet C({\cal W}). 
\end{eqnarray*} 
On the other hand, it is straightforward to check that 
\[ \tau_{+}({\cal G} * {\cal W}) 
    = {\cal G} * \tau_{+}({\cal W}) - \tau_{+}({\cal W}) \] 
for any vector field ${\cal W}$. Therefore 
\begin{eqnarray*} 
\qpcc{W}{X} 
& = & {\cal S} \bullet \{{\cal G} * \tau_{+}({\cal W}) \} 
    -b_{\mu} \gwii{{\cal S} \, \tau_{+}({\cal W}) \, \gum} \gm 
    + {\cal S} \bullet C({\cal W}) \\ 
& = & {\cal S} \bullet \{{\cal G} * T({\cal W}) + 
        {\cal G} * ({\cal S} \bullet \tau_{+}({\cal W}))\} 
    - {\cal G}* ({\cal S} \bullet \tau_{+}({\cal W})) 
    + {\cal S} \bullet C({\cal W}). 
\end{eqnarray*} 
The middle terms in the last expression are cancelled since ${\cal 
G}* ({\cal S} \bullet \tau_{+}({\cal W}))$ is a primary vector 
field. This proves the theorem. $\Box$

\begin{cor} \label{cor:Requiv} 
If ${\cal W} \sim {\cal V}$, then $R({\cal W}) \sim R({\cal V})$. 
\end{cor} 
{\bf Proof}: 
$R({\cal W}) \sim R({\cal V})$ if and only if 
${\cal S} \bullet R({\cal W}) = {\cal S} \bullet R({\cal V})$, 
which indeed follows from Theorem~\ref{thm:Xprod}. 
$\Box$ 
 
 The following formula will be useful later. 
\begin{lem} \label{lem:Rtau-} 
For any vector field ${\cal W}$, 
\[ R(\tau_{-}({\cal W})) = \tau_{-}(R({\cal W})) - 
        {\cal G} * \bvw - {\cal W}.\] 
\end{lem} 
{\bf Proof}: 
It is straightforward to check that 
\[ \tau_{-}({\cal G}*{\cal W}) = {\cal G}* \tau_{-}({\cal W}) 
            + \tau_{-}({\cal W}) \] 
and 
\[ \tau_{-}(C({\cal W})) = C(\tau_{-}({\cal W})) \] 
for any vector field ${\cal W}$. 
Therefore 
\begin{eqnarray*} 
\tau_{-}(R({\cal W}))& = &{\cal G}* {\cal W} + {\cal W} 
            + C(\tau_{-}({\cal W}) ) \\ 
& = &{\cal G}* (T(\tau_{-}({\cal W})) + \bvw) + {\cal W} 
            + C(\tau_{-}({\cal W}) ) \\ 
& = &R(\tau_{-}({\cal W})) + {\cal G}* \bvw + {\cal W}. 
\end{eqnarray*} 
The lemma follows. 
$\Box$ 
 
Since $T({\cal S}) = \vd$, it is straightforward to check that 
\begin{equation} \label{eqn:RS} 
    R({\cal S}) = {\cal X} + (b_{1} + 1) {\cal D} 
\end{equation} 
So by Theorem~\ref{thm:Xprod} and Lemma~\ref{lem:Rtau-}, 
\begin{eqnarray*} 
 {\cal X} \bullet \tau_{-}({\cal S}) 
    & = & {\cal S} \bullet R(\tau_{-}({\cal S})) \\ 
    & = &  {\cal S} \bullet (\tau_{-}(R({\cal S})) 
        - {\cal G} * {\cal S}^{2} - {\cal S}) \\ 
    & = &  {\cal S} \bullet (\tau_{-}({\cal X} + (b_{1} + 1) {\cal D}) 
        - {\cal G} * {\cal S}^{2} - {\cal S}) \\ 
    & = & {\cal S} \bullet \tau_{-}({\cal X}) + b_{1} {\cal S}^{2} 
        - {\cal G} * {\cal S}^{2}. 
\end{eqnarray*} 
Equivalently, we have 
\begin{equation} \label{eqn:tau-X} 
 {\cal S} \bullet \tau_{-}({\cal X}) 
    =  {\cal X} \bullet \tau_{-}({\cal S}) - b_{1} {\cal S}^{2} 
        + {\cal G} * {\cal S}^{2}. 
\end{equation} 
 
Another useful formula for the operator $R$ is the following 
\begin{lem} 
\label{lem:derR} 
For any vector fields ${\cal W}$ and ${\cal V}$, 
\[ \nabla_{\cal V} (R({\cal W})) = R(\nabla_{\cal V} {\cal W}) 
    - {\cal G} * (\qpcc{V}{W}). \] 
\end{lem} 
{\bf Proof}: 
It is straightforward to check that 
\begin{equation} \label{eqn:derGC} 
 \nabla_{\cal V} ({\cal G} * {\cal W}) = 
        {\cal G} * (\nabla_{\cal V} {\cal W}) 
  \hspace{20pt} {\rm and} \hspace{20pt} 
 \nabla_{\cal V} (C({\cal W})) = C(\nabla_{\cal V} {\cal W}) 
\end{equation} 
for any vector fields ${\cal W}$ and ${\cal V}$. 
Therefore 
\begin{eqnarray*} 
\nabla_{\cal V} (R({\cal W})) 
& = & \nabla_{\cal V} ({\cal G}* T({\cal W}) + C({\cal W})) \\ 
& = & {\cal G} * \nabla_{\cal V} T({\cal W}) 
    + C( \nabla_{\cal V} {\cal W}) \\ 
& = & {\cal G} * (T(\nabla_{\cal V} {\cal W}) - \qpcc{V}{W}) 
    + C( \nabla_{\cal V} {\cal W}) \\ 
& = & R(\nabla_{\cal V} {\cal W}) - {\cal G} * (\qpcc{V}{W}). 
\end{eqnarray*} 
$\Box$ 
 
The quasi-homogeneity equation for all genera 
 implies the following rule for 
eliminating the Euler vector field from a correlation function: 
\begin{lem} \label{lem:rmX} 
\begin{eqnarray*} 
\gwiig{\vx \, \vw_{1} \, \cdots \, \vw_{k}} 
& = & \sum_{i=1}^{k} 
    \gwiig{\vw_{1} \, \cdots \, \left\{ \tau_{-} R (\vw_{i}) \right\} 
    \, \cdots \, \vw_{k}}  \\ 
&& - \left\{ (2g+k-2) b_{1} + 2 (g+k-1) \right\} 
        \gwiig{\vw_{1} \, \cdots \, \vw_{k}} \\ 
&& + \delta_{g, 0} \nabla^{k}_{\vw_{1}, \cdots, \vw_{k}} \left( 
    \frac{1}{2} {\cal C}_{\alpha \beta} t_{0}^{\alpha} t_{0}^{\beta} 
    \right) 
\end{eqnarray*} 
for any vector fields $\vw_{1}, \ldots, \vw_{k}$ and $k \geq 1$. 
Here $\nabla^{k}_{\vw_{1}, \cdots, \vw_{k}}$ is the $k$-th covariant 
derivative. 
Note that if $g>0$ or $k > 2$, the last term on the right hand side 
vanishes. 
\end{lem} 
{\bf Proof}: It is straightforward to check that 
\begin{eqnarray*} 
 \nabla_{\vw} \vx & = & - \vg * \vw + (b_{1} + 1) \vw - C(\tau_{-}(\vw)) \\ 
&=& - R (\tau_{-} (\vw)) - \vg* \bvw + (b_{1}+1) \vw 
\end{eqnarray*} 
for any vector field $\vw$. By Lemma~\ref{lem:Rtau-}, we also have 
\begin{equation} \label{eqn:derX} 
\nabla_{\vw} \vx = - \tau_{-} R(\vw) + (b_{1} + 2) \vw. 
\end{equation} 
It suffices to prove the lemma for parallel vector fields $\vw_{i}$ 
since both sides of the equation are tensors in these vector fields. 
We observed that if $\vw$ is parallel, i.e. $\nabla_{\vv} \vw =0$ 
for all $\vv$, then $\tau_{-}R(\vw)$ is also parallel. 
In fact, by Lemma~\ref{lem:DerT} and Lemma~\ref{lem:derR} 
\[ \nabla_{\vv} \, \tau_{-}R(\vw) = \tau_{-} R (\nabla_{\vv} \vw) 
    - \tau_{-}(\vg * (\vv \bullet \vw)) = 0\] 
since $\vw$ is parallel and $\vg * (\vv \bullet \vw)$ is a primary 
vector field. 
The lemma then follows from repeatedly taking derivatives of 
\eqref{eqn:quasiallgenus} 
and using \eqref{eqn:DerCorr}. 
$\Box$

In case that $\vw$ is a primary vector field, $\tau_{-}(\vw) = 0$. 
Therefore by Lemma~\ref{lem:Rtau-}, 
\[ \tau_{-}R(\vw) = \vg * \vw + \vw. \] 
Hence the formula in Lemma~\ref{lem:rmX} can be written as 
\begin{eqnarray} 
 \gwiig{{\cal X} \, \vw_{1} \, \ldots \, \vw_{k}} 
&=&   \sum_{i=1}^{k}  \gwiig{ \vw_{1} \, \ldots \, 
        \{ \vg * \vw_{i}\} \, \ldots   \, \vw_{k}} 
\nonumber \\ 
&&     - \left\{k-2 + 2g + (2g+k-2) b_{1} \right\} 
   \gwiig{\vw_{1} \, \ldots \, \vw_{k}} 
\nonumber \\ 
&& + \delta_{g, 0} \nabla^{k}_{\vw_{1}, \cdots, \vw_{k}} \left( 
    \frac{1}{2} {\cal C}_{\alpha \beta} t_{0}^{\alpha} t_{0}^{\beta} 
    \right) 
\label{eqn:rmXprim} 
\end{eqnarray} 
for all primary fields $\vw_{1}, \ldots, \vw_{k}$ and $k \geq 2$. 
 
An immediate consequence of Lemma~\ref{lem:rmX} is the following 
\begin{lem} \label{lem:XXX4pt} 
\begin{eqnarray*} 
 \gwii{{\cal X} \, {\cal X}  \, {\cal X} \, \gua} \ga 
& = & 2 \vx^{3} \bullet \tau_{-}(\vs) - 3b_{1} \vx^{2} \\ 
&&  + 2 \vx^{2} \bullet (\vg * \bvs) 
    + 2 \vx \bullet (\vg * \bvx) 
    - \vg * \vx^{2}, \\ 
 \gwii{{\cal X} \, {\cal X}  \, \bvx^{k} \, \gua} \ga 
& = & \vx^{k+2} \bullet \tau_{-}(\vs) - 2b_{1} \bvx^{k+1} 
    + \vx^{k+1} \bullet (\vg * \bvs)  \\ 
&&  + \vx^{k} \bullet (\vg * \bvx) 
    + \vx \bullet (\vg * \bvx^{k}) 
    - \vg * \vx^{2}, \\ 
 \gwii{{\cal X} \, \bvx^{m} \, \bvx^{k}  \, \gua} \ga 
& = &  - b_{1} \bvx^{m+k} - \vg * \bvx^{m+k} \\ 
&&  +  \vx^{k} \bullet (\vg * \bvx^{m}) 
    + \vx^{m} \bullet (\vg * \bvx^{k}), 
\end{eqnarray*} 
where $m$, $k \geq 0$. 
\end{lem} 
{\bf Proof}: The last formula is a direct consequence of 
\eqref{eqn:rmXprim}. Since by \eqref{eqn:WT}, 
\[ \vx = T(\tau_{-}(\vx)) + \bvx, \] 
corollary~\ref{cor:4ptT} implies that 
\begin{eqnarray*} 
&&  \gwii{{\cal X} \, {\cal X}  \, \bvx^{k} \, \gua} \ga \\ 
&=&  \gwii{{\cal X} \, T(\tau_{-}({\cal X}))  \, 
    (\bvx^{k}) \, \gua} \ga 
    +  \gwii{{\cal X} \, \bvx  \, 
            \bvx^{k} \, \gua} \ga \\ 
&=& \vx^{k+1} \bullet \tau_{-}(\vx) 
    +  \gwii{{\cal X} \, \bvx  \, \bvx^{k} \, \gua} \ga 
\end{eqnarray*} 
and 
\begin{eqnarray*} 
&&  \gwii{{\cal X} \, {\cal X}  \, \vx \, \gua} \ga \\ 
&=&  \gwii{{\cal X} \, T(\tau_{-}({\cal X}))  \,  \vx \, \gua} \ga 
    +  \gwii{{\cal X} \, \bvx  \, \vx \, \gua} \ga \\ 
&=& \vx^{2} \bullet \tau_{-}(\vx) 
    +  \gwii{{\cal X} \, \bvx  \, \vx \, \gua} \ga \\ 
&=& 2 \vx^{2} \bullet \tau_{-}(\vx) 
    +  \gwii{{\cal X} \, \bvx  \, \bvx \, \gua} \ga 
\end{eqnarray*} 
Therefore the first two equations in the lemma follows from the 
third equation and \eqref{eqn:tau-X}. 
$\Box$ 
 
If a 4-point function does not contain $\vx$, but contains 
$\vx^{k}$ with $k \geq 2$, it can be computed by using the 
following lemma to create 4-point functions involving $\vx$. 
\begin{lem} \label{lem:Wk4pt} 
For any vector field $\vw$, $\vv_{1}$, $\vv_{2}$, $\vv_{3}$, 
 and $k \geq 2$, 
\begin{eqnarray*} 
 \gwii{ \vw^{k} \, \vv_{1} \, \vv_{2} \, \vv_{3}} 
& = & - \sum_{i=1}^{k-1} 
    \gwii{\vw \, \,  \vw^{k-i} \, \, 
    \{\vv_{1} \bullet \vv_{2} \bullet \vw^{i-1}\} 
        \, \,  \vv_{3}} \\ 
&&  + \sum_{i=2}^{k-1} 
    \gwii{\vw \, \, \{\vw^{k-i} \bullet \vv_{1}\} \, \, 
        \{\vv_{2} \bullet \vw^{i-1}\} 
        \,\,  \vv_{3}}\\ 
&&  + \gwii{\vw \, \, \{\vw^{k-1} \bullet \vv_{1}\} \, \, 
        \vv_{2} \,\,  \vv_{3}} 
          + \gwii{\vw \, \, \vv_{1} \, \, 
        \{\vv_{2} \bullet \vw^{k-1}\} 
        \,\,  \vv_{3}}. 
\end{eqnarray*} 
\end{lem} 
This lemma follows from the following form of the 
 first derivatives of the generalized 
WDVV equation 
\begin{eqnarray} 
\gwii{ \{\vw_{1} \bullet \vw_{2} \} \vw_{3} \vw_{4} \vw_{5}} 
& = & \gwii{ \{\vw_{1} \bullet \vw_{3} \} \vw_{2} \vw_{4} \vw_{5}} 
\nonumber \\ 
&&+ \gwii{\vw_{1} \vw_{3} \{\vw_{2} \bullet \vw_{4}\} \vw_{5}} 
\nonumber \\ 
&& -\gwii{ \vw_{1} \vw_{2} \{\vw_{3} \bullet \vw_{4}\} \vw_{5}} 
\label{eqn:derWDVV} 
\end{eqnarray} 
for any vector fields $\vw_{1}, \ldots, \vw_{5}$. 
A similar lemma on the small phase space was 
proved in \cite{L1}. The same proof also works here on the 
big phase space. 
\begin{cor} \label{cor:Xk4pt} 
For all primary vector fields $\vw$, $\vv$, 
 and $k \geq 0$, 
\begin{eqnarray*} 
  \gwii{ \bvx^{k} \, \vw \, \vv \, \gua} \ga 
& = & - \bvx^{k} \bullet \vw \bullet \vv \bullet \tau_{-}(\vs)  \\ 
&&  \hspace{-40pt} + \sum_{i=1}^{k} \left\{ 
        (\vg * (\bvx^{k-i} \bullet \vw)) \bullet 
        \vv \bullet \bvx^{i-1} 
    + \bvx^{k-i} \bullet \vw \bullet 
        (\vg * (\vv \bullet \bvx^{i-1})) 
     \right. \\ 
&&  \hspace{-20pt} \left. - (\vg * \bvx^{k-i}) \bullet 
    \vw \bullet \vv \bullet \bvx^{i-1} 
    - \bvx^{k-i} \bullet 
    (\vg * (\vw \bullet \vv \bullet \bvx^{i-1})) 
             \right\} 
\end{eqnarray*} 
\end{cor} 
{\bf Proof}: Since $ \bvs = \vs - T(\tau_{-}(\vs)),$ by 
Corollary~\ref{cor:4ptT} and Lemma~\ref{lem:String}, 
\begin{eqnarray*} 
&& \gwii{\bvs \, \vw \, \vv \, \gua} \ga \\ 
&=&\gwii{\vs \, \vw  \, \vv   \, \gua} \ga 
    - \tau_{-}(\vs) \bullet \vw \bullet \vv  \\ 
&=&- \tau_{-}(\vs) \bullet \vw \bullet \vv 
\end{eqnarray*} 
if $\vw$ and $\vv$ are primary fields. This proves the corollary 
for $k=0$. Similarly since $\bvx = \vx - T(\tau_{-}(\vx))$, by 
Corollary~\ref{cor:4ptT}, 
\begin{eqnarray*} 
&& \gwii{\bvx \,\vw \, \vv  \, \gua} \ga \\ 
&=&\gwii{\vx \, \vw \, 
    \vv  \, \gua} \ga 
    - \tau_{-}(\vx) \bullet \vw \bullet \vv. 
\end{eqnarray*} 
By \eqref{eqn:rmXprim} and \eqref{eqn:tau-X}, 
\begin{eqnarray*} 
&& \gwii{\bvx \, \vw \, \vv \, \gua} \ga \\ 
& = &  -  \bvx \bullet \vw \bullet \vv \bullet \tau_{-}(\vs) 
    - \vw \bullet \vv \bullet (\vg * \bvs) 
    - \vg * (\vw \bullet \vv) \\ 
&&  + \vw \bullet (\vg * \vv) 
    + \vv \bullet (\vg * \vw). 
\end{eqnarray*} 
This proves the corollary for $k=1$. For $k \geq 2$, the corollary 
follows from Lemma~\ref{lem:Wk4pt} and \eqref{eqn:rmXprim}. 
 $\Box$ 
 
In particular, we have 
\begin{cor} \label{cor:Epower4pt} 
\begin{eqnarray*} 
&& \gwii{\bvx^{n} \, \bvx^{m} \, \bvx^{k}  \, \gua} \ga \\ 
& = &  -  \vx^{n+m+k} \bullet \tau_{-}(\vs) \\ 
&&  - \sum_{i=0}^{k-1} \vx^{i} \bullet 
        (\vg * \bvx^{n+m+k-i-1}) 
    - \sum_{i=n+m}^{n+m+k-1} \vx^{i} \bullet 
        (\vg * \bvx^{n+m+k-i-1}) \\ 
&&  + \sum_{i=m}^{m+k-1} \vx^{i} \bullet 
        (\vg * \bvx^{n+m+k-i-1}) 
    + \sum_{i=n}^{n+k-1} \vx^{i} \bullet 
        (\vg * \bvx^{n+m+k-i-1}) 
\end{eqnarray*} 
for $n, m, k \geq 0$. 
\end{cor}

\section{Virasoro vector fields} 
\label{sec:VirVF} 
 
In this section, we study a sequence of vector fields on the big 
phase space obtained from the string vector field by recursively 
applying the operator $R$ defined in Definition~\ref{def:R}. We 
will show that this sequence of vector fields satisfy the 
Virasoro bracket relation and arise naturally in the study of 
higher genus Virasoro conjecture. For this reason, we call these 
vector fields Virasoro vector fields. 
 
Define the $k$-th {\it Virasoro vector field} by 
\begin{equation} \label{eqn:VirVF} 
\vl_{k} := R^{k+1}(-\vs) 
\end{equation} 
for $k \geq -1$. Here $\vl_{-1}$ should be understood as $- \vs$. 
By Theorem~\ref{thm:Xprod}, these vector fields satisfy 
the following interesting property 
\begin{lem} \label{lem:Virg0} 
\[ {\cal L}_{n} \sim - {\cal X}^{n+1} \] 
for all $n \geq -1$. 
\end{lem} 
Moreover 
Lemma~\ref{lem:derR} and Lemma~\ref{lem:String} make 
it simple to compute the covariant 
derivatives of these vector fields. In fact, we have the following 
\begin{lem} \label{lem:derLk} 
\[ \nabla_{\vw} \vl_{k} = \tau_{-}R^{k+1}(\vw) - (k+1) R^{k}(\vw) \] 
for any vector field $\vw$. 
\end{lem} 
{\bf Proof}: By Lemma~\ref{lem:String} 
\[ \nabla_{\vw} \vl_{-1} = - \nabla_{\vw} \vs = \tau_{-}(\vw). \] 
This proves the lemma for $k=-1$. We now prove the lemma by induction 
on $k$. Suppose that the lemma holds for $k=n$. Then 
by Lemma~\ref{lem:derR} and Lemma~\ref{lem:Rtau-} 
\begin{eqnarray*} 
 \nabla_{\vw} \vl_{n+1} 
&= & \nabla_{\vw} R(\vl_{n}) 
\, = \, R(\nabla_{\vw} \vl_{n}) - \vg * (\vw \bullet \vl_{n}) \\ 
& = & R\left\{\tau_{-}R^{n+1}(\vw) - (n+1) R^{n}(\vw)\right\} 
    - \vg * (\vw \bullet \vl_{n}) \\ 
&=& \tau_{-}R^{n+2}(\vw) - \vg * \overline{R^{n+1}(\vw)} - (n+2)R^{n+1}(\vw) 
    - \vg * (\vw \bullet \vl_{n}). 
\end{eqnarray*} 
The lemma follows since by Theorem~\ref{thm:Xprod} 
\[ \overline{R^{n+1}(\vw)} = \vx^{n+1} \bullet \vw 
    = - \vw \bullet \vl_{n}.\] 
$\Box$ 
 
Since covariant derivatives commute with the operator $\tau_{-}$, we have 
\begin{cor} \label{cor:dertau-mLk} 
\[ \nabla_{\vw} \, \, \tau_{-}^{m}(\vl_{k}) 
     = \tau_{-}^{m+1} R^{k+1}(\vw) - (k+1) \tau_{-}^{m}R^{k}(\vw) \] 
for any vector field $\vw$. 
\end{cor} 
A special case of this corollary is the following 
\[ \nabla_{\tau_{-}^{n}(\vl_{j})} \, \, \tau_{-}^{m}(\vl_{k}) 
     = \tau_{-}^{m+1} R^{k+1} \tau_{-}^{n}(\vl_{j}) 
     - (k+1) \tau_{-}^{m}R^{k}\tau_{-}^{n}(\vl_{j}). \] 
We can use this formula to compute the brackets of vector fields 
$\tau_{-}^{m} \vl_{k}$. For example, when $m=n=0$, we have 
\begin{equation} 
\label{eqn:DerLk} 
 \nabla_{\vl_{j}} \, \, \vl_{k} 
     = \tau_{-} (\vl_{k+j+1}) - (k+1) \vl_{k+j}. 
\end{equation} 
Since 
$ [\vl_{j}, \, \vl_{k}] = \nabla_{\vl_{j}}  \, \vl_{k} 
        - \nabla_{\vl_{k}} \, \vl_{j}$, 
we have 
\begin{cor} \label{cor:BracketVir} 
\[  [\vl_{j}, \, \vl_{k}] = (j-k) \vl_{k+j} \] 
for $k, j \geq -1$. 
\end{cor} 
In other words, the sequence of vector fields 
$\{\vl_{k} = R^{k+1}(-\vs) \mid k \geq -1\}$ 
form a half branch of the Virasoro algebra. 
More generally, in view of Lemma~\ref{lem:Rtau-}, the more general 
class of vector fields 
$\{ \tau_{-}^{m}(\vl_{k}) \mid m \geq 0, \,\, k \geq -1 \}$ seems generate a 
$W$-type algebra. 
 
In some cases, it is more convenient to use the following formula 
for covariant derivatives of the Virasoro vector fields. 
\begin{lem} \label{lem:derLkv2} 
\[ \nabla_{\vw} \vl_{k} = R^{k+1}\tau_{-}(\vw) 
    + \sum_{i=0}^{k} R^{i}(\vg*(\vx^{k-i} \bullet \vw)) \] 
for any vector field $\vw$. 
\end{lem} 
{\bf Proof}: This follows from Lemma~\ref{lem:derLk} by interchanging 
positions of $\tau_{-}$ and $R^{k+1}$ using Lemma~\ref{lem:Rtau-}. 
$\Box$ 
 
Since $T(\vw) \, \bullet \, \vv = 0$ for any vector field $\vv$, 
one application of Lemma~\ref{lem:derLkv2} is the following 
\begin{equation} \label{eqn:DerTLk} 
 \nabla_{T(\vw)} \vl_{k} = R^{k+1}(\vw) 
\end{equation} 
for any vector field $\vw$. In particular, 
\[ \nabla_{T(\vl_{m})} T(\vl_{k}) = T(\nabla_{T(\vl_{m})} \vl_{k}) 
     = T(\vl_{m+k+1}). \] 
Therefore 
\[ [T(\vl_{m}), \, T(\vl_{k})] = 0\] 
for all $m, k \geq -1$. 
Moreover, the two sequences of vector fields 
$\{ \vl_{k} \mid k \geq -1\}$ 
and 
$\{ T(\vl_{k}) \mid k \geq -1\}$ 
together form a Lie algebra isomorphic to the Lie 
algebra spanned by $\{t^{m} \partial_{t}, \, \, \, t^{k} \mid m, k \geq 0\}$ 
as operators on the space of functions on the unit circle, 
where $t$ is the standard coordinate on the circle. 
The isomorphism between these two Lie algebras is given by 
the map 
\[ \vl_{k} \longmapsto - t^{k+1} \partial_{t}, \hspace{30pt} 
  T(\vl_{k}) \longmapsto - t^{k+1} \] 
for $k \geq -1$. To verify this statement, we only need to check 
the bracket relation 
\[ [T(\vl_{m}), \, \vl_{k}] = (m+1) T(\vl_{m+k}).\] 
This follows from the fact that 
$\nabla_{T(\vl_{m})} \vl_{k} = \vl_{m+k+1}$ by \eqref{eqn:DerTLk} 
and 
\begin{eqnarray*} 
\nabla_{\vl_{k}} T(\vl_{m}) & = & T( \nabla_{\vl_{k}} \vl_{m}) 
     - \vl_{k} \bullet \vl_{m} \\ 
&=& T(\tau_{-}(\vl_{m+k+1}) - (m+1) T(\vl_{m+k})) + \overline{\vl_{m+k+1}} \\ 
&=& \vl_{m+k+1} - (m+1) T(\vl_{m+k}) 
\end{eqnarray*} 
by \eqref{eqn:DerLk}, Lemma~\ref{lem:Virg0} and \eqref{eqn:WT}. 
Similar computations show that 
\[ [\vl_{k}, \, T^{n}(\vl_{j})] = \left( T^{n-1} R^{k+1} - R^{k+1}T^{n-1} 
    -(j+1) T^{n}R^{k}\right)(\vl_{j}) \] 
and 
\[ [T^{m}(\vl_{k}), \, T^{n}(\vl_{j})] 
    = \left( T^{n}R^{j+1}T^{m-1} R^{k+1} - 
    T^{m}R^{k+1}T^{n-1}R^{j+1} \right)(\vl_{-1}) \] 
for all $m, n \geq 1$ and $j, k \geq -1$. Therefore the set of 
vector fields $\{T^{m}(\vl_{k}) \mid m \geq 0, \, k \geq -1\}$ is 
closed under the Lie bracket (and consequently forms a Lie algebra) due to the 
following fact 
\begin{lem} 
\[ R T = T R + T^{2}.\] 
\end{lem} 
{\bf Proof}: For any vector field $\vw$, 
\[ R T(\vw) = T \tau_{-} R T(\vw) + \overline{R T(\vw)} \] 
by \eqref{eqn:WT}. Since $\overline{R T(\vw)} = \vx \bullet 
T(\vw) = 0$, by Lemma~\ref{lem:Rtau-}, 
\[ R T(\vw) = T \left( R \tau_{-}  T(\vw) + \vg * \overline{T(\vw)} + T(\vw) \right) 
    = T R(\vw) + T^{2}(\vw). \] 
The lemma is thus proved. $\Box$ 
 
Let $\partial_{t}^{-1}$ be the integral operator on the space of functions 
on the unit circle, where $t$ is the standard coordinate on the circle. Then 
as operators, 
\[ t \partial_{t}^{-1} = \partial_{t}^{-1} t + \left(\partial_{t}^{-1} \right)^{2}. \] 
Therefore the map 
\[ T^{m}(\vl_{k}) \longmapsto - \left(\partial_{t}^{-1}\right)^{m} t^{k+1} \partial_{t} \] 
for $k \geq -1$ and $m \geq 1$ defines an isomorphism between the Lie 
algebra $\{ T^{m}(\vl_{k}) \mid m \geq 1, \, k \geq -1\}$ and 
the Lie algebra of integral operators 
$\left\{ t^{k} \left(\partial_{t}^{-1} \right)^{m} \mid m, k \geq 0 \right\}$. 
Since $\vl_{k} = R^{k+1}(\vl_{-1})$, we see that under this isomorphism, $R$ 
corresponds to multiplying by $t$, and $T$ corresponds to composition by 
$\partial_{t}^{-1}$ on the space of pseudo-differential operators on the 
unit circle. As a left inverse of $T$, $\tau_{-}$ might be expected to correspond to 
$\partial_{t}$. However, the commutator of $\tau_{-}$ and $R$ 
 is not quite the same as the commutator of $\partial_{t}$ 
and $t$ due to the twisting given by $\vg*(\vs \bullet)$ in 
Lemma~\ref{lem:Rtau-}. 
 
To see the connection between vector fields $\vl_{n}$ with the 
Virasoro conjecture, we explicitly compute the first several 
Virasoro vector fields and obtain the following 
\begin{thm} \label{thm:LkRec} 
\begin{eqnarray} 
{\cal L}_{-1} & = & - {\cal S},    \nonumber \\ 
{\cal L}_{0} & = & - {\cal X} - (b_{1}+1) {\cal D}, \nonumber \\ 
{\cal L}_{1} & = & \sum_{m, \alpha} (m+b_{\alpha})(m+b_{\alpha}+1) 
        \tilde{t}^{\alpha}_{m} 
                \grava{m+1} 
         + \sum_{m, \alpha , \beta} (2m+2b_{\alpha}+1) 
                {\cal C}_{\alpha}^{\beta} 
        \tilde{t}^{\alpha}_{m} 
                \gravb{m}  \nonumber \\ 
       &&  + \sum_{m, \alpha, \beta} 
                ({\cal C}^{2})_{\alpha}^{\beta} 
        \tilde{t}^{\alpha}_{m} \gravb{m-1} 
    - \sum_{\alpha} b_{\alpha} (b_{\alpha} -1) \gwii{\gua} \ga, 
       \nonumber \\ 
{\cal L}_{2} & = & 
         \sum_{m, \alpha} (m+b_{\alpha})(m+b_{\alpha}+1) 
                (m+b_{\alpha}+2) 
        \tilde{t}^{\alpha}_{m} \grava{m+2}  \nonumber \\ 
        && + \sum_{m, \alpha, \beta} 
        \left\{ 3(m+b_{\alpha})^{2} + 6(m+b_{\alpha}) +2 \right\} 
     {\cal C}_{\alpha}^{\beta} 
        \tilde{t}^{\alpha}_{m} \gravb{m + 1} \nonumber \\ 
        && + \sum_{m, \alpha, \beta} 3(m+b_{\alpha}+1) 
                ({\cal C}^{2})_{\alpha}^{\beta} 
        \tilde{t}^{\alpha}_{m} \gravb{m} 
         + \sum_{m, \alpha, \beta} 
                ({\cal C}^{3})_{\alpha}^{\beta} 
        \tilde{t}^{\alpha}_{m} \gravb{m-1}   \nonumber \\ 
    && - \sum_{\alpha} b_{\alpha}(b_{\alpha}^{2} -1) 
        \{\gwii{\grava{1}} \gua + \gwii{\gua} \grava{1} \} 
     - \sum_{\alpha, \, \beta} (3 b_{\beta}^{2} -1) 
        {\cal C}_{\beta}^{\alpha} \gwii{\ga} \gub . 
    \nonumber 
\end{eqnarray} 
\end{thm} 
{\bf Proof}: The first equation is just the definition of $\vl_{-1}$. 
The second equation is precisely 
\eqref{eqn:RS}. 
The genus-0 quasi-homogeneity equation and the dilaton equation implies 
(cf. \cite[Lemma 1.2 and 1.4]{LT}) 
\[ \gwii{ {\cal L}_{0} \, \gua} = -(1-b_{\alpha}) \gwii{\gua} 
            - {\cal C}_{\beta}^{\alpha}\, t_{0}^{\beta}. \] 
The genus-0 $L_{1}$-constraint (cf. \cite[Formula (23)]{LT}) 
implies 
\[ \gwii{ {\cal L}_{1} \, \gua} 
    = -(1-b_{\alpha})(2-b_{\alpha}) \gwii{\gravua{1}} 
        - (3-2b_{\alpha}){\cal C}_{\beta}^{\alpha} \gwii{\gub} 
        - ({\cal C}^{2})_{\beta}^{\alpha} \, t_{0}^{\beta}. \] 
Using these formulas and the definitions of operators $T$ and $R$, 
it is straightforward to check the last two equations in the lemma. 
$\Box$ 
 
From the formulas is Theorem~\ref{thm:LkRec}, we see that 
$\vl_{-1}$ and $\vl_{0}$ are the first derivative parts of the 
first two Virasoro operators defined in \cite{EHX}. These vector 
fields were considered in \cite{LT}. The first derivative parts 
of the Virasoro operators $L_{1}$ and $L_{2}$ are also considered 
in \cite{LT} and played important role in the proof of the 
genus-0 Virasoro conjecture. These vector fields are the same as 
the linear parts of formulas for $\vl_{1}$ and $\vl_{2}$ 
(omitting the finitely many terms containing genus-0 1-point 
functions). The formulas for $\vl_{1}$ and $\vl_{2}$ as given in 
Theorem~\ref{thm:LkRec} do coincide with the corresponding vector 
fields used in \cite{DZ2} and \cite{G2}. These vector fields 
arise naturally when studying Virasoro conjecture of genus bigger 
than 0. More explicitly, for $g \geq 1$, the genus-$g$ 
$L_{n}$-constraint just computes $\gwiig{\vl_{n}}$ in terms of 
data with genus less than $g$. This is true not only for $-1 \leq 
n \leq 2$. For $n > 2$, this follows from 
Corollary~\ref{cor:BracketVir} since the first four Virasoro 
vector fields generate all others by taking iterating brackets. 
 
We can now also give a new interpretation for Lemma~\ref{lem:Virg0}. 
Using the formulas in Theorem~\ref{thm:LkRec}, we see that 
the most important cases of this lemma, i.e. 
\begin{equation} \label{eqn:VirOpEuler} 
{\cal L}_{1} \sim -{\cal X}^{2} \, \, \, \, {\rm and} \, \, \, \, 
{\cal L}_{2} \sim -{\cal X}^{3}, 
\end{equation} 
 are special cases of \cite[Equation (19) and (26)]{LT}, 
which are crucial steps in the proof of the genus-0 Virasoro conjecture. 
Because of equation~(\ref{eqn:equivPrimary}) 
these relations  are equivalent to 
\begin{equation} \label{eqn:VirOpEuler2} 
\overline{\vl_{1}} = -{\cal X}^{2} \, \, \, \, {\rm and} \, \, \, \, 
\overline{\vl_{2}} = -{\cal X}^{3}, 
\end{equation} 
which are special cases of \cite[Lemma 3.1 and Lemma 4.2]{LT} 
due to equation~(\ref{eqn:StringProd}) and \cite[Lemma 1.4]{LT}. 
We notice that Lemma~\ref{lem:Virg0} also follows from \eqref{eqn:VirOpEuler} 
due to Theorem~\ref{thm:VirEuler}. 
As we see from \cite{LT}, \eqref{eqn:VirOpEuler} is equivalent 
to the second derivatives of the genus-0 $L_{1}$ and $L_{2}$ constraints. 
Therefore, we can interpret Lemma~\ref{lem:Virg0} as the second 
derivatives of the genus-0 Virasoro conjecture. As explained in\cite{LT}, 
the second derivatives of the genus-0 Virasoro constraints imply 
the genus-0 Virasoro conjecture because of the dilaton equation. We notice 
that in the proof 
of Theorem~\ref{thm:LkRec}, the formula for $L_{k}$ is derived only using 
genus-0 $L_{k-1}$-constraint. Therefore one can modify the 
proof of Theorem~\ref{thm:LkRec} to give a new proof to the genus-0 
Virasoro conjecture using Lemma~\ref{lem:Virg0}. 
 
 To prove higher genera Virasoro conjecture, one needs to compute 
$\gwiig{\vl_{n}}$. Applying \eqref{eqn:WTW} to $\vl_{n}$ and consider the 
genus-$g$ topological recursion relation, we see that 
it is important to understand vector fields 
$\overline{ \tau_{-}^{m}({\cal L}_{n})}$. 
For $m=0$, 
$\overline{\tau_{-}^{m}(\vl_{n})}=\bvl_{n} = - \bvx^{n+1}$ by 
Lemma~\ref{lem:Virg0}. 
We will see that when $m \geq 1$, all such vector fields can be expressed 
in terms of certain twisted quantum powers of 
the Euler vector field. 
 We have the following 
\begin{thm} 
\label{thm:Stau-Lk} 
\[ \overline{\tau_{-}^{m}(\vl_{n+1})} = \vx \bullet \overline{\tau_{-}^{m}(\vl_{n})} 
 + m \, \overline{\tau_{-}^{m-1}(\vl_{n})} + \vg * \overline{\tau_{-}^{m-1}(\vl_{n})} 
\] for all $m \geq 1$ and $n \geq -1$. 
\end{thm} 
{\bf Proof}: 
 Using Lemma~\ref{lem:Rtau-}, we can prove 
inductively that 
\[ \tau_{-}^{m}(\vl_{n+1}) = R(\tau_{-}^{m}(\vl_{n})) + m \, 
\tau_{-}^{m-1}(\vl_{n}) + \vg * \overline{\tau_{-}^{m-1}(\vl_{n})} 
\] for all $m \geq 1$ and $n \geq -1$. The theorem then follows from 
Theorem~\ref{thm:Xprod}. 
$\Box$ 
 
Recursively applying Theorem~\ref{thm:Stau-Lk} and Lemma~\ref{lem:Virg0}, 
we can express 
$\overline{\tau_{-}^{m}({\cal L}_{n})}$ in terms of 
twisted quantum powers of ${\cal X}$ for any $m \geq 0$. Here the 
twisting is given by the operation $\vg *$. Such twisting is actually very 
important since it forces the sequence of vector fields 
$\{ \overline{\tau_{-}^{m}(\vl_{n})} \mid n \geq -1\}$ 
obey different linear relations for 
different $m$. This would make the span 
of vector fields 
$\{\sum_{m=0}^{k} T^{m}(\overline{\tau_{-}^{m}(\vl_{n})}) \mid n \geq -1\} $ 
large enough 
for each fixed $k$, and consequently make the Virasoro conjecture 
more interesting. 
We omit the explicit formulas for 
$\overline{\tau_{-}^{m}({\cal L}_{n})}$ 
for $m>1$ as it is not needed in this paper. 
To study the genus-2 Virasoro conjecture, we need the formula for 
$m=1$. 
In this case Theorem~\ref{thm:Stau-Lk} have the following form 
\[ \overline{\tau_{-}({\cal L}_{n+1})} 
= {\cal X} \bullet \overline{\tau_{-}({\cal L}_{n})} - \bvx^{n+1} 
    - {\cal G} * \bvx^{n+1} \] 
for $n \geq -1$. 
Recursively applying this formula, we obtain 
\begin{cor} \label{cor:Stau-Lk} 
\[ \overline{\tau_{-}({\cal L}_{n})} = 
-\bvx^{n+1} \bullet \tau_{-}({\cal S}) - (n+1)\bvx^{n} 
    - \sum_{j=0}^{n} \bvx^{j} \bullet 
        ({\cal G} * \bvx^{n-j}) \] 
for all $n \geq 0$. 
\end{cor} 
Note that for $m \geq 1$, $\tau_{-}^{m}(\vs)$ is zero when restricted 
to the small phase. Therefore, on the small phase space, the first 
term on the right hand side of this formula will disappear. 
 
It will be useful later to have a formula for  covariant 
derivatives of $\overline{\tau_{-}({\cal L}_{k})}$. 
\begin{cor} \label{cor:dertau-Lk} 
\[ \nabla_{T(\vw)} \overline{\tau_{-}(\vl_{k})} 
= \bvx^{k+1} \bullet \tau_{-}(\vw) + (k+1) \bvx^{k} \bullet \vw + 
    \sum_{i=0}^{k} \vx^{i} \bullet (\vg * (\bvx^{k-i} \bullet \vw)) \] 
for any vector field $\vw$ and $k \geq -1$. 
\end{cor} 
{\bf Proof}: 
By Lemma~\ref{lem:TWEulerPower}, Corollary~\ref{cor:TderProd}, and 
\eqref{eqn:derGC}, 
\[ \nabla_{T(\vw)} 
    \{\vx^{k} \bullet (\vg*\bvx^{m})\} 
= - \vx^{k} \bullet (\vg*(\vw \bullet \bvx^{m})) \] 
for any vector field $\vw$. The corollary then follows from 
Corollary~\ref{cor:Stau-Lk} and Corollary~\ref{cor:TWderTau-S}. 
$\Box$

\section{Applications to the Virasoro conjecture} 
\label{sec:AppVir} 
 
\subsection{The Virasoro Conjecture} 
\label{sec:VirConj} 
 
We will not describe the original version of the Virasoro conjecture as given in \cite{EHX}. 
Instead we will use the following formulation of this conjecture 
\[ \gwiig{\vl_{n}} = \rho_{g, n} \] 
for $g \geq 0$ and $n \geq -1$.  One of the advantages of this formulation is that the 
sequence of vector fields $\vl_{n}$ has a simple recursive definition 
(see \eqref{eqn:VirVF}). 
The $L_{-1}$-constraint is just the string equation. $L_{0}$-constraint 
was discovered by Hori and is a combination of the 
quasi-homogeneity equation and the dilaton equation.  Therefore for $n=-1$ and $0$, it is easy to figure out $\rho_{g,n}$  from these equations. 
For the $L_{1}$-constraint, we have 
\[ \rho_{0, 1} = \sum_{\alpha} \frac{1}{2} b_{\alpha}(1-b_{\alpha}) 
        \gwii{\ga} \gwii{\gua} - 
            \sum_{\alpha, \beta} \frac{1}{2} 
            \left({\cal C}^{2} \right)_{\alpha \beta} t_{0}^{\alpha} t_{0}^{\beta}, \] 
and for genus $g \geq 1$ 
\[ \rho_{g, 1} = -\sum_{\alpha} \frac{1}{2} b_{\alpha}(1-b_{\alpha}) 
        \left\{ \gwiih{g-1}{\ga \, \gua} 
                + \sum_{h=1}^{g-1} \gwiih{h}{\ga} \gwiih{g-h}{\gua} \right\}.\] 
For the $L_{2}$-constraint, we have 
\begin{eqnarray*} 
 \rho_{0, 2} & = & \sum_{\alpha} b_{\alpha}(1-b_{\alpha}^{2}) 
        \gwii{\grava{1}} \gwii{\gua} \\ 
    && +  \sum_{\alpha, \beta} \frac{1}{2} (1-3b_{\alpha}^{2}) {\cal C}_{\alpha}^{\beta} 
            \gwii{\gb} \gwii{\gua}   -  \frac{1}{2} 
       \left({\cal C}^{3} \right)_{\alpha \beta} t_{0}^{\alpha} t_{0}^{\beta}, 
\end{eqnarray*} 
and for genus $g \geq 1$ 
\begin{eqnarray*} 
\rho_{g, 2} &=& 
    - \sum_{\alpha} b_{\alpha}(1-b_{\alpha}^{2}) \left\{ \gwiih{g-1}{\grava{1} \, \gua} 
            + \sum_{h=1}^{g-1}\gwiih{h}{\grava{1}} \gwiih{g-h}{ \gua} \right\}  \\ 
  &&    + \sum_{\alpha, \beta} \frac{1}{2}(1-3b_{\alpha}^{2}) {\cal C}_{\alpha}^{\beta} 
        \left\{ \gwiih{g-1}{\gb \, \gua} +  \sum_{h=1}^{g-1}\gwiih{h}{\gb} 
                        \gwiih{g-h}{\gua} \right\}. 
\end{eqnarray*} 
Since $L_{k}$-constraint is generated by $L_{1}$ and $L_{2}$ 
constraints  for $k \geq 1$, this information suffices to determine 
the entire Virasoro conjecture. 
 
It may be desirable to have a description of $\rho_{g, n}$ in terms of 
recursive operators. For this purpose, we define 
\[ R_{+}(\vw) := \vg * \tau_{+}(\vw) + C(\vw) \] 
for any vector field $\vw$. 
Then for $n \geq 1$, 
\begin{eqnarray*} 
 \rho_{0, n} &=& \frac{1}{2} \sum_{i=0}^{n-1} \sum_{j=0}^{n-1-i} 
         \gwii{ \left\{ R_{+}^{j}(\vg * C^{i}(\ga)) \right\}} 
         \gwii{ \left\{ R_{+}^{n-1-i-j}(\vg * \gua) \right\}} \\ 
    &&  - \frac{1}{2} ({\cal C}^{n+1})_{\alpha \beta} 
            t_{0}^{\alpha} t_{0}^{\beta} 
\end{eqnarray*} 
and for $g \geq 1$, 
\begin{eqnarray*} 
 \rho_{g, n} &=& - \frac{1}{2} \sum_{i=0}^{n-1} \sum_{j=0}^{n-1-i} 
    \left\{  \gwiih{g-1}{ \left\{ R_{+}^{j}(\vg * C^{i}(\ga)) \right\} 
           \, \, 
        \left\{ R_{+}^{n-1-i-j}(\vg * \gua) \right\}} \right. \\ 
    &&  \hspace{80pt} + \sum_{h=1}^{g-1} \left. 
        \gwiih{h}{ \left\{ R_{+}^{j}(\vg * C^{i}(\ga)) \right\}} 
         \gwiih{g-h}{ \left\{ R_{+}^{n-1-i-j}(\vg * \gua) \right\}} 
            \right\}. 
\end{eqnarray*} 
Note that for $g >0$, $\rho_{g, n}$ depends only on the 
data with genus less than $g$. Therefore the genus-$g$ Virasoro conjecture just computes 
$\gwiig{\vl_{n}}$ in terms of data with genus less than $g$. 
 
One might also formulate the following {\it weak Virasoro conjecture}: 
$\gwiig{\vl_{n}}$ can be explicitly expressed in terms of data with 
genus less than $g$ for $g \geq 1$ and $n \geq -1$. This weak version 
of the Virasoro conjecture makes sense for all compact symplectic manifolds 
while the original version requires non-trivial topological 
conditions. From the computational point of view, once the weak Virasoro conjecture 
is proved, it will have the same computational power as the original Virasoro conjecture. 
This weak Virasoro conjecture can be generalized to the following {\it weak $W$-type constraints}: 
$\gwiig{\tau_{-}^{m}(\vl_{n})}$ can be explicitly expressed in terms of data with 
genus less than $g$ for $g \geq 1$, $n \geq -1$ and $m \geq 0$. 
The computation of $\overline{\tau_{-}^{m}(\vl_{n})}$ in terms of twisted quantum powers 
of the Euler vector field in Section~\ref{sec:VirVF} 
might be thought of as the genus-0 part of the $W$-constraints. 
While the work in \cite{DZ2} implies the genus-1 part and results in this paper imply 
the genus-2 part of the weak $W$-constraints for manifolds with semisimple quantum cohomology 
in a rather trivial way, the $W$-type constraints would be interesting only when the Virasoro conjecture 
can not completely determine the generating functions.

The relationship between the quantum powers of the Euler vector 
field and the Virasoro vector fields as revealed in 
Lemma~\ref{lem:Virg0} may be thought of as an interpretation 
of the second derivatives of the genus-0 Virasoro conjecture.  Vector fields ${\cal X}^{k+1}$ are also 
closely related to the genus-1 Virasoro conjecture. 
In fact, equation (\ref{eqn:WT}) 
and (\ref{eqn:TRRg1}) implies 
\begin{eqnarray*} 
\gwiione{{\cal L}_{k}} 
& = & \gwiione{T(\tau_{-}({\cal L}_{k}))} + \gwiione{ \overline{{L}_{k}} } \\ 
& = & \frac{1}{24} \gwii{\tau_{-}({\cal L}_{k}) \, \gua \, \ga} 
    - \gwiione{ \bvx^{k+1}}. 
\end{eqnarray*} 
Therefore the genus-1 $L_{k}$-constraint can be written as 
\[ \gwiione{\bvx^{k+1}} = \frac{1}{24} \gwii{\tau_{-}({\cal L}_{k}) \, \gua \, \ga} - \rho_{1, k}. \] 
It can be shown by using the genus-0 topological recursion relation and 
Lemma~\ref{lem:Virg0} that 
\[ \rho_{1, k} = - \frac{k+1}{8} \gwii{\vx^{k} \, \ga \, \gua} 
    + \frac{1}{4} \sum_{i=0}^{k} \ba \bb \gwii{\ga \, \vx^{i} \, \gub} 
            \gwii{\gb \, \vx^{k-i} \, \gua}. 
\] 
Together with Corollary~\ref{cor:Stau-Lk}, 
this explains the mysterious formulas in \cite{L1} for 
$\gwiione{\bvx^{k+1}}$ on the small phase space.

Since the genus-1 topological recursion relation is very powerful, 
the genus-1 Virasoro conjecture can be studied just using the 
quantum product on the small phase space. However, it seems that 
it is necessary to use the quantum product on the big phase space 
to study the genus-2 Virasoro conjecture since the genus-2 topological 
recursion relations are not strong enough. 
In fact, by \eqref{eqn:WT}, 
\begin{eqnarray*} 
{\cal L}_{k} 
& = & \bvl_{k} + T(\overline{\tau_{-}({\cal L}_{k})}) 
    + T^{2}(\tau_{-}^{2}({\cal L}_{k})). 
\end{eqnarray*} 
So by the genus-2 topological recursion relation~(\ref{eqn:TRR1}), 
\begin{eqnarray} 
\gwiitwo{{\cal L}_{k}} 
&= &  - \gwiitwo{ \bvx^{k+1} } + 
    \gwiitwo{T(\overline{\tau_{-}({\cal L}_{k})})} 
    + A_{1}(\tau_{-}^{2}({\cal L}_{k})). 
\label{eqn:Virg2} 
\end{eqnarray} 
The second term on the right hand side 
contains descendant vector fields which seems can not 
be reduced to primary vector fields by the known topological recursion 
relations. 
It is expected that 
the situation would become worse at higher genera. Therefore we believe that 
the study of structures on the quantum product on the big phase space 
is essential in the study of the higher genera Virasoro conjecture. 
 
Due to Corollary~\ref{cor:Stau-Lk}, the second term on the 
right hand side of \eqref{eqn:Virg2} can be expressed in terms 
of quantum powers of the Euler vector field. This is also true 
for the last term on the 
right hand side of \eqref{eqn:Virg2} because of Theorem~\ref{thm:Stau-Lk} 
and the following formula 
\begin{eqnarray} 
A_{1}(\vw) &= & A_{1}(\bvw) 
    + \frac{1}{20} \gwiione{ \left\{ \overline{\tau_{-}(\vw)} 
        \bullet \ga \bullet \gua \right\} } 
\nonumber \\ 
&&  + \frac{1}{480} \gwii{ \left\{ \overline{\tau_{-}(\vw)} 
        \bullet \ga \right\} \, \gua \, \gb \, \gub } 
    + \frac{1}{1152} \gwii{ \overline{\tau_{-}(\vw)} \, 
        \left\{ \ga \bullet \gua \right\} 
            \, \gb \, \gub} 
\nonumber \\ 
&&  + \frac{1}{1152} \gwii{ \left\{ \overline{\tau_{-}^{2}(\vw)} 
        \bullet \ga \bullet \gua \right\} 
        \, \gb \, \gub }. 
\label{eqn:A1WWbar} 
\end{eqnarray} 
This equation follows from the definition of $A_{1}$, 
Corollary~\ref{cor:4ptT},  equations (\ref{eqn:5ptT}) and 
(\ref{eqn:derTRRg1}). 
 
In the rest of this section we will study the genus-2 Virasoro 
conjecture. We also solve the genus-2 generating function in terms 
of genus-0 and genus-1 data for manifolds whose quantum cohomology 
is not too degenerate. 
 
\subsection{Reduce the genus-2 Virasoro conjecture to the genus-2 
$L_{1}$-constraint} 
\label{sec:L1->L2} 
 
 In this subsection, we 
prove that the genus-2 Virasoro conjecture can be reduced 
to the $L_{1}$-constraint provided that the genus-1 $L_{1}$-constraint 
is satisfied. 
Note that the genus-1 $L_{1}$-constraint holds if and only if 
 the genus-1 Virasoro 
conjecture holds (cf. \cite{L1}). 
For a discussion of some sufficient conditions for the 
genus-1 Virasoro conjecture, see \cite{L1} and \cite{L2}. 
 
Define 
\begin{equation} 
\label{eqn:psi} \psi_{k} := \gwiitwo{\bvx^{k}} - 
        \gwiitwo{T(\overline{\tau_{-}(\vl_{k-1})})}. 
\end{equation} 
Then genus-2 $L_{k}$-constraint have the following form 
\[ \psi_{k+1} = A_{1}(\tau_{-}^{2}(\vl_{k})) - \rho_{2, k}. 
\] 
 
We first apply \eqref{eqn:BP} to the case $W_{i} = \vx$ and obtain 
\begin{eqnarray*} 
B(\vx, \vx, \vx) 
&= & 2 \gwiitwo{\vx^{3}} - 2 \gwii{\vx \, \vx \, \vx \, \gua} \gwiitwo{T(\ga)} \\ 
&& - 3 \gwiitwo{T(\vx) \, \, \vx^{2}} + 3 \gwiitwo{\vx \, \, T(\vx^{2})}. 
\end{eqnarray*} 
Since 
\begin{eqnarray*} 
 \gwiitwo{T(\vx) \, \, \vx^{2}} 
= T(\vx) \gwiitwo{\vx^{2}} - \gwiitwo{\left\{ \nabla_{T(\vx)} \vx^{2} \right\} } 
 = T(\vx) \gwiitwo{\vx^{2}} + \gwiitwo{\vx^{3} } 
\end{eqnarray*} 
and 
\begin{eqnarray*} 
 \gwiitwo{ \vx \, \, T(\vx^{2}) } 
& = & T(\vx^{2}) \gwiitwo{\vx} - \gwiitwo{\left\{ \nabla_{T(\vx^{2})} \vx \right\} } \\ 
& = & -(3b_{1}+2) \gwiitwo{T(\vx^{2})} + \gwiitwo{T(\vg * \vx^{2})}+ \gwiitwo{\vx^{3} }, 
\end{eqnarray*} 
we have 
\begin{eqnarray} 
B(\vx, \vx, \vx) 
=  2 \gwiitwo{\vx^{3}} - 3 T(\vx) \gwiitwo{\vx^{2}} -  \gwiitwo{T(\vv_{1})} 
\label{eqn:BPXXX} 
\end{eqnarray} 
where 
\[ \vv_{1} = 2\gwii{\vx \, \vx \, \vx \, \gua} \ga + 3 (3b_{1}+2) \vx^{2} - 3\, \vg * \vx^{2}. \] 
Using this formula, we can prove the following 
\begin{lem} \label{lem:L1L2} 
If the genus-2 $L_{1}$-constraint holds, then 
\begin{eqnarray*} 
\psi_{3} &=& \frac{3}{2} \left\{ A_{1}(\nabla_{T(\vx)} \tau_{-}^{2}(\vl_{1})) 
        + A_{2}(\vx, \tau_{-}(\vl_{1})) \right\} 
        + \frac{1}{2} B(\vx, \vx, \vx) 
        - \frac{3}{2} T(\vx) \rho_{2,1}. 
\end{eqnarray*} 
\end{lem} 
{\bf Proof}: 
We first observe that by Lemma~\ref{lem:XXX4pt}, Corollaries \ref{cor:Stau-Lk} and \ref{cor:dertau-Lk}, 
and \eqref{eqn:tau-X}, 
the vector field $\vv_{1}$ defined after \eqref{eqn:BPXXX} also satisfies 
\[ \vv_{1} = 2 \, \overline{\tau_{-}(\vl_{2})} - 9 \, \vx \bullet \overline{\tau_{-}(\vl_{1})} 
    - 3 \, \nabla_{T(\vx)} \overline{\tau_{-}(\vl_{1})}. \] 
Secondly, by \eqref{eqn:T}, 
\[ T(\tau_{-}(\vw)) = \tau_{+}(\tau_{-}(\vw)) - \overline{\tau_{+}(\tau_{-}(\vw))} 
    = \vw - \bvw \] 
for any vector field $\vw$. Therefore by Corollary~\ref{cor:A1A2}, 
\begin{eqnarray*} 
T(\vx) \, A_{1}(\tau_{-}^{2}(\vl_{1})) 
&=& A_{2}\left(\vx, T(\tau_{-}^{2}(\vl_{1}))\right)  + 
    A_{1}( \nabla_{T(\vx)} \tau_{-}^{2}(\vl_{1}))  \\ 
&=& A_{2}\left(\vx, \tau_{-}(\vl_{1}) - \overline{\tau_{-}(\vl_{1})}\right) + 
    A_{1}( \nabla_{T(\vx)} \tau_{-}^{2}(\vl_{1})). 
\end{eqnarray*} 
Moreover, by \eqref{eqn:TRR2}, 
\begin{eqnarray*} 
T(\vx)\gwiitwo{ T(\overline{\tau_{-}(\vl_{1})}) } 
&=& \gwiitwo{ T(\vx) \, \, T(\overline{\tau_{-}(\vl_{1})}) } + 
    \gwiitwo{ \left\{ \nabla_{T(\vx)} T(\overline{\tau_{-}(\vl_{1})}) \right\}}  \\ 
&=& A_{2}(\vx, \overline{\tau_{-}(\vl_{1})}) + 
    3 \gwiitwo{ T(\vx \bullet \overline{\tau_{-}(\vl_{1})}) } + 
    \gwiitwo{ T(\nabla_{T(\vx)} \overline{\tau_{-}(\vl_{1})})}. 
\end{eqnarray*} 
If genus-2 $L_{1}$-constraint is satisfied, then 
\[ \gwiitwo{\vx^{2}} = \gwiitwo{ T(\overline{\tau_{-}(\vl_{1})}) } + 
    A_{1}(\tau_{-}^{2}(\vl_{1})) 
    - \rho_{2, 1}. 
    \] 
The lemma then follows by plugging this equation into \eqref{eqn:BPXXX}.  $\Box$ 
 
\begin{lem} \label{lem:dertau-2L1} 
\[ \nabla_{T(\vx)} \tau_{-}^{2}(\vl_{1}) = - \tau_{-}^{2}(\vl_{2}) + 
    (b_{1} +1) \tau_{-}(\vl_{1}) + 2(b_{1}+1) \vl_{0} - 2 (b_{1} + 1) \vd.\] 
\end{lem} 
{\bf Proof}: 
By Lemma~\ref{lem:derLkv2} and the fact that covariant derivatives 
commute with $\tau_{-}$, 
\[ \nabla_{T(\vx)} \tau_{-}^{2}(\vl_{1}) =  \tau_{-}^{2}R^{2}(\vx) 
    = \tau_{-}^{2}R^{2}(-\vl_{0}-(b_{1}+1) \vd) 
    = -\tau_{-}^{2}(\vl_{2}) - (b_{1}+1) \tau_{-}^{2}R^{2}(\vd) . \] 
Since $\bvd = 0$, by Lemma~\ref{lem:Rtau-}, 
\[ \tau_{-}R(\vd) = R\tau_{-}(\vd) + \vd = R(\vs) + \vd = - \vl_{0}+\vd.\] 
So 
\[ \tau_{-}R\tau_{-}R(\vd) = - \tau_{-}(\vl_{1}) - \vl_{0} + \vd\] 
 and 
\[\tau_{-}^{2}R^{2}(\vd) = \tau_{-}R\tau_{-}R(\vd) + 
    \tau_{-}(\vg*\overline{R(\vd)}) + \tau_{-}R(\vd) 
    = - \tau_{-}(\vl_{1}) - 2 \vl_{0} + 2 \vd. \] 
The lemma follows. 
$\Box$ 
 
\begin{lem} \label{lem:A1L0-D} 
\begin{eqnarray*} 
A_{1}(\vd) &=& \frac{1}{20} \gwiione{\left\{ \gua \bullet \ga \right\}} 
 + \frac{1}{480} \gwii{\ga \, \gua \, \gb \, \gub}, 
\end{eqnarray*} 
and 
\begin{eqnarray*} 
A_{1}(\vl_{0}) 
&=& - \frac{7}{10}  \gwiione{\left\{\vx \bullet \gua \right\}} \gwiione{\ga} 
 - \frac{1}{10}  \gwiione{ \left\{\vx \bullet \gua \right\}  \, \ga} \\ 
&& + \frac{1}{120} (7b_{\beta} -13)  \gwii{\gua \, \ga \, \gub} \gwiione{\gb} \\ 
&& - \frac{1}{480} \gwii{\ga \, \gua \, \gb \, \gub}. 
\end{eqnarray*} 
\end{lem} 
{\bf Proof}: The formula for $A_{1}(\vd)$ follows from derivatives 
of the dilaton equation 
\[ \gwiig{\vd \, \vw_{1} \cdots \, \vw_{k}} 
= (k+2g-2) \gwiig{\vw_{1} \, \cdots \vw_{k}} \] 
for any vector fields $\vw_{1}, \ldots, \vw_{k}$. 
The formula for $A_{1}(\vl_{0})$ is obtained by applying derivatives of 
the $L_{0}$-constraint 
\begin{eqnarray*} 
 \gwiig{\vl_{0} \, \vw_{1} \cdots \, \vw_{k}} 
&=& - \sum_{i=1}^{k} \gwiig{\vw_{1} \, \cdots \, 
    \left\{ \tau_{-}R(\vw_{i}) \right\} 
    \, \cdots \, \vw_{k}}  \\ 
&& +k \gwiig{\vw_{1} \, \cdots \, \vw_{k}} 
    - \delta_{g, 0} \nabla^{k}_{\vw_{1}, \cdots, \vw_{k}} 
    \left( \frac{1}{2} {\cal C}_{\alpha \beta} 
    t_{0}^{\alpha} t_{0}^{\beta} \right) 
\end{eqnarray*} 
to remove $\vl_{0}$ from genus-0 correlation functions with more than 3 points 
and genus-1 correlation functions. 
$\Box$ 
 
To compute $B(\vx, \vx, \vx)$, we need the following 
\begin{lem} \label{lem:X2g1} 
If genus-1 $L_{1}$-constraint holds, then 
\begin{eqnarray*} 
\gwiione{\vx^{2} \, \gua} &=& (1-b_{\alpha}+b_{\beta}) 
    \gwii{\vx \, \gua \, \gub}  \gwiione{\gb} \\ 
    && + \left\{\frac{1}{2} b_{\beta}(1-b_{\beta}) + 
    \frac{1}{24}(1-b_{\alpha})(2-b_{\alpha}) \right\} 
    \gwii{\gua \, \gb \, \gub} \\ 
    &&+ \gwiione{ \left\{ \tau_{-}(\vl_{1}) \bullet \gua \right\}} 
    + \frac{1}{24} \gwii{ \tau_{-}(\vl_{1}) \, \gua \, \gb \, \gub} 
\end{eqnarray*} 
and 
\begin{eqnarray*} 
\gwiione{\vx^{2} \, \gua \, \ga} &=& 2(1-b_{\alpha}+b_{\beta}) 
    \gwii{\vx \, \gua \, \gub}  \gwiione{\ga \, \gb} \\ 
    && + \left\{(1- b_{\alpha}+b_{\beta})(2- b_{\alpha}-b_{\beta}) + 
    b_{\alpha}(b_{\alpha} + 1) \right\} 
    \gwii{\gua \, \ga \, \gub} \gwiione{\gb} \\ 
    &&+ \left\{\frac{1}{2} b_{\beta}(1-b_{\beta}) + 
        \frac{1}{24}(1-b_{\alpha})(2-b_{\alpha}) 
        +\frac{1}{24}b_{\alpha}(b_{\alpha}+1) \right\} 
    \gwii{\ga \, \gua \, \gb \, \gub} \\ 
    &&+ 2\gwiione{ \left\{ \tau_{-}(\vl_{1}) \bullet \gua \right\} 
            \, \ga} 
     + \gwii{\tau_{-}(\vl_{1}) \, \gua \, \ga \, \gub} \gwiione{\gb} \\ 
    &&+ \frac{1}{24} \gwii{ \tau_{-}(\vl_{1}) \, \ga \, \gua \, \gb \, \gub}. 
\end{eqnarray*} 
\end{lem} 
{\bf Proof}: First observe that 
\[ \vx^{2} = - \bvl_{1} = - \vl_{1} + T(\tau_{-}(\vl_{1})) \] 
and 
\begin{eqnarray*} 
\nabla_{\gm} \vl_{1} &=& \nabla_{\gm} R(\vl_{0}) 
\, = \, R(\nabla_{\gm} \vl_{0}) - \vg * (\vl_{0} \bullet \gm) \\ 
&=& \vx \bullet (b_{\mu} \gm) + T(\tau_{-}R(b_{\mu} \gm)) +\vg * (\vx \bullet \gm) \\ 
&=& (b_{\mu} + b_{\alpha}) \gwii{\vx \, \gm \, \gua} \ga + b_{\mu}(b_{\mu}+1) T(\gm). 
\end{eqnarray*} 
The formulas in the lemma are obtained by applying equations 
(\ref{eqn:derTRRg1}) and (\ref{eqn:2derTRRg1}) and taking derivatives of 
the genus-1 $L_{1}$-constraint which has the following 
form 
\[ \gwiione{\vl_{1}} = - 
    \sum_{\alpha} \frac{1}{2} b_{\alpha}(1-b_{\alpha}) 
        \gwii{\ga \, \gua}. 
    \] 
$\Box$

\begin{lem} \label{lem:3A2+BXXX} 
If genus-1 $L_{1}$-constraint holds, then 
\begin{eqnarray*} 
&& 3 A_{2}(\vx, \tau_{-}(\vl_{1})) + B(\vx, \vx, \vx) \\ 
&=& 5 A_{1}(\tau_{-}^{2}(\vl_{2})) - 3(b_{1} + 1) A_{1}(\tau_{-}(\vl_{1})) \\ 
&& + \frac{1}{5} \left\{ 5 b_{\alpha}(b_{\alpha}+ b_{\beta}) - 5 b_{\alpha} + 
        21(b_{1}+1) \right\} \gwii{\vx \, \gua \, \gub} \gwiione{\ga} \gwiione{\gb} \\ 
&& + \frac{1}{10} \left\{ 10 b_{\alpha}(b_{\alpha}+ b_{\beta}) - 10 b_{\alpha} + 
        6(b_{1}+1) \right\} \gwii{\vx \, \gua \, \gub} \gwiione{\ga \, \gb} \\ 
&& + \frac{1}{120} \left\{ -5 b_{\beta}^{3} - 15 b_{1} b_{\beta}^{2} 
    -27 b_{1} b_{\beta} - 37 b_{\beta} + 114(b_{1}+1)   \right. \\ 
&& \hspace{60pt}   \left.  +180 (b_{\beta}+b_{1}+1)b_{\alpha}(1-b_{\alpha}) \right\} 
        \gwii{\gua \, \ga \, \gub} \gwiione{\gb} \\ 
&& + \frac{1}{40}(b_{1}+1)(-5 b_{\alpha}^{2} + 5 b_{\alpha} +1) 
        \gwii{\ga \, \gua \, \gb \, \gub}. 
\end{eqnarray*} 
\end{lem} 
{\bf Proof}: We first observe that 
\begin{eqnarray*} 
 \tau_{-}R(\vx) & = & \tau_{-}R(-\vl_{0}-(b_{1}+1)\vd) 
    \, = \, -\tau_{-}(\vl_{1}) - (b_{1}+1)\tau_{-}R(\vd) \\ 
    &=& -\tau_{-}(\vl_{1}) + (b_{1}+1)\vl_{0} - (b_{1}+1)\vd. 
\end{eqnarray*} 
Similarly, we have 
\[ \tau_{-}R\tau_{-}R(\vx) = 
    - \tau_{-}^{2}(\vl_{2}) + (b_{1}+2) \tau_{-}(\vl_{1}) 
    + (b_{1}+1)\vl_{0} - (b_{1}+1)\vd \] 
and 
\[ \tau_{-}R(\tau_{-}(\vl_{1})) = 
    \tau_{-}^{2}(\vl_{2}) - \tau_{-}(\vl_{1}). \] 
Remove $\vx$ from genus-0 correlation functions 
with more than 3 points and genus-1 correlation functions 
in $A_{2}(\vx, \tau_{-}(\vl_{1}))$ and $B(\vx, \vx, \vx)$ using 
Lemma~\ref{lem:rmX} and the above formulas. Then use 
Lemma~\ref{lem:X2g1} to replace $\gwiione{\vx^{2} \, \gua}$ and 
$\gwiione{\vx^{2} \, \gua \, \ga}$. 
The lemma then follows from simplifying the resulting expression. 
$\Box$ 
 
The last term of the formula in Lemma~\ref{lem:L1L2} can be computed in the 
following way 
\begin{lem} \label{lem:TXg1L1} 
\begin{eqnarray*} 
-T(\vx) \rho_{2,1} 
&=& b_{\alpha}(1-b_{\alpha}) \left\{ 
    \gwiione{\left\{\vx \bullet \gua \right\}} \gwiione{\ga} 
    +  \gwiione{ \left\{\vx \bullet \gua \right\}  \, \ga} \right\} \\ 
&& + \frac{1}{2} \left\{ b_{\alpha}(1-b_{\alpha}) 
        + \frac{1}{12}b_{\beta}(1-b_{\beta}) \right\} 
 (1-b_{1}-b_{\beta})  \gwii{\gua \, \ga \, \gub} \gwiione{\gb} \\ 
&& - \frac{1}{24} b_{1} b_{\alpha}(1-b_{\alpha}) 
    \gwii{\ga \, \gua \, \gb \, \gub}. 
\end{eqnarray*} 
\end{lem} 
{\bf Proof}: 
This formula is obtained by applying equations (\ref{eqn:derTRRg1}) 
and (\ref{eqn:2derTRRg1}) first, then remove $\vx$ from genus-0 
correlation functions with more than 3 points. 
$\Box$ 
 
To write the prediction of the genus-2 $L_{2}$ constraint 
in a form consistent with the above calculations, we need the following 
\begin{lem} \label{lem:3bto2b} 
\[ \sum_{\alpha} b_{\alpha}(1-b_{\alpha}^{2}) 
    \gwiig{\ga \, \gua \, \vw_{1} \, \cdots \vw_{k}} 
= \frac{3}{2} \, \sum_{\alpha}  b_{\alpha}(1-b_{\alpha}) 
    \gwiig{\ga \, \gua \, \vw_{1} \, \cdots \vw_{k}} \] 
for any vector fields $\vw_{1}, \ldots, \vw_{k}$. 
\end{lem} 
{\bf Proof}: The difference of the right hand side and the left hand 
side of this equation is 
\begin{eqnarray*} 
&& \frac{1}{2} \, \sum_{\alpha}  b_{\alpha}(1-b_{\alpha})(1-2b_{\alpha}) 
    \gwiig{\ga \, \gua \, \vw_{1} \, \cdots \vw_{k}} \\ 
& = & \frac{1}{2} \, \sum_{\alpha, \beta} 
     b_{\alpha}(1-b_{\alpha})(1-2b_{\alpha}) \eta^{\alpha \beta} 
    \gwiig{\ga \, \gb \, \vw_{1} \, \cdots \vw_{k}} \\ 
&=&  \frac{1}{2} \, \sum_{\alpha, \beta} 
     (1-b_{\beta})b_{\beta}(-1+2b_{\beta}) \eta^{\alpha \beta} 
    \gwiig{\ga \, \gb \, \vw_{1} \, \cdots \vw_{k}} \\ 
&=&  - \frac{1}{2} \, \sum_{\beta} 
     b_{\beta}(1-b_{\beta})(1-2b_{\beta}) 
    \gwiig{\gub \, \gb \, \vw_{1} \, \cdots \vw_{k}}. 
\end{eqnarray*} 
Here we have used the fact that $b_{\alpha} + b_{\beta}=1$ if 
$\eta^{\alpha \beta} \neq 0$. Comparing the two sides of this equation, 
both of them must be zero. 
$\Box$ 
 
Now we can rewrite the prediction of $-\gwiitwo{\vl_{2}}$ given by the Virasoro conjecture 
in the following form 
\begin{lem} \label{lem:g2L2v2} 
\begin{eqnarray*} 
- \rho_{2, 2} 
&=& \frac{1}{2} \left\{ -2 b_{\alpha}^{2} + 2b_{\alpha}+ b_{\alpha}b_{\beta} \right\} 
    \gwii{\vx \, \gua \, \gub}  \left\{ \gwiione{\ga} \gwiione{\gb} 
            + \gwiione{\ga \, \gb} \right\} \\ 
&& +  \left\{ \frac{3}{2} b_{\alpha}(1-b_{\alpha}) 
        + \frac{1}{24}b_{\beta}(1-b_{\beta})(2-b_{\beta}) \right\} 
        \gwii{\gua \, \ga \, \gub} \gwiione{\gb} \\ 
&& + \frac{1}{16}b_{\alpha}(1- b_{\alpha}) 
        \gwii{\ga \, \gua \, \gb \, \gub}. 
\end{eqnarray*} 
\end{lem} 
{\bf Proof}: 
Applying equation (\ref{eqn:derTRRg1}), (\ref{eqn:2derTRRg1}) and 
Lemma~\ref{lem:3bto2b}, we have 
\begin{eqnarray*} 
&&\sum_{\alpha} b_{\alpha}(1-b_{\alpha}^{2}) 
    \left\{ \gwiione{\grava{1} \, \gua} 
            + \gwiione{\grava{1}} \gwiione{ \gua} \right\}  \\ 
&=& \sum_{\alpha, \beta} b_{\alpha}(1-b_{\alpha}^{2}) 
    \gwii{\ga \, \gb} \left\{ \gwiione{\gua \, \gub} 
            + \gwiione{\gua} \gwiione{ \gub} \right\}  \\ 
&& +  \left\{ \frac{3}{2} b_{\alpha}(1-b_{\alpha}) 
        + \frac{1}{24}b_{\beta}(1-b_{\beta})(2-b_{\beta}) \right\} 
        \gwii{\gua \, \ga \, \gub} \gwiione{\gb} \\ 
&& + \frac{1}{16}b_{\alpha}(1- b_{\alpha}) 
        \gwii{\ga \, \gua \, \gb \, \gub}. 
\end{eqnarray*} 
By the symmetry of $\alpha$ and $\beta$, the first term on the 
right hand side can be written as 
\begin{eqnarray*} 
&& \sum_{\alpha, \beta} \, \frac{1}{2} 
    \left\{ b_{\alpha}(1-b_{\alpha}^{2}) + b_{\beta}(1-b_{\beta}^{2}) 
        \right\} 
    \gwii{\ga \, \gb} \left\{ \gwiione{\gua \, \gub} 
            + \gwiione{\gua} \gwiione{ \gub} \right\}  \\ 
&=&  \sum_{\alpha, \beta} \, \frac{1}{2} 
    \left(1- b_{\alpha}^{2} + b_{\alpha}b_{\beta} -b_{\beta}^{2} \right) 
    \left\{ (b_{\alpha} + b_{\beta}) \gwii{\ga \, \gb} \right\} 
    \left\{ \gwiione{\gua \, \gub} 
            + \gwiione{\gua} \gwiione{ \gub} \right\} 
\end{eqnarray*} 
By \eqref{eqn:X3pt}, 
\begin{eqnarray*} 
(b_{\alpha} + b_{\beta}) \gwii{\ga \, \gb} 
&=& \gwii{ \vx \, \ga \, \gb} - {\cal C}_{\alpha \beta}. 
\end{eqnarray*} 
Since ${\cal C}_{\alpha \beta} \neq 0$ implies $b_{\alpha} = - b_{\beta}$, 
\[ \left(1- b_{\alpha}^{2} + b_{\alpha}b_{\beta} -b_{\beta}^{2} \right) 
  {\cal C}_{\alpha \beta} = (1-3b_{\alpha}^{2}) {\cal C}_{\alpha \beta}.\] 
The lemma then follows from interchanging upper indices and lower indices and 
using the symmetry of $\alpha$ and $\beta$. 
$\Box$ 
 
We can now prove 
\begin{thm} \label{thm:L1->L2} 
For any manifold which satisfies the genus-1 $L_{1}$-constraint, 
if the genus-2 $L_{1}$-constraint holds, then the genus-2 
Virasoro conjecture holds. 
\end{thm} 
{\bf Proof}: Combining the results in Lemmas 
\ref{lem:L1L2} - \ref{lem:A1L0-D}, \ref{lem:3A2+BXXX} and \ref{lem:TXg1L1}, 
we obtain a formula for $\psi_{3}$. 
Due to Lemma~\ref{lem:g2L2v2}, this formula coincides with the prediction 
of the genus-2 $L_{2}$-constraint. Since $L_{k}$-constraint is 
generated by $L_{-1}$ and $L_{2}$-constraints, the theorem is thus proved. 
$\Box$ 
 
\subsection{A recursive formula for $\psi_{k}$} 
 
In this subsection, we prove a recursive formula for the function 
 $\psi_{k}$ which was defined 
in \eqref{eqn:psi}. 
We first apply \eqref{eqn:BP} to $\vw_{i} = 
\bvx^{m_{i}}$, for $i=1, 2, 3$ and $m_{i} \geq 0$ and obtain 
\begin{eqnarray} 
&& 2 \gwiitwo{\bvx^{m}} 
- 2 \gwii{\bvx^{m_{1}} \, 
    \bvx^{m_{2}} \, \bvx^{m_{3}} \, \gua} 
        \gwiitwo{T(\ga)} 
\nonumber \\ 
&& + \sum_{i=1}^{3} \left\{ \gwiitwo{\bvx^{m_{i}} \, 
        T(\bvx^{m-m_{i}})} - 
\gwiitwo{ T(\bvx^{m_{i}}) \, \bvx^{m-m_{i}} } 
            \right\} 
\nonumber \\ 
& = & B(\bvx^{m_{1}}, \bvx^{m_{2}}, \bvx^{m_{k}}), 
\label{eqn:BPE} 
\end{eqnarray} 
where $m=m_{1}+m_{2}+m_{3}$. 
 
To simplify this equation, we notice that, by 
Lemma~\ref{lem:TWEulerPower}, 
\begin{eqnarray*} 
\gwiitwo{\bvx^{n} \, T(\bvx^{k})} 
&=& T(\bvx^{k}) \gwiitwo{ \bvx^{n}} 
        - \gwiitwo{\{\nabla_{T(\bvx^{k})} 
            \bvx^{n}\}} \\ 
&=& T(\bvx^{k}) \gwiitwo{ \bvx^{n}} 
        + \gwiitwo{\bvx^{n+k} } 
\end{eqnarray*} 
for any $n, \, k \geq 0$. Hence \eqref{eqn:BPE} can be written as 
\begin{eqnarray} 
&& 2 \gwiitwo{\bvx^{m} } 
- 2 \gwii{\bvx^{m_{1}} \, \bvx^{m_{2}} \, \bvx^{m_{3}} \, \gua} 
        \gwiitwo{T(\ga)} 
\nonumber \\ 
&& +\sum_{i=1}^{3} \left\{ 
T(\bvx^{m-m_{i}}) \gwiitwo{\bvx^{m_{i}}} 
-T(\bvx^{m_{i}}) \gwiitwo{\bvx^{m-m_{i}}} 
\right\} 
\nonumber \\ 
& = & B(\bvx^{m_{1}}, \bvx^{m_{2}}, \bvx^{m_{3}}). 
\label{eqn:BPE2} 
\end{eqnarray} 
To simplify this equation further, we need the following lemma. 
\begin{lem} 
\label{lem:TXderTg2} 
\begin{eqnarray*} 
&& T(\bvx^{k}) \gwiitwo{T(\overline{\tau_{-}(\vl_{n-1})})} 
-T(\bvx^{n}) \gwiitwo{T(\overline{\tau_{-}(\vl_{k-1})})} \\ 
&=& 2(k-n) \gwiitwo{T(\bvx^{k+n-1})} 
    +2 \sum_{i=0}^{k-1} \gwiitwo{T(\vx^{i} \bullet 
        (\vg * \bvx^{k+n-1-i}))}   \\ 
&&  -  2 \sum_{i=0}^{n-1} \gwiitwo{T(\vx^{i} \bullet 
        (\vg * \bvx^{k+n-1-i}))} \\ 
&& + A_{2}(\bvx^{k},  \overline{\tau_{-}(\vl_{n-1})}) 
    -A_{2}(\bvx^{n},  \overline{\tau_{-}(\vl_{k-1})}) 
\end{eqnarray*} 
for $k, n \geq 0$. 
\end{lem} 
{\bf Proof}: 
By \eqref{eqn:DerCorr}, 
\begin{eqnarray*} 
&& T(\bvx^{k}) \gwiitwo{T(\overline{\tau_{-}(\vl_{n-1})})} \\ 
&=& \gwiitwo{T(\bvx^{k}) \, \, T(\overline{\tau_{-}(\vl_{n-1})})} 
    + \gwiitwo{ \left\{\nabla_{T(\bvx^{k})} 
    T(\overline{\tau_{-}(\vl_{n-1})}) \right\}} . 
\end{eqnarray*} 
By \eqref{eqn:TRR2} and Lemma~\ref{lem:DerT}, 
\begin{eqnarray*} 
&& T(\bvx^{k}) \gwiitwo{T(\overline{\tau_{-}(\vl_{n-1})})} \\ 
&=& 3 \gwiitwo{T(\vx^{k} \bullet \overline{\tau_{-}(\vl_{n-1})})} 
    + A_{2}(\bvx^{k},  \overline{\tau_{-}(\vl_{n-1})}) \\ 
&&  + \gwiitwo{ T(\nabla_{T(\bvx^{k})} 
            \overline{\tau_{-}(\vl_{n-1})}))} . 
\end{eqnarray*} 
The lemma then follows from Corollary~\ref{cor:Stau-Lk} 
and Corollary~\ref{cor:dertau-Lk}. 
$\Box$ 
 
\begin{thm} 
\label{thm:psi} 
For $m_{1}, m_{2}, m_{3} \geq 0$ and $m=m_{1}+m_{2}+m_{3}$, 
\begin{eqnarray*} 
&& 2 \psi_{m} 
 + \sum_{i=1}^{3} \left\{ T(\bvx^{m-m_{i}})\psi_{m_{i}} 
     -T(\bvx^{m_{i}}) \psi_{m-m_{i}} \right\} 
\nonumber \\ 
& = & B(\bvx^{m_{1}}, \bvx^{m_{2}}, \bvx^{m_{3}}) \\ 
&& + \sum_{i=1}^{3} \left\{ 
    A_{2}(\bvx^{m_{i}},  \, 
            \overline{\tau_{-}(\vl_{m-m_{i}-1})}) 
    -A_{2}(\bvx^{m-m_{i}}, \, 
            \overline{\tau_{-}(\vl_{m_{i}-1})}) 
    \right\}. 
\end{eqnarray*} 
\end{thm} 
{\bf Proof}: Plugging $\gwiitwo{\bvx^{k}} = \psi_{k} + 
        \gwiitwo{T(\overline{\tau_{-}(\vl_{k-1})})}$ 
into \eqref{eqn:BPE2} and applying Lemma~\ref{lem:TXderTg2}, we obtain 
\begin{eqnarray*} 
&& 2 \psi_{m} + 2 \gwiitwo{ T({\cal Y})} 
 + \sum_{i=1}^{3} \left\{ T(\bvx^{m-m_{i}})\psi_{m_{i}} 
     -T(\bvx^{m_{i}}) \psi_{m-m_{i}} \right\} \\ 
&&  + \sum_{i=1}^{3} \left\{ 
        A_{2}(\bvx^{m-m_{i}}, \, 
            \overline{\tau_{-}(\vl_{m_{i}-1})}) 
    -A_{2}(\bvx^{m_{i}},  \, 
            \overline{\tau_{-}(\vl_{m-m_{i}-1})}) \right\} 
\nonumber \\ 
& = & B(\bvx^{m_{1}}, \bvx^{m_{2}}, \bvx^{m_{3}}). 
\end{eqnarray*} 
where 
\begin{eqnarray*} 
{\cal Y} &=& \overline{\tau_{-}(\vl_{m-1})} 
    -\gwiitwo{\bvx^{m_{1}} \, \bvx^{m_{2}}  \, \bvx^{m_{3}} \, \gua} \ga \\ 
&& + \sum_{i=1}^{3} \left\{ (m-2m_{i}) \bvx^{m-1} 
    +  \sum_{j=0}^{m-m_{i}-1} \vx^{j} \bullet 
        (\vg * \bvx^{m-1-j}) 
            \right.  \\ 
&& \hspace{50pt} \left. - \sum_{j=0}^{m_{i}-1} \vx^{j} \bullet 
        (\vg * \bvx^{m-1-j}) 
    \right\}. 
\end{eqnarray*} 
By Corollary~\ref{cor:Epower4pt} and Corollary~\ref{cor:Stau-Lk}, 
\[ {\cal Y} = 0. \] 
The theorem is thus proved. 
$\Box$ 
 
\begin{cor} \label{cor:Recpsi} 
For any manifold, 
\begin{eqnarray*} 
&&  2(k-1) \psi_{k+1} - (k+1) T(\bvx) \psi_{k} \\ 
&=& (k+1) A_{2}(\bvx,  \, 
                        \overline{\tau_{-}(\vl_{k-1})}) 
        - (k+1) A_{2}(\bvx^{k}, \, 
                        \overline{\tau_{-}(\vl_{0})}) 
 - (k+1)T(\bvx^{k})A_{1}(\tau_{-}^{2}(\vl_{0})) \\ 
&&  - \delta_{k, 0}  B(\bvx, \bvs, \bvs) 
    + \sum_{j=1}^{k-1} B(\bvx, \bvx^{j}, \bvx^{k-j}) 
\end{eqnarray*} 
for $k \geq 0$. Here the last summation should be 
understood as 0 for $k=0$ and $k=1$. 
\end{cor} 
{\bf Proof}: 
For $k \geq 2$, the formula is obtained by 
first applying Theorem~\ref{thm:psi} for $m_{1}=1$, $m_{2}=j$, $m_{3}=k-j$, 
and summing over $j = 1, \ldots, k-1$, then using 
the genus-2 $L_{-1}$ and $L_{0}$ constraints: 
\[ \psi_{0} = - A_{1}(\tau_{-}^{2}(\vs)) 
 \hspace{20pt} {\rm and} \hspace{20pt} 
 \psi_{1} = A_{1}(\tau_{-}^{2}(\vl_{0})). \] 
For $k=1$, the formula is trivial. For $k=0$, it follows from 
Theorem~\ref{thm:psi} for $m_{1}=1$, $m_{2}=0$, $m_{3}=0$. 
$\Box$ 
 
{\bf Remark}: The right hand side of the formula in 
Corollary~\ref{cor:Recpsi}  only depends on 
genus-0 and genus-1 data. Therefore once $\psi_{2}$ is known, 
 we can compute $\psi_{k}$ recursively from this formula 
for all $k \geq 3$. This is the main reason why we should 
expect the result in Theorem~\ref{thm:L1->L2}.

\subsection{Solving the genus-2 generating function} 
\label{sec:F2} 
 
In this subsection, we will show that if the quantum cohomology is not 
too degenerate, then the recursive relation in Corollary~\ref{cor:Recpsi} 
not only determines all functions $\psi_{k}$, but also determines 
the genus-2 generating function $F_{2}$. In fact, we can give a formula 
for $F_{2}$ in terms of genus-0 and genus-1 data. 
 
Since at each point, the space of primary vectors is finite 
dimensional, in a neighborhood of a generic point, there exists an 
integer $n$ such that $\{ \bvx^{k} \mid k = 0, \cdots, n\}$ are 
linearly independent and 
\begin{equation} \label{eqn:EulerWrap} 
\bvx^{n+1} = \sum_{i=0}^{n} f_{i} \bvx^{i} 
\end{equation} 
where $f_{i}$ are functions in an open subset of the big phase 
space. Multiplying both sides by $\bvx^{k}$, we have 
\begin{equation} \label{eqn:EulerWrapk} 
\bvx^{n+1+k} = \sum_{i=0}^{n} f_{i} \bvx^{i+k} 
\end{equation} 
for any $k \geq 0$. 
\begin{lem} \label{lem:TWf} 
\[ T(\vw) f_{i} = 0 \] 
for any vector field $\vw$ and $i=0, \ldots, n$. 
\end{lem} 
{\bf Proof}: Taking derivative of \eqref{eqn:EulerWrap} along the 
direction $T(\vw)$ and using Lemma~\ref{lem:TWEulerPower}, we have 
\[ - \bvx^{n+1} \bullet \vw = \sum_{i=0}^{n} \left\{ T(\vw) f_{i} 
                \right\} \bvx^{i} 
     - \sum_{i=0}^{n} f_{i} \bvx^{i} \bullet \vw. \] 
The lemma then follows from \eqref{eqn:EulerWrap}. $\Box$ 
 
Define a sequence of vector fields 
\begin{equation} \label{eqn:compctVF} 
\vy_{k} := \sum_{i=0}^{k-1} \bvx^{i} \bullet 
    ((\vg- \frac{1}{2}\vz)* \bvx^{k-1-i}) 
\end{equation} 
where $\vz$ is defined by \eqref{eqn:starID}. 
As explained in \cite{L1} \cite{L2}, in order for the genus-1 Virasoro conjecture 
to hold, the sequence of genus-0 functions representing $\gwiione{\bvx^{k}}$ 
as predicted by the genus-1 Virasoro conjecture must be compatible with 
the linear relation (\ref{eqn:EulerWrap}) 
(this condition was called the algebraic 
compatibility condition in \cite{L2}). Considering this fact and 
 the precise formula for $\gwiione{\bvx^{k}}$ as given in \cite{L1}, it is reasonable 
 to make the following {\it assumption}: 
\begin{equation} \label{eqn:compct} 
\vy_{n+1} = \sum_{i=0}^{n} f_{i} \, \vy_{i}. 
\end{equation} 
This assumption is satisfied 
for all manifolds with semisimple quantum cohomology (cf. 
\cite{L2}). Moreover this assumption implies the algebraic compatibility condition 
for the genus-1 Virasoro conjecture. To see this, we define 
for any vector field $\vw$, 
\[ y_{k}(\vw) = \sum_{i=0}^{k-1} \bvx^{i} \bullet 
    ((\vg- \frac{1}{2}\vz)* (\vw \bullet \bvx^{k-1-i})). \] 
Then $y_{k}(\vw) = - \nabla_{T(\vw)} \vy_{k}$ by Corollary~\ref{cor:TderProd}, 
Lemma~\ref{lem:TWEulerPower} and \eqref{eqn:derGC}. 
Therefore taking covariant derivative with respect to $T(\vw)$ on both sides of 
\eqref{eqn:compct} and using Lemma~\ref{lem:TWf}, we obtain 
\begin{equation} \label{eqn:compW} 
y_{n+1}(\vw) = \sum_{i=0}^{n} f_{i} \, y_{i}(\vw). 
\end{equation} 
Moreover, since 
\[ \vy_{k+1} = y_{k}(\bvx) + \bvx^{k} \bullet 
    ((\vg- \frac{1}{2}\vz)*\bvs), 
\] 
\eqref{eqn:compW} and \eqref{eqn:EulerWrapk} imply that 
\begin{equation} \label{eqn:compctk} 
\vy_{n+1+k} = \sum_{i=0}^{n} f_{i} \, \vy_{i+k} 
\end{equation} 
for all $k \geq 0$. Therefore we also have 
\begin{equation} \label{eqn:compctWk} 
y_{n+1+k}(\vw) = \sum_{i=0}^{n} f_{i} \, y_{i+k}(\vw) 
\end{equation} 
for all $k \geq 0$. 
When restricted to the small phase space, the formula for $\gwiione{\bvx^{k}}$ 
as given in \cite{L1} is a linear combination of 
the trace of the map $\vw \mapsto \vw \bullet \vy_{k}$ and the trace of the 
map $\vw \mapsto y_{k}((\vg - \frac{1}{2}\vz) * \vw)$. Therefore the genus-1 algebraic compatibility 
condition is satisfied. However, we will not use the genus-1 Virasoro 
constraints in this subsection.

Now we come back to the genus-2 Virasoro conjecture. 
By Corollary \ref{cor:Stau-Lk}, 
\begin{equation} \label{eqn:Ytau-L} 
\overline{\tau_{-}(\vl_{k-1})} + \frac{3}{2}k \bvx^{k-1} 
= - \bvx^{k} \bullet \tau_{-}(\vs) - \vy_{k}. 
\end{equation} 
Define 
\[ \widetilde{\psi}_{k} := \psi_{k} - \frac{3k}{2} 
    \gwiitwo{T(\bvx^{k-1})}. \] 
Then by \eqref{eqn:Ytau-L}, 
\[ \widetilde{\psi}_{k} = \gwiitwo{\bvx^{k}} 
    + \gwiitwo{T(\bvx^{k} \bullet \tau_{-}(\vs))} + \gwiitwo{T(\vy_{k})}. \] 
Therefore \eqref{eqn:compctk} implies 
\begin{equation} \label{eqn:psicomp} 
\widetilde{\psi}_{n+1+k} = \sum_{i=0}^{n} f_{i} 
    \widetilde{\psi}_{i+k} 
\end{equation} 
for every $k \geq 0$. Repeatedly applying this equation, we have, 
\begin{equation} \label{eqn:psiwrap} 
\widetilde{\psi}_{n+1+k} = \sum_{i=1}^{n+1} b_{k, i} 
    \widetilde{\psi}_{i} 
\end{equation} 
where $b_{k, i}$ is given by the recursion relation 
\begin{eqnarray*} 
b_{k+1, i} = \left\{ \begin{array}{l} 
     \sum_{j=0}^{k-1} f_{j+n-k+1} b_{j+1, i}  \,\,\,\, {\rm for} \, \, \, \, 
    1 \leq i \leq k,  \\  \\ 
     f_{i-k-1} + \sum_{j=0}^{k-1} 
    f_{j+n-k+1} b_{j+1, i} \,\,\,\, {\rm for} \, \, \, \, 
    k+1 \leq i \leq n+1 \end{array} \right. 
\end{eqnarray*} 
and 
\[ b_{1, i} = f_{i-1} \] 
for $1 \leq i \leq n+1$. 
 
\begin{lem} \label{lem:lineqnpsi} 
For every $k \geq 0$, 
\[ \sum_{i=0}^{n} \left( \frac{n+k}{n+k+2} - \frac{i+k-1}{i+k+1} 
            \right) f_{i} \, \widetilde{\psi}_{i+k+1} = g_{k} \] 
where 
\[ g_{k} :=  \sum_{j=1}^{n+k} \frac{B(\bvx, \bvx^{j}, \bvx^{n+1+k-j})}{2(n+k+2)} 
 - \sum_{i=0}^{n}  \sum_{j=1}^{i+k-1} \frac{f_{i}B(\bvx, \bvx^{j}, \bvx^{i+k-j})}{2(i+k+1)} + \delta_{k, 0} \frac{f_{0}}{2} B(\bvx, \bvs, \bvs). \] 
\end{lem} 
{\bf Proof}: By \eqref{eqn:DerCorr}, 
\[ T(\bvx) \gwiitwo{T(\bvx^{k-1})} 
= \gwiitwo{T(\bvx) T(\bvx^{k-1})} + 
    \gwiitwo{\nabla_{T(\bvx)} T(\bvx^{k-1})}. 
     \] 
By \eqref{eqn:TRR2}, Lemma~\ref{lem:DerT} and 
Lemma~\ref{lem:TWEulerPower}, 
\[ T(\bvx) \gwiitwo{T(\bvx^{k-1})} 
= 2 \gwiitwo{T(\bvx^{k})} + A_{2}(\bvx, \bvx^{k-1}). 
     \] 
Hence by Corollary~\ref{cor:Recpsi}, 
\begin{equation} 
     T(\bvx)\widetilde{\psi}_{k} \, \, = \, \,  \frac{2(k-1)}{k+1} 
            \, \widetilde{\psi}_{k+1} - 3 \gwiitwo{ T(\bvx^{k})} - h_{k} 
\label{eqn:Tpsitilde} 
\end{equation} 
where 
\begin{eqnarray} 
 h_{k} &:=& A_{2}\left(\bvx, \, \, \frac{3}{2} k \bvx^{k-1} + 
    \overline{\tau_{-}(\vl_{k-1})}\right) 
        - A_{2}\left(\bvx^{k}, \, \overline{\tau_{-}(\vl_{0})}\right) 
        - T(\bvx^{k}) \, A_{1}(\tau_{-}^{2}(\vl_{0}))  \nonumber \\ 
 &&      + \frac{1}{k+1} 
    \left\{ -\delta_{k, 0} B(\bvx, \bvs, \bvs) 
    + \sum_{j=1}^{k-1} B(\bvx, \bvx^{j}, \bvx^{k-j}) \right\}. 
\label{eqn:defhk} 
\end{eqnarray} 
 
On the other hand, taking derivative of \eqref{eqn:psicomp} along 
the direction of $T(\bvx)$ and using Lemma~\ref{lem:TWf}, we have 
\[ T(\bvx) \widetilde{\psi}_{n+k+1} = \sum_{i=0}^{n} f_{i} 
        \, T(\bvx) \widetilde{\psi}_{i+k}.\] 
Applying \eqref{eqn:Tpsitilde} to both sides 
of this equation and using \eqref{eqn:psicomp} and 
\eqref{eqn:EulerWrapk}, we obtain 
\[ \sum_{i=0}^{n} \left( \frac{n+k}{n+k+2} - \frac{i+k-1}{i+k+1} 
            \right) f_{i} \, \widetilde{\psi}_{i+k+1} 
= \frac{1}{2} \left(h_{n+k+1} - \sum_{i=0}^{n} f_{i} h_{i+k} \right). \] 
Moreover, equations (\ref{eqn:Ytau-L}), (\ref{eqn:EulerWrapk}) and (\ref{eqn:compctk}) 
 imply that 
\[ \frac{1}{2} \left(h_{n+k+1} - \sum_{i=0}^{n} f_{i} h_{i+k} \right) = g_{k}. \] 
The lemma is thus proved. 
$\Box$ 
 
Applying \eqref{eqn:psiwrap} to replace every $\widetilde{\psi}_{k}$ with $k > n+1$ in 
Lemma~\ref{lem:lineqnpsi} by linear combinations of 
$\widetilde{\psi}_{1}, \ldots, \widetilde{\psi}_{n+1}$,  we obtain the following 
\begin{cor} 
\label{cor:lineqnpsi} For every $k \geq 0$, 
\[ \sum_{i=1}^{n+1} c_{k, i} \, \, \widetilde{\psi}_{i} = g_{k} \] 
where $c_{k, i}$ are given by the recursion relations 
\begin{eqnarray*} 
c_{k, i} = \left\{ \begin{array}{l} 
    - \frac{2}{n+k+2} b_{k+1, i} + \sum_{j=0}^{k-1} 
    \frac{2}{j+n+2} f_{j+n-k+1} b_{j+1, i} \,\,\,\, {\rm for} \, \, \, \, 
    1 \leq i \leq k,  \\  \\ 
    - \frac{2}{n+k+2} b_{k+1, i} + \frac{2}{i} f_{i-k-1} + \sum_{j=0}^{k-1} 
    \frac{2}{j+n+2} f_{j+n-k+1}b_{j+1, i} \,\,\,\, {\rm for} \, \, \, \, 
    k+1 \leq i \leq n+1 \end{array} \right. 
\end{eqnarray*} 
and 
\[ c_{0, i} = \left( \frac{n}{n+2} - \frac{i-2}{i} \right) 
f_{i-1} \] for $1 \leq i \leq n+1$. 
\end{cor} 
 
The following lemma will be proved in the appendix. 
\begin{lem} \label{lem:equiZg1} 
The matrix $(c_{k, i})$ with $0 \leq k \leq n$ and $1 \leq i \leq 
n+1$ is invertible if and only if the matrix 
\[ \left( \begin{array} {llll} 
        \bvx f_{0}, & \bvx^{2} f_{0}, & \ldots, &\bvx^{n+1} f_{0} 
        \\ 
        \bvx f_{1}, & \bvx^{2} f_{1}, & \ldots, &\bvx^{n+1} f_{1} 
        \\ 
        \ldots, & \ldots, & \ldots, & \ldots \\ 
        \bvx f_{n}, & \bvx^{2} f_{n}, & \ldots, &\bvx^{n+1} f_{n} 
        \\ 
        \end{array} \right) 
        \] 
is invertible. 
\end{lem} 
{\bf Remark}: The columns of the last matrix in this lemma 
are given by the coefficients of 
representing the vector fields $Z_{1}, \ldots, Z_{n+1}$ defined in 
\cite{L1} and \cite{L2} in terms of $\bvx^{0}, \ldots, \bvx^{n}$. 
 
\begin{thm} \label{thm:VirNondeg} 
Assume that \eqref{eqn:compct} holds. If the polynomial 
\[ p(x) = x^{n+1} - \sum_{i=0}^{n} f_{i} x^{i} \] 
does not have repeated roots at generic points, then 
the generating function for genus-2 Gromov-Witten invariants is given 
by 
\[ F_{2} = \frac{1}{2} A_{1}(\tau_{-}(\vs)) 
    + \frac{1}{3} A_{1}(\tau_{-}^{2}(\vl_{0})) - \frac{1}{3} \sum_{k=0}^{n} \lambda_{1, k} g_{k} \] 
where $(\lambda_{i, k})$ is the inverse of 
the matrix $(c_{k, i})$ with $0 \leq k \leq n$ and $1 \leq i \leq 
n+1$ and $g_{k}$ is defined in Lemma~\ref{lem:lineqnpsi}. Moreover for $2 \leq i \leq n+1$, 
\[ \psi_{i} = (i-1) \sum_{k=0}^{n} \lambda_{i, k} g_{k} 
    - \frac{i}{2} T(\bvx)\sum_{k=0}^{n} \lambda_{i-1, k} g_{k}  - \frac{i}{2} h_{i-1} \] 
where $h_{i}$ is defined by \eqref{eqn:defhk}. 
\end{thm} 
{\bf Proof}: If the polynomial $p(x)$ does not have repeated 
roots, then by Lemma~\ref{lem:equiZg1} and \cite[Lemma 3.4]{L2}, 
the matrix $(c_{k, i})$ with $0 \leq k \leq n$ and $1 \leq i \leq 
n+1$ is invertible. Let $(\lambda_{i, k})$ be the inverse of this matrix. 
Then by Corollary~\ref{cor:lineqnpsi}, 
\[ \widetilde{\psi}_{i} = \sum_{k=0}^{n} \lambda_{i, k} g_{k} \] 
for $i=1, \ldots, n+1$. 
By definition of $\widetilde{\psi}_{1}$ and the genus-2 $L_{0}$-constraint 
\[ \gwiitwo{T(\bvs)} = \frac{2}{3} \left( \psi_{1} - \widetilde{\psi}_{1} \right) 
    = \frac{2}{3} \left( A_{1}(\tau_{-}^{2}(\vl_{0})) - \sum_{k=0}^{n} \lambda_{1, k} g_{k} \right). \] 
Since 
\[ \vd = T(\vs) = T(\bvs) + T^{2}(\tau_{-}(\vs)), \] 
by the genus-2 dilaton equation and \eqref{eqn:TRR1} 
\[ F_{2} = \frac{1}{2} \gwiitwo{\vd} 
    = \frac{1}{2} \left\{ \gwiitwo{T(\bvs)} + A_{1}(\tau_{-}(\vs)) 
        \right\}.\] 
Therefore we obtain the desired formula for $F_{2}$. 
Moreover, by 
\eqref{eqn:Tpsitilde}, we have 
\[ \gwiitwo{ T(\bvx^{k})} \, =  \, - \frac{1}{3} T(\bvx)\widetilde{\psi}_{k} +  \frac{2(k-1)}{3(k+1)} 
            \, \widetilde{\psi}_{k+1} - \frac{1}{3} h_{k} 
\] 
for $1 \leq k \leq n$. The formula for $\psi_{i}$ is then obtained by using the 
definition of $\widetilde{\psi}_{i}$. The theorem is proved. $\Box$ 
 
{\bf Remark}: (1) 
If the quantum cohomology of the underlying manifold 
is semisimple, the conditions in 
Theorem~\ref{thm:VirNondeg} are satisfied (cf. \cite{L1} \cite{L2}). So in 
this case, we obtained an explicit solution for the generating function of genus-2 
Gromov-Witten invariants. 
(2) The conditions in Theorem~\ref{thm:VirNondeg} 
is a sufficient condition. With careful analysis of the equation 
in Corollary~\ref{cor:lineqnpsi}, one expects that the condition 
can be weakened to similar conditions as in \cite{L1} and \cite{L2}.

\appendix 
\section{Appendix: \hspace{20pt} Proof of Lemma~\ref{lem:equiZg1}} 
 
Define 
\[ B_{k+1} = \left( \begin{array}{c} 
    b_{k+1, 1} \\ b_{k+1, 2} \\ \vdots \\ b_{k+1, n+1} 
        \end{array} \right), \hspace{20pt} 
    C_{k} = \left( \begin{array}{c} 
    c_{k, 1} \\ c_{k, 2} \\ \vdots \\ c_{k, n+1} 
        \end{array} \right), \hspace{20pt} 
    D_{k} = \left( \begin{array}{c} 
    0 \\ \vdots \\ 0 \\ f_{0} \\ \vdots \\ f_{n-k} 
        \end{array} \right) 
        \] 
for $k=0, \ldots, n$. Let 
\[ H ={\rm Diag}(1, 2, \ldots, n+1) 
\] 
be the $(n+1) \times (n+1)$ diagonal matrix whose diagonal entries 
are $1, 2, \ldots, n+1$ and 
\[ H^{-1} ={\rm Diag}(1, \frac{1}{2}, \ldots, \frac{1}{n+1}). 
\] 
 Then $B_{k}$ and 
$C_{k}$ satisfy the recursion relation 
\[ B_{k+1} = D_{k} + \sum_{i=0}^{k-1} f_{n-i} B_{k-i} \] 
and 
\begin{equation} \label{eqn:CDB} 
 C_{k} = 2 H^{-1} D_{k} 
    -\frac{2}{n+k+2}B_{k+1} + 
    \sum_{j=0}^{k-1} \frac{2}{j+n+2} f_{j+n-k+1} B_{j+1}. 
\end{equation} 
 Let 
\[ M = \left( \begin{array}{lllllll} 
            -1 & f_{n} & f_{n-1} & f_{n-2} & \cdots & f_{2} & f_{1} \\ 
            0 & -1 & f_{n} & f_{n-1} & \cdots & f_{3} & f_{2} \\ 
            0 & 0 & -1 & f_{n} &  \cdots & f_{4} & f_{3} \\ 
            \cdots & \cdots & \cdots & \cdots & \cdots & \cdots & \cdots \\ 
            0 & 0 & 0 & 0 & \cdots & -1 & f_{n} \\ 
            0 & 0 & 0 & 0 & \cdots & 0 & -1 
                \end{array} 
    \right). 
\] 
Then the inverse of $M$ has the following form 
\[ M^{-1} = \left( \begin{array}{lllllll} 
            a_{0} & a_{1} & a_{2} & a_{3} & \cdots & a_{n-1} & a_{n} \\ 
            0 & a_{0} & a_{1} & a_{2} & \cdots & a_{n-2} & a_{n-1} \\ 
            0 & 0 & a_{0} & a_{1} &  \cdots & a_{n-3} & a_{n-2} \\ 
            \cdots & \cdots & \cdots & \cdots & \cdots & \cdots & \cdots \\ 
            0 & 0 & 0 & 0 & \cdots & a_{0} & a_{1} \\ 
            0 & 0 & 0 & 0 & \cdots & 0 & a_{0} 
                \end{array} 
    \right) 
\] 
where $a_{0} = -1$ and for $1 \leq k \leq n$, 
\[  a_{k} = \sum_{i=0}^{k-1} a_{i} f_{n-k+i+1}.\] 
Equation~(\ref{eqn:CDB}) implies that 
\begin{equation} 
\frac{1}{n+1+k} B_{k} = \sum_{i=0}^{k-1} a_{i}\left( 
        \frac{1}{2} C_{k-1-i} - H^{-1} D_{k-1-i} \right). 
\label{eqn:BCDa} 
\end{equation} 
 
 Let 
\[ \widetilde{C}_{0} = \frac{1}{2} C_{0}, \] 
\[ \widetilde{C}_{k} = \frac{1}{2} C_{k} - \frac{1}{2}\sum_{i=0}^{k-1} 
        \left( \sum_{j=i}^{k-1} f_{j+n-k+1} a_{j-i} \right) C_{i} 
        \] 
for $1 \leq k \leq n$ and 
\[ V_{0} = (n+2) H \widetilde{C}_{0}, \] 
\[ V_{k} = H \left\{ (n+k+2) \widetilde{C}_{k} - 
            \sum_{i=1}^{k} (n+k-i+2) f_{n-i+1} \widetilde{C}_{k-i} 
            \right\} 
\] for $1 \leq k \leq n$. Then 
\begin{eqnarray} 
&&  \det \left(C_{0}, C_{1},\ldots, C_{n}\right) \neq 0  \nonumber \\ 
& \Leftrightarrow & 
    \det \left(\widetilde{C}_{0}, 
    \widetilde{C}_{1},\ldots, \widetilde{C}_{n}\right) \neq 0 \nonumber \\ 
& \Leftrightarrow & 
    \det \left(V_{0}, V_{1},\ldots, V_{n}\right) \neq 0. 
\label{eqn:CEquiV} 
\end{eqnarray} 
Moreover, using \eqref{eqn:BCDa}, we obtain 
\begin{equation} 
 V_{k} = (n+k+2) D_{k} - H \cdot D_{k} 
    - \sum_{j=0}^{k-1} p_{k-j} D_{j} 
\label{eqn:VD} 
\end{equation} 
where 
\[ p_{k} := \sum_{i=1}^{k} i f_{n-i+1} a_{k-i}. \] 
 
On the other hand, the recursion relation in \cite[Lemma 2.3]{L2} 
can be written as 
\begin{equation} 
 \bvx D_{k} = (n+k+2) D_{k} - H  D_{k}, \label{eqn:XDk} 
\end{equation} 
 and 
\[ \bvx^{k} D_{0} = \left( \bvx^{k-1} f_{n} \right) D_{0} + 
                \bvx^{k-1} D_{1}. \] 
Since $D_{k}$ is obtained from $D_{0}$ by a simple shift, 
recursively applying this formula, we obtain 
\begin{equation} 
\bvx^{k} D_{0} = \sum_{i=0}^{k-2} 
        \left( \bvx^{k-1-i} f_{n} \right) D_{i} + 
                \bvx D_{k-1}. 
\label{eqn:RecXkD} 
\end{equation} 
\begin{lem} \label{lem:pkXkfn} 
\[ \bvx^{k} f_{n} = -p_{k} \] 
for $k \geq 1$. 
\end{lem} 
{\bf Proof}: The lemma is true for $k=1$ since $\bvx f_{n} = f_{n} 
= -p_{1}$. Assume the lemma holds for $1 \leq k \leq m$. By 
\eqref{eqn:RecXkD}, 
\[\bvx^{m+1} f_{n} = \bvx f_{n-m} + \sum_{i=0}^{m-1} f_{n-i} \bvx^{m-i} f_{n}. \] 
By the induction hypothesis and \eqref{eqn:XDk}, 
\begin{eqnarray*} 
\bvx^{m+1} f_{n} & = & (m+1)f_{n-m} - \sum_{i=0}^{m-1} f_{n-i} 
        \sum_{j=1}^{m-i} j f_{n-j+1} a_{m-i-j} \\ 
& = &(m+1)f_{n-m} - \sum_{j=1}^{m} j f_{n-j+1} 
        \sum_{i=0}^{m-j} f_{n-i}  a_{m-i-j}. 
\end{eqnarray*} 
By the recursion relation for $a_{k}$, 
\[ 
\bvx^{m+1} f_{n}  = (m+1)f_{n-m} - \sum_{j=1}^{m} j f_{n-j+1} 
                a_{m-j+1} = - p_{m+1}. 
\] 
So the lemma is proved by induction on $k$. $\Box$ 
 
Lemma~\ref{lem:pkXkfn} and \eqref{eqn:VD}, (\ref{eqn:XDk}), 
(\ref{eqn:RecXkD}) implies that 
\[ V_{k} = \bvx^{k+1} D_{0}\] 
for $0 \leq k \leq n$. Therefore Lemma~\ref{lem:equiZg1} follows 
from \eqref{eqn:CEquiV}.

 

\vspace{30pt} 
\noindent 
Department of Mathematics  \\ 
University of Notre Dame \\ 
Notre Dame,  IN  46556, USA \\ 
 
\vspace{10pt} 
\noindent 
E-mail address: {\it xliu3@nd.edu} 
 
\end{document}